\pdfoutput=1
\documentclass[11pt]{amsart}
\usepackage[margin=1.15in,marginparwidth=0.9in,marginparsep=0.1in]{geometry}
\usepackage[T1]{fontenc}
\usepackage{lmodern}
\usepackage{amssymb,amsmath}
\usepackage{enumitem}
\usepackage{booktabs}
\usepackage{array}
\usepackage{tikz}
\usepackage[hidelinks]{hyperref}
\hypersetup{pdftitle={Toric gravitational instantons with four fixed points},pdfauthor={Lars Andersson, Mattias Dahl, Gustav Nilsson}}

\newtheorem{theorem}{Theorem}[section]
\newtheorem{lemma}[theorem]{Lemma}
\newtheorem{proposition}[theorem]{Proposition}
\newtheorem{corollary}[theorem]{Corollary}

\newtheorem{mainthm}{Theorem}

\theoremstyle{definition}

\theoremstyle{remark}
\newtheorem{remark}[theorem]{Remark}
\numberwithin{equation}{section}

\newcommand{\R}{\mathbb R}
\newcommand{\Q}{\mathbb Q}
\newcommand{\C}{\mathbb C}
\newcommand{\Z}{\mathbb Z}
\newcommand{\F}{\mathbb F}
\newcommand{\CP}{\mathbb{CP}}
\newcommand{\CPb}{\overline{\mathbb{CP}}{}}
\newcommand{\Pl}{\mathcal P}

\title[Toric instantons with four fixed points]{Toric gravitational instantons with four fixed points}

\author[L. Andersson]{Lars Andersson}
\address[Lars Andersson]{Department of Mathematics, KTH Royal Institute of Technology, 100 44 Stockholm, Sweden}
\email{larand@kth.se}

\author[M. Dahl]{Mattias Dahl}
\address[Mattias Dahl]{Department of Mathematics, KTH Royal Institute of Technology, 100 44 Stockholm, Sweden}
\email{dahl@kth.se}

\author[G. Nilsson]{Gustav Nilsson}
\address[Gustav Nilsson]{Mathematisches Institut, Universit\"at M\"unster,\newline Einsteinstrasse 62, 48149 M\"unster, Germany}
\email{gustav.nilsson@uni-muenster.de}

\date{\today}

\begin{document}
\begin{abstract}
We classify the simply connected toric gravitational instantons whose torus action has four fixed points. They are the hyperk\"ahler multi-Eguchi--Hanson and multi-Taub--NUT metrics with four collinear centres, the
asymptotically flat instantons on Li and Sun's manifold $X_3$, up to scale one for each twist at infinity in $(0,1)$ for a fixed
orientation of the torus, and, up to scale, a single further instanton with an ALF end, which is not hyperk\"ahler. The last two are Hermitian for
neither orientation. In particular, the toric asymptotically flat
instanton on $X_3$ with a given twist is unique up to scale, as conjectured by Li and Sun. The proof combines the solitonic form of Geroch's conjecture, which turns smoothness into polynomial equations for the soliton
data, with a convexity principle adapted from Aksteiner. At a fixed asymptotic parameter, the cone errors of the conical metrics with
a given rod structure are the gradient of a strictly convex function of the rod lengths, and the strict convexity gives uniqueness.
Existence follows from an explicit solution for the ALF end and, for the AF ends, from Li and Sun's existence theorem. We also outline how a chart of Aksteiner identifies these AF instantons with those found by Teo.
\end{abstract}
\maketitle
\setcounter{tocdepth}{1}
\tableofcontents
\enlargethispage{\baselineskip}

\section{Introduction}\label{sec:intro}
A \emph{gravitational instanton} is a complete, non-compact Ricci-flat Riemannian $4$-manifold whose curvature has finite energy.
Gravitational instantons were introduced by Hawking in Euclidean quantum gravity \cite{Hawking1977}. The classical examples, among
them the Euclidean Schwarzschild and Kerr metrics, the Eguchi--Hanson and Taub--NUT metrics, and the multi-centre metrics of Gibbons and
Hawking with collinear centres, are \emph{toric}: the torus $T^2$ acts on them effectively by isometries. A natural measure of the
complexity of a toric instanton is the number of fixed points of the torus action. This paper classifies the simply connected toric
gravitational instantons with four fixed points, the smallest number of fixed points for which such an instanton can be Hermitian for
neither orientation. Here a metric is called \emph{Hermitian} for an orientation if it admits a complex structure that is
compatible with the metric and induces this orientation.

The geometry of a toric instanton $(M,g)$ is encoded in its orbit space. If $M$ is simply connected and the action has $n\ge1$ fixed
points, the orbit space $M/T^2$ is a closed half-plane, and the fixed points map to $n$ points of its boundary line, the \emph{axis}. We
call them the \emph{turning points}. They cut the axis into $n-1$ bounded intervals, the finite \emph{rods}, whose lengths
$L_1,\dots,L_{n-1}$ are the \emph{rod lengths}, and two semi-infinite rods. Over each rod the stabiliser of the action is a circle
subgroup of $T^2$. It is generated by a primitive vector of the integer lattice in the Lie algebra of $T^2$, the \emph{rod vector} of the
rod, which is determined up to sign. The turning points, the lines spanned by the rod vectors, and the type of the end, described below,
form the \emph{rod structure} of the instanton \cite{LS2025}. Li and Sun draw it like the moment polygon of a toric symplectic
$4$-manifold, with each rod as an edge orthogonal to its rod vector (Figure~\ref{fig:rods}). Number the rods from left to right, with rod vectors
$v_0,\dots,v_n$. Since the action is smooth, the rod vectors of the two rods that meet at a turning point form a basis of the lattice
\cite{OR1970}. Consequently, for suitable signs, the sum $v_{j-1}+v_{j+1}$ of the neighbours of the rod vector $v_j$ of each finite rod is an
integer multiple $\kappa_jv_j$ of it. The integers $\kappa_1,\dots,\kappa_{n-1}$ form the \emph{chain}, which determines the topology.
The preimage of the closure of each finite rod is an embedded $2$-sphere whose self-intersection number is the corresponding entry of the
chain, and $M$ is the interior of a linear plumbing of disc bundles over these spheres. For example, the
multi-Eguchi--Hanson metrics with four collinear centres have, in their hyperk\"ahler orientation, the chain $(-2,-2,-2)$. The spheres
are the $(-2)$-curves of the minimal resolution of $\C^2/\Z_4$.

At infinity a toric instanton is modelled on one of three types of ends \cite{LS2025}. An \emph{ALE} end is asymptotic to the end of a
quotient $\R^4/\Gamma$ of flat space by a finite group. An \emph{ALF} end is asymptotic, up to a finite cover, to the end of the
Taub--NUT metric, a circle bundle over the complement of a ball in $\R^3$ whose fibres have asymptotically constant length. It is
ALF$^+$ or ALF$^-$ according as the orientation of $M$ corresponds to the hyperk\"ahler orientation of the Taub--NUT metric or to the
opposite one. An \emph{AF} end is asymptotic to the end of the quotient of $\R^3\times\R$ by an isometry that translates $\R$ and rotates
$\R^3$ about an axis through the angle $2\pi\beta$. The number $\beta$, defined modulo $\Z$, is the \emph{twist} of the end, and its sign
depends on an orientation of the torus. Geodesic balls of radius $\mathfrak r$ have volume of order $\mathfrak r^4$ for ALE ends and of
order $\mathfrak r^3$ for ALF and AF ends. The end is AF if and only if the two semi-infinite rods have the same stabiliser.

The classification rests on the soliton representation of toric instantons proved in \cite{Gerochcf}. With two commuting Killing
fields, the Ricci-flat equations reduce to an equation for the Gram matrix of the Killing fields which is the integrability condition
of a linear system with a spectral parameter. The soliton transformations of Belinski and Zakharov \cite{BZ1978,BZ1979} add finitely
many poles in the spectral parameter to a known solution and produce the new metric by linear algebra. Geroch conjectured that an
infinite-dimensional group of such transformations generates all these metrics from flat space \cite{Geroch1971,Geroch1972}. In
\cite{Gerochcf} a solitonic form of this conjecture is proved for toric instantons: each of them is obtained from a flat metric by
finitely many soliton transformations with poles on the axis. Each such instanton is given by a closed formula in finitely many \emph{soliton data}:
a positive weight at each turning point, a weight at infinity, which is zero for ALE ends and positive otherwise, and an integer, the degree index.
Up to a common factor, the components of the rod vectors are Hankel determinants of measures built from the soliton data, hence
polynomials in them. Smoothness, which is the
relation between consecutive rod vectors expressed by the chain, becomes a system of polynomial equations with integer coefficients.
Conversely, every solution with positive weights, the weight at infinity included, and a degree index between $1$ and $n-1$ defines
a smooth complete instanton, whose end is ALF or AF. The classification of toric
instantons with a given number of turning points is thus a question about the positive real solutions of polynomial systems, one for
each chain, type of end, and degree index.

We say that two metrics agree \emph{up to scale} if one is isometric to a constant positive multiple of the other. The main result is the
following.

\begin{mainthm}\label{thm:mainA}
Let $(M,g)$ be a simply connected toric gravitational instanton whose torus action has four fixed points. Then $(M,g)$ is one of the following.
\begin{enumerate}[leftmargin=*,label=\textup{(\alph*)}]
\item A hyperk\"ahler multi-Eguchi--Hanson metric with four collinear centres. Its end is ALE.
\item A hyperk\"ahler multi-Taub--NUT metric with four collinear centres. Its end is ALF.
\item Up to scale, one further instanton $g_{\rm ALF}$, which has an ALF end. It is Hermitian for neither orientation.
\item An instanton on Li and Sun's manifold $X_3=(S^2\times\R^2)\#\CP^2\#\CPb^2$. It is Hermitian for neither
orientation, and its end is AF. In the orientation of the torus in which its chain is $(-1,0,1)$, its twist $\beta$, taken in $[0,1)$,
is not $0$.
For each $\beta\in(0,1)$ there is, up to scale, exactly one simply connected toric gravitational instanton with four fixed points
and an AF end that has, in one orientation of the torus, the chain $(-1,0,1)$ and the twist $\beta$.
\end{enumerate}
Conversely, each of these metrics is a simply connected toric gravitational instanton with four fixed points.
\end{mainthm}

Theorem~\ref{thm:main} in \S\ref{subsec:main} is the detailed version of Theorem~\ref{thm:mainA}: it gives the chains, the soliton data
of $g_{\rm ALF}$, which lie in $\Q(\sqrt2,\sqrt5)$, the types of the ends, and the conventions for the twist.
Parts (a) and (b) are known: the toric hyperk\"ahler instantons are the Gibbons--Hawking metrics with collinear centres, by the
classification of hyperk\"ahler instantons (\S\ref{subsec:literature}) or directly by Li and Sun \cite{LS2025}. The new content concerns the
instantons that are hyperk\"ahler for neither orientation: none of them has an ALE end, only one has an ALF end, and those with an AF
end live on $X_3$, one up to scale for each twist. The instanton $g_{\rm ALF}$ was found in \cite[Rem.~11.5(a)]{Gerochcf} and
independently by Jones \cite{Jones2026}, and by Theorem~\ref{thm:mainA} it is, up to scale, the only simply connected
toric instanton with four fixed points and an ALF end that is hyperk\"ahler for neither orientation. Part (d) proves Li and Sun's
conjecture \cite[Conj.~8.5]{LS2025} for $X_3$: for each twist the rod lengths, and with them the instanton up to scale, are unique. The
AF instantons of~(d) are, as outlined in \S\ref{subsec:AF-teo}, the metrics found by Teo \cite{Teo2026}, which are therefore smooth complete instantons
for every twist in $(0,1)$.

Li and Sun's degree (\S\ref{subsec:rods}) links the rod structure to the Hermitian geometry. For a given orientation, a toric instanton
has degree $0$ if and only if it is hyperk\"ahler, and degree $1$ if and only if it is Hermitian non-K\"ahler, that is, Hermitian with
nonvanishing self-dual Weyl curvature for this orientation \cite[Thm.~7.4]{LS2025}, \cite[\S2]{Li2023}. The degrees for the two
orientations add up to $n-1$ for ALE ends and to $n$ for ALF and AF ends, where $n$ is the number of fixed points
\cite[Rem.~2.8]{Gerochcf}. An instanton that is Hermitian for neither orientation has degree at least $2$ for both, and the first case in
which this is possible is that of four fixed points and an ALF or AF end. Li and Sun proved that such instantons exist. For every $m\ge1$
and every twist $\beta\in(0,1)$ their manifold $X_m$ (\S\ref{subsec:rods}) carries a toric AF instanton with the twist $\beta$
\cite[Thm.~1.1]{LS2025}, which is Hermitian for neither orientation when $m\ge3$, and they conjecture that it is unique for each
twist.

\emph{Overview of this paper.} Section~\ref{sec:setting} gives background and preliminaries and ends with the detailed main
theorem and an outline of its proof (\S\ref{subsec:main}). Section~\ref{sec:reduction} carries out the reduction, and Sections~\ref{sec:ALF}
and~\ref{sec:AF} treat the ALF and the AF chain. Section~\ref{sec:explicit} gives the metrics explicitly, and
Section~\ref{sec:outlook} describes what is known beyond four turning points.
The appendices contain the fold argument (Appendix~\ref{app:ALF}) and the convexity argument (Appendix~\ref{app:convex}). The role of
computation and of an AI assistant is described in a note at the end of the main text.

\section{Background}\label{sec:setting}
A toric instanton has two descriptions that the proof combines: its rod structure, which records the stabilisers of the torus action
along the boundary of the orbit space, and its soliton data, which give the metric in closed form \cite{Gerochcf}.

\subsection{Rod structures, chains, and ends}\label{subsec:rods}
A \emph{gravitational instanton} is a complete, non-compact, oriented Ricci-flat Riemannian $4$-manifold $(M,g)$ whose curvature has
finite energy, $\int_M|\mathrm{Rm}_g|^2\,d\mathrm{vol}_g<\infty$. Following Li and Sun \cite[Def.~2.1]{LS2025} and the companion
paper \cite[Defs.~2.1, 2.2]{Gerochcf}, we call it \emph{toric} if it carries an effective isometric action of the
torus $T^2$ in which every point has a connected stabiliser. By
Lemma~\ref{lem:stab} below, every effective isometric action of $T^2$ on a simply connected gravitational instanton has connected
stabilisers, so the instanton with this action is toric.

We write $\mathfrak t$ for the Lie algebra of $T^2$ and $\Lambda\subset\mathfrak t$ for the lattice of elements $X$ with
$\exp(2\pi X)=1$, so that $T^2=\mathfrak t/\Lambda$. Orlik and Raymond list the possible stabilisers of an effective action of $T^2$ on an orientable $4$-manifold
together with their slice representations \cite[\S1.4]{OR1970}. For closed manifolds $M$ with $H_1(M;\Z)=0$, Orlik and Raymond also
exclude finite non-trivial stabilisers \cite[Lemma~5.2]{OR1970}. The following lemma treats non-compact manifolds.

\begin{lemma}\label{lem:stab}
Let $(M,g)$ be a simply connected, non-compact Riemannian manifold of dimension~$4$ with an effective isometric action of $T^2$. Then the
stabiliser of each point of $M$ is trivial, a circle subgroup of $T^2$, or $T^2$. In particular, it is connected.
\end{lemma}

\begin{proof}
Let $x\in M$, let $H\subset T^2$ be its stabiliser and $O=T^2\cdot x$ its orbit, and let $V\subset T_xM$ be the orthogonal complement of
$T_xO$. For small $\epsilon>0$ and the open ball $V_\epsilon\subset V$ of radius $\epsilon$ about $0$, the map
\begin{equation}\label{eq:tube}
T^2\times_HV_\epsilon\to M,\qquad [t,v]\mapsto t\cdot\exp_x(v),
\end{equation}
is an equivariant diffeomorphism onto an open neighbourhood $U$ of $O$ \cite[Ch.~VI, Thm.~2.2, Cor.~2.4]{Bredon1972}. Here $H$ acts on $V$ by
the differentials of its elements at $x$. Since $T^2$ is abelian, these differentials are the identity on $T_xO$. An isometry of the
connected manifold $M$ that fixes $x$ and whose differential at $x$ is the identity is the identity map. Since the action is effective,
$H$ therefore acts faithfully on $V$. Every element of $H$ is joined to the identity by a path in $T^2$ and therefore preserves the
orientation of $M$, which is orientable since it is simply connected. Its differential at $x$ thus has positive determinant, and $H$ is a
closed subgroup of $SO(V)$ with $\dim V=2+\dim H$.

If $\dim H=2$, then $H=T^2$. If $\dim H=1$, then the identity component $H_0$ of $H$ is the group of rotations of the
$3$-dimensional space $V$ about a line $V_0$. Since $H$ is abelian, each element of $H$ maps $V_0$, the fixed line of $H_0$, onto
itself. A rotation of $V$ with this property is a rotation about $V_0$ or a half-turn about a line orthogonal to $V_0$. A half-turn of
this kind conjugates each rotation about $V_0$ to its inverse, so it does not commute with $H_0$. Hence $H=H_0$ is a circle. If $\dim
H=0$, then $H$ is a finite subgroup of $SO(2)$, hence cyclic of some order $k$, and it remains to show that $k=1$.

Suppose that $k\ge2$. Then $H$ acts on $V\cong\C$ by multiplication with the $k$-th roots of unity. Let $h\in H$ be the
element that acts by $e^{2\pi i/k}$, and choose $\zeta\in\mathfrak t$ with $\exp(2\pi\zeta)=h$. The lattice $\Lambda_x=\Lambda+\Z\zeta$ contains $\Lambda$ with $\Lambda_x/\Lambda\cong
H$, and $X\mapsto\exp(2\pi X)\cdot x$ identifies $\mathfrak t/\Lambda_x$ with $O$. Let $N=M/T^2$ be
the orbit space, $\pi\colon M\to N$ the orbit map, and $y=\pi(x)$. By \eqref{eq:tube}, $D=\pi(U)$ is homeomorphic to $V_\epsilon/H$,
hence to an open disc by $v\mapsto v^k$, with $y$ as its centre. By the slice theorem at all points, $N$ is a surface with boundary
\cite[\S1.4]{OR1970}. It is connected, it is non-compact since $\pi$ is proper, and it is simply connected since
$\pi_*\colon\pi_1(M)\to\pi_1(N)$ is onto \cite[Ch.~II, Cor.~6.3]{Bredon1972}. Gluing $\partial N\times[0,1)$ to $N$ along $\partial N$
gives a connected non-compact surface without boundary which deformation retracts onto $N$, so $H_2(N;\Z)=0$
\cite[Prop.~3.29]{Hatcher2002}. Since also $H_1(N;\Z)=0$ and $H_1(D;\Z)=0$, the Mayer--Vietoris sequence of the open cover of $N$ by $N\setminus\{y\}$ and $D$ shows that
the inclusion induces an isomorphism
\begin{equation}\label{eq:H1N}
\Z\cong H_1(D\setminus\{y\};\Z)\to H_1(N\setminus\{y\};\Z),
\end{equation}
where the first isomorphism is the winding number about the centre in the coordinate $v^k$ of $D$.

The open sets $U$ and $W=M\setminus O$ cover $M$. Both are connected, $W$ because $O$ has codimension $2$, and so is
$U\cap W=U\setminus O$.
Since $h\cdot\exp_x(v)=\exp_x(e^{2\pi i/k}v)$, the universal covering
\begin{equation}\label{eq:tube-cover}
p\colon\mathfrak t\times\R\times(0,\epsilon)\to U\setminus O,\qquad p(X,\theta,s)=\exp(2\pi X)\cdot\exp_x(se^{i\theta}),
\end{equation}
satisfies $p(X+\zeta,\theta,s)=p(X,\theta+2\pi/k,s)$. Its deck transformations are the translations by the lattice
$\Pi\subset\mathfrak t\times\R$ generated by $\Lambda\times\{0\}$, $(0,2\pi)$, and $\nu=(-\zeta,2\pi/k)$, so $\pi_1(U\setminus O)\cong\Pi$.
Since $k\nu+(k\zeta,0)=(0,2\pi)$ and $k\zeta\in\Lambda$, and since the second component of $c\nu$ vanishes only for $c=0$,
$\Pi=(\Lambda\times\{0\})\oplus\Z\nu$.

The set $U$ deformation retracts onto $O\cong\mathfrak t/\Lambda_x$, so $\pi_1(U)\cong\Lambda_x$, and the inclusion of $U\setminus O$
in $U$ induces the map $(X,\theta)\mapsto X$ from $\Pi$ to $\Lambda_x$. Let $\varphi_U\colon\pi_1(U)\to\Z_k$ be the quotient map
$\Lambda_x\to\Lambda_x/\Lambda$ followed by the isomorphism $\Lambda_x/\Lambda\cong\Z_k$ that maps the class of $-\zeta$ to $1$. Let
$\varphi_W\colon\pi_1(W)\to\Z_k$ be the map induced by $\pi\colon W\to N\setminus\{y\}$, followed by the Hurewicz map to
$H_1(N\setminus\{y\};\Z)$, the inverse of \eqref{eq:H1N}, and the reduction modulo $k$. The loop of $U\setminus O$ that corresponds to
$(X,\theta)\in\Pi$ is mapped by $\pi$ to a loop with the winding number $k\theta/2\pi$ about $y$. So both $\varphi_U$ and $\varphi_W$
map $\Lambda\times\{0\}\subset\pi_1(U\setminus O)$ to $0$ and $\nu$ to $1$. By the theorem of Seifert and van Kampen they define a
homomorphism $\pi_1(M)\to\Z_k$, which is onto since $\varphi_U$ is. This contradicts $\pi_1(M)=1$.
\end{proof}

Throughout, $M$
is simply connected and the action has finitely many fixed points, $n\ge1$ of them. The orbit space $M/T^2$ is then a closed
half-plane, with a coordinate $z$ along its boundary line, the \emph{axis}, and the fixed points map to $n$ points $w_1<\dots<w_n$ of
the axis, the \emph{turning points} \cite[Prop.~3.8]{Gerochcf}. They divide the axis into the \emph{rods} $I_0=(-\infty,w_1)$,
$I_j=(w_j,w_{j+1})$ for $1\le j\le n-1$, and $I_n=(w_n,\infty)$. Over each rod the stabiliser of the action is a circle subgroup of
$T^2$.

Fix an orientation of $\mathfrak t$. The stabiliser of the rod $I_j$ is generated by a primitive vector $v_j\in\Lambda$, the \emph{rod vector}, which is
determined up to sign. The \emph{rod structure} of $g$ consists of the turning points, the \emph{rod lines} $\R v_j$ (called rod vectors
in \cite{LS2025}), and the type of the end, defined below \cite[Def.~3.1]{LS2025}, \cite[Def.~2.5]{Gerochcf}. For an instanton the rod lines determine
the rod vectors up to sign. The numbers $L_j=w_{j+1}-w_j$, $1\le j\le n-1$, are the \emph{rod lengths}. At the
fixed point over $w_j$ the vectors $v_{j-1},v_j$ form a $\Z$-basis of $\Lambda$ \cite{OR1970}. With $\det$ the determinant in a positively
oriented $\Z$-basis of $\Lambda$, the signs of the rod vectors can therefore be chosen, uniquely up to one common sign, such that
\begin{equation}\label{eq:det}
\det(v_{j-1},v_j)=-1,\qquad 1\le j\le n .
\end{equation}
Since $(v_{j-1},v_j)$ and $(v_j,v_{j+1})$ are then bases of $\Lambda$ with the same determinant, $v_{j+1}+v_{j-1}$ is an integer multiple of
$v_j$:
\begin{equation}\label{eq:chain}
v_{j+1}+v_{j-1}=\kappa_j\,v_j,\qquad \kappa_j=-\det(v_{j-1},v_{j+1})\in\Z,\qquad 1\le j\le n-1 .
\end{equation}
We call $\kappa=(\kappa_1,\dots,\kappa_{n-1})$ the \emph{chain} of $g$.

The orientations of $M$ and of $\mathfrak t$ are related as follows. By \cite[Prop.~3.8]{Gerochcf}, $g$ has global Weyl coordinates $(\rho,z)$, in which the
orbit space is the closed half-plane $\{\rho\ge0\}$. They are computed in a positively oriented $\Z$-basis of $\Lambda$, with dual angles
$\phi^1,\phi^2$. Here $\rho=(\det G)^{1/2}$ for the Gram matrix $G$ of the corresponding Killing fields,
and $z$ is the harmonic conjugate of $\rho$ for which $d\rho\wedge dz\wedge d\phi^1\wedge d\phi^2$ is positive
\cite[\S2.1]{Gerochcf}. This fixes the orientation of $M$ relative to that of $\mathfrak t$.

The chain is the combinatorial part of the rod structure and the lattice. It is determined by the rod vectors through \eqref{eq:chain}, and conversely
it determines them up to an automorphism of $\Lambda$: every $\Z$-basis $(v_0,v_1)$ with $\det(v_0,v_1)=-1$ extends by
$v_{j+1}=\kappa_jv_j-v_{j-1}$ to rod vectors with the chain $\kappa$. The rod vectors in turn determine $M$ with its torus action up to
equivariant diffeomorphism \cite{OR1970}, \cite[proof of Thm.~1]{Nilsson2024}. The rod lengths and the soliton data of \S\ref{subsec:hankel} determine
the metric. The chain also has a direct geometric meaning. The preimage of the closure of a finite rod $I_j$ is an embedded $2$-sphere
$S_j\subset M$. The spheres $S_j$ and $S_{j+1}$ meet transversally in the fixed point over $w_{j+1}$, spheres that are not consecutive
are disjoint, and the classes of $S_1,\dots,S_{n-1}$ form a basis of $H_2(M;\Z)$ \cite[Thm.~1]{Nilsson2024}. The
\emph{self-intersection number} $S_j\cdot S_j$ is the Euler number of the normal bundle of $S_j$, that is, the intersection number of
$S_j$ with a small generic deformation of itself. For the orientation of $M$ just fixed
\begin{equation}\label{eq:selfint}
S_j\cdot S_j=\kappa_j,\qquad 1\le j\le n-1
\end{equation}
(see below). So $M$ is the interior of a linear plumbing $\Pl(\kappa_1,\dots,\kappa_{n-1})$ of disc bundles over $S^2$ with the
Euler numbers $\kappa_1,\dots,\kappa_{n-1}$. Its intersection form in the basis $(S_j)$ is the tridiagonal matrix
$\mathrm{tridiag}(\kappa_j;1)$. For example, the multi-Eguchi--Hanson metrics with four collinear centres have, in their hyperk\"ahler
orientation, the chain $(-2,-2,-2)$, the spheres $S_j$ being the $(-2)$-curves of the minimal resolution of $\C^2/\Z_4$. With the
opposite orientation the chain is $(2,2,2)$. The chain depends on the orientations: reversing the orientation of $\mathfrak t$ replaces
$\kappa$ by $(\kappa_{n-1},\dots,\kappa_1)$, and reversing the orientation of $M$ replaces it by $-(\kappa_{n-1},\dots,\kappa_1)$
(see below). Statements about chains refer to their classes under these two operations.

To see \eqref{eq:selfint}, we use the notation of Li and Sun. They orient $M$ in the
same way and choose rod vectors $u_j$ with $\det(u_{j-1},u_j)=1$. Then the self-intersection number of $S_j$ is
$-\det(u_{j-1},u_{j+1})$ \cite[\S A.7]{LS2026b}, and the substitution $u_j=(-1)^jv_j$ turns this into $\kappa_j$. Remark~\ref{rem:hk} confirms that the hyperk\"ahler orientation is the one in which the spheres are $(-2)$-curves and the chain is
$(-2,\dots,-2)$.

The two orientations enter as follows. If the orientation of $\mathfrak t$ is reversed and that of $M$ kept, then $z$ changes sign, the
rods are numbered in the opposite order, and determinants change sign. So $v'_j=v_{n-j}$ satisfies \eqref{eq:det}, and $\kappa$ is replaced
by $(\kappa_{n-1},\dots,\kappa_1)$. If the orientation of $M$ is reversed and that of $\mathfrak t$ kept, then $z$ changes sign and
determinants do not. So $v'_j=(-1)^jv_{n-j}$ satisfies \eqref{eq:det}, and $\kappa$ is replaced by $-(\kappa_{n-1},\dots,\kappa_1)$. For
$n=4$, \eqref{eq:chain} gives
\begin{equation}\label{eq:v4}
v_4=(1-\kappa_2\kappa_3)\,v_0-p\,v_1,\qquad p=\det(v_0,v_4)=\kappa_1+\kappa_3-\kappa_1\kappa_2\kappa_3 .
\end{equation}

The type of the end is read off from the outer rods and the asymptotics of the metric \cite[Def.~4.30]{LS2025},
\cite[Def.~2.9]{Gerochcf}. The end is \emph{ALE} if it is asymptotic to the end of $\R^4/\Gamma$ for a finite group $\Gamma$. It is
\emph{ALF}$^+$, respectively \emph{ALF}$^-$, if a finite cover of it is asymptotic to the end of the Taub--NUT metric with its
hyperk\"ahler orientation, respectively with the opposite orientation. The end of the Taub--NUT metric is a circle bundle over the complement of a ball in $\R^3$ whose
fibres have asymptotically constant length. The end is \emph{AF}$_\beta$ if it is asymptotic, through an equivariant isometry, to the end of
$(\R^3\times\R)/\Z$ with its torus action, where the generator rotates $\R^3$ by the angle $2\pi\beta$ about an axis and translates $\R$.
The number $\beta$, defined modulo $\Z$ and represented in $[0,1)$, is the \emph{twist} of the end \cite[\S1, Def.~4.30]{LS2025}. The
sense of the rotation is fixed as in Li and Sun's model \cite[(3.8)]{LS2025}, whose basis $(\partial_{\phi_1},\beta\partial_{\phi_1}+\partial_{\phi_2})$
is positively oriented. Concretely, if $\xi\in\mathfrak t$ is the Killing field of the model that translates $\R$, scaled to the period of the
model, and the rod vector $v_0$ of the outer rods is taken with the sign $\det(v_0,\xi)=1$, then $\Lambda=\Z v_0+\Z(\beta v_0+\xi)$. The
twist depends on the orientation of
$\mathfrak t$, not on that of $M$. Reversing the orientation of $\mathfrak t$ replaces $\beta$ by $-\beta$. The reflection of $\R^3$ in
the plane through the centre orthogonal to the axis of rotation reverses the orientation of the model and commutes with the generator
and with the torus action. In \cite[Rem.~4.31]{LS2025}, where the orientation of $M$ determines that of $\mathfrak t$, reversing the
orientation replaces AF$_\beta$ by AF$_{-\beta}$ up to an automorphism of $T^2$. This automorphism reverses the orientation of
$\mathfrak t$.
Geodesic balls of radius $\mathfrak r$ have volume of
order $\mathfrak r^4$ for ALE ends and of order $\mathfrak r^3$ for ALF and AF ends, as for the model ends. The end is AF if and only if
the two outer rods have the same stabiliser, that is, when $p=\det(v_0,v_n)$ vanishes \cite[Def.~3.1]{LS2025}. ALE and ALF ends have
$p\ne0$.

The AF instantons of the main theorem live on manifolds introduced by Li and Sun,
\begin{equation}\label{eq:Xn}
\begin{aligned}
X_m&=(S^2\times\R^2)\#k\,\CP^2\#k\,\CPb^2\qquad &&(m=2k+1),\\
X_m&=(S^2\times\R^2)\#k\,\CP^2\#(k+1)\,\CPb^2\qquad &&(m=2k+2),
\end{aligned}
\end{equation}
which carry torus actions with $m+1$ fixed points \cite[(1.2), \S6.2]{LS2025}. In particular $X_3=(S^2\times\R^2)\#\CP^2\#\CPb^2$ has
four, and its chain is in the class of $(1,0,-1)$. In Li and Sun's orientation of $\mathfrak t$, that of their model, it is
$(-1,0,1)$ \cite[\S6.2]{LS2025}. Aksteiner's twist \cite[\S5.4]{Aksteiner2026} is minus the twist defined above: for the rod
structures of the manifolds $X_m$ the twists he allows form the interval $(-1,0)$, and Li and Sun's form $(0,1)$
\cite[Thm.~1.1]{LS2025}.

The \emph{Li--Sun degree} of a rod structure with the rod lines $[v_0],\dots,[v_n]$ \cite[(3.2), Def.~3.11]{LS2025} is defined as
follows. For $1\le j\le n$ rotate $[v_{j-1}]$ counter-clockwise, that is, in the sense of the orientation of $\mathfrak t$, through the
angle in $(0,\pi)$ that takes it to $[v_j]$. If $[v_n]\ne[v_0]$, rotate $[v_n]$ in the same sense through the angle in $(0,\pi)$ that
takes it to $[v_0]$. The sum of these angles, measured in a positively oriented basis of $\mathfrak t$, is $\mathsf e\pi$. The integer
$\mathsf e\ge0$ is the winding number of the resulting loop in $\mathbb{RP}^1$ and does not depend on the basis. The degree is $d=\mathsf e-1$ for ALE and
ALF$^+$ ends and $d=\mathsf e$ for ALF$^-$ and AF ends. If $\bar d$ is the degree of the rod structure of $M$ with the opposite orientation, then
\begin{equation}\label{eq:degid}
d+\bar d=n-1\quad\text{(ALE)},\qquad d+\bar d=n\quad\text{(ALF, AF)}
\end{equation}
\cite[Rem.~2.8]{Gerochcf}.

\subsection{Results from the literature}\label{subsec:literature}
Two results of Li and Sun relate the degree to the Hermitian geometry and bound the chain (Theorems~\ref{thm:LSdegree}
and~\ref{thm:LSbound}). The first of them involves the Gibbons--Hawking metrics, which we recall. Let
$\mathbf y_1,\dots,\mathbf y_n$ be distinct points on a line in $\R^3$, $U_\infty\ge0$, $U=U_\infty+\sum_{k=1}^n\frac1{2|\mathbf y-\mathbf y_k|}$, and $\eta_{\rm GH}$ a $1$-form with
$d\eta_{\rm GH}=*dU$. The \emph{Gibbons--Hawking metric}
\begin{equation}\label{eq:GH}
g=U^{-1}(d\psi+\eta_{\rm GH})^2+U\,d\mathbf y\cdot d\mathbf y,\qquad \psi\in\R/2\pi\Z,
\end{equation}
with the collinear centres $\mathbf y_1,\dots,\mathbf y_n$ is hyperk\"ahler. It is the multi-Eguchi--Hanson metric for $U_\infty=0$ and the
multi-Taub--NUT metric for $U_\infty>0$ \cite{GH1978}. It is a toric gravitational instanton with $n$ turning points, the centres
\cite[\S11.2]{Gerochcf}.

\begin{theorem}[{\cite[Thm.~7.4]{LS2025}}]\label{thm:LSdegree}
Let $(M,g)$ be a toric gravitational instanton, with $M$ oriented as in \S\ref{subsec:rods}, and let $d$ be the degree of its rod
structure. Then $d=0$ if and only if $g$ is hyperk\"ahler for this orientation, and then $g$ is a Gibbons--Hawking metric with
collinear centres. Further, $d=1$ if and only if $g$ is Hermitian non-K\"ahler for this orientation, and $d\ge2$ if and only if $g$ is not
Hermitian for this orientation.
\end{theorem}

\begin{theorem}[{\cite[Lemma~A.20]{LS2026b}}]\label{thm:LSbound}
Let $(M,g)$ be as in \S\ref{subsec:rods}, and suppose that $g$ is hyperk\"ahler for neither orientation. Then $\kappa_j\in\{-1,0,1\}$
for $1\le j\le n-1$.
\end{theorem}
Both theorems are stated for proper toric gravitational instantons \cite[Def.~3.14]{LS2025}, \cite[Defs.~A.1, A.2]{LS2026b}. Besides
properness, \cite[Def.~A.1]{LS2026b} requires the metric to be integrable and regular and to have quadratic curvature decay. The instantons
of \S\ref{subsec:rods} have these four properties. The curvature decays quadratically, since for a complete Ricci-flat
$4$-manifold finite energy is equivalent to quadratic curvature decay \cite[Rem.~2.3]{LS2026}. The distribution orthogonal to the orbits is integrable since the action has fixed points \cite[\S2.1]{Gerochcf}. Finally, $(\rho,z)$ is a
homeomorphism of the orbit space onto $\{\rho\ge0\}$ and a conformal diffeomorphism of its interior onto the open half-plane
\cite[Prop.~3.8(ii)]{Gerochcf}. So $d\rho\ne0$ in the interior, which is regularity, and the instanton is proper
\cite[Def.~A.2]{LS2026b} (see also the proof of \cite[Prop.~3.8(iii)]{Gerochcf}). Theorem~7.4 of \cite{LS2025} is stated for the harmonic map
strongly tamed by a rod structure \cite[Defs.~4.9, 4.10, Thm.~4.24]{LS2025}. By \cite[Prop.~3.7]{Gerochcf} the harmonic map of $g$ is
locally tamed by the rod structure of $g$ and tends to the model of its end. By the uniqueness lemma \cite[Lemma~3.1]{Gerochcf} it is
therefore the harmonic map strongly tamed by this rod structure. The hyperk\"ahler gravitational instantons are classified \cite[\S4]{LS2026}. That a toric hyperk\"ahler
instanton is a Gibbons--Hawking metric with collinear centres \cite{GH1978} is the toric case of this classification, which Li and Sun
also prove directly \cite[\S7.1.1 and the proof of Thm.~7.7]{LS2025}. Below we need only that degree $0$ characterises the hyperk\"ahler metrics. The
last statement of Theorem~\ref{thm:LSdegree} follows from the first two. Indeed, a metric that is Hermitian for an orientation and not
Hermitian non-K\"ahler for it has vanishing self-dual Weyl curvature. Since $g$ is Ricci-flat and $M$ is simply connected, it is then
hyperk\"ahler for this orientation.

\subsection{Soliton data and Hankel determinants}\label{subsec:hankel}
The soliton representation of \cite{Gerochcf} describes a toric instanton by finitely many soliton data, and its rod vectors by Hankel
determinants of measures built from them.

Geroch conjectured that the local solutions of the vacuum Einstein equations with two commuting Killing fields and vanishing
twist constants are generated from flat space by an infinite-dimensional group of transformations \cite{Geroch1971,Geroch1972}.
Hauser and Ernst proved this for stationary axisymmetric spacetimes near the axis \cite{HE1981}. The solutions generated from flat
space by the transformations with finitely many poles in the spectral parameter are the gravitational solitons of Belinski and
Zakharov \cite{BZ1978,BZ1979,BV2001} (Section~\ref{sec:intro}). In Riemannian
signature a solitonic form of the conjecture holds for toric instantons \cite[Thm.~1]{Gerochcf}. The end of $M$ is of one of the types ALE, ALF$^\pm$, or
AF$_\beta$ of \S\ref{subsec:rods}, and $g$ is a gravitational soliton. It is obtained from a flat metric, on $\R^4$ for an ALE end and
on $\R^3\times S^1$ otherwise, by finitely many soliton transformations with poles on the axis. The proof gives the metric in closed
form \cite[Thm.~10.1, (10.1)]{Gerochcf}. The closed formula depends on the \emph{soliton data}: a positive measure
\begin{equation}\label{eq:omega}
\omega=\sum_{k=1}^na_k\,\delta_{w_k},\qquad a_k>0,
\end{equation}
on the turning points, a number $\gamma\ge0$, the weight of an atom at infinity, which vanishes if and only if the end is ALE, and an
integer $r$, the \emph{degree index} \cite[\S3.3]{Gerochcf}. The Gram matrix of the Killing fields and the conformal factor of the
metric on the orbit space are Hankel determinants of measures built from $\omega$ and from the distances to the turning points
(\S\ref{subsec:geroch}). On the axis the formula becomes algebraic. For a measure $\tilde\omega=\sum_k\tilde\omega_k\delta_{w_k}$ with
real weights let $H_s(\tilde\omega)=\det(h_{i+l})_{i,l=0}^{s-1}$ be the Hankel determinant of its moments $h_q=\sum_k\tilde\omega_kw_k^q$.
Let $H^{(\gamma)}_{r+1}(\tilde\omega)$ be $H_{r+1}(\tilde\omega)$ with $\gamma$ added to the moment of order $2r-1$
\cite[(5.22)]{Gerochcf}. Put
\begin{equation}\label{eq:V}
V_j=\bigl(H_r(\omega_{(j)}),\,H^{(\gamma)}_{r+1}(\omega_{(j)})\bigr)\in\R^2,\qquad
\omega_{(j)}=-\sum_{k\le j}a_k\,\delta_{w_k}+\sum_{k>j}a_k\,\delta_{w_k},\qquad 0\le j\le n,
\end{equation}
so that $\omega_{(j)}$ is $\omega$ with the signs of the atoms below the rod $I_j$ reversed. Then, in a suitable basis of the Lie
algebra of $T^2$, the generator of the stabiliser of $I_j$ is a fixed multiple of $V_j$, up to sign (Theorem~\ref{thm:soliton}).

Heine's formula expresses these Hankel determinants as sums over sets of atoms. For a measure
$\tilde\omega=\sum_{k=1}^n\tilde\omega_k\delta_{w_k}$ with real weights on $w_1<\dots<w_n$ and a set $S\subset\{1,\dots,n\}$ put
\begin{equation}\label{eq:PS}
P_S(\tilde\omega)=\prod_{k\in S}\tilde\omega_k\prod_{k<l\in S}(w_k-w_l)^2,
\end{equation}
with empty products equal to $1$. We call them the \emph{Heine terms}. By Heine's formula \cite[Lemma~5.1]{Gerochcf},
\begin{equation}\label{eq:heine}
H_s(\tilde\omega)=\sum_{|S|=s}P_S(\tilde\omega)\qquad(s\ge1),\qquad H_0=1,\qquad H_{-1}=0 .
\end{equation}
For $r=2$ the shifted determinant is \cite[(5.20), (5.22)]{Gerochcf}
\begin{equation}\label{eq:H3gamma}
H^{(\gamma)}_3=H_3-2\gamma\,\Sigma_2-\gamma^2H_1,\qquad \Sigma_s(\tilde\omega)=\sum_{|S|=s}P_S(\tilde\omega)\sum_{k\in S}w_k ,
\end{equation}
so that $\Sigma_2(\tilde\omega)=\sum_{k<l}P_{\{k,l\}}(\tilde\omega)\,(w_k+w_l)$. The vectors \eqref{eq:V} form the \emph{Hankel family} of $(\omega,r,\gamma)$ \cite[Def.~5.7]{Gerochcf}. For $\gamma=0$ the second component is
$H_{r+1}(\omega_{(j)})$.

\subsection{The soliton representation and its converse}\label{subsec:geroch}
The metric of a toric instanton is given in closed form by its soliton data, and conversely positive data that satisfy the smoothness
relations define instantons (Theorems~\ref{thm:soliton} and~\ref{thm:converse}). A normalisation of the data changes the metric only by a
constant factor (Lemma~\ref{lem:normalise}). The closed formula \cite[(10.1)]{Gerochcf} uses, for a point $(\rho,z)$ off the
axis, the distances $R_k=(\rho^2+(z-w_k)^2)^{1/2}$ to the turning points and the unimodular numbers $\varpi_k=(w_k-z-i\rho)/R_k$. It also uses the measures
$\varpi\omega$, $\mathbf R\omega$, and $\omega/\mathbf R$ on the turning points with the weights $\varpi_ka_k$, $R_ka_k$, and $a_k/R_k$. Hankel
determinants of complex measures are defined by the same determinants. In a basis $(f_0,f_1)$ of $\mathfrak t$ with dual angles
$\psi^0,\psi^1$, put
\begin{equation}\label{eq:closed}
\begin{gathered}
g=\frac{\bigl|\mathcal Y\,d\psi^0-\mathcal X\,d\psi^1\bigr|^2}{F_1F_2}+\mathfrak c\,F_1F_2\,(d\rho^2+dz^2),\\
\begin{aligned}
\mathcal X&=H_r(\varpi\omega),&\qquad \mathcal Y&=H^{(\gamma)}_{r+1}(\varpi\omega)-i\gamma\rho\,H_r(\varpi\omega),\\
F_1&=H_r(\mathbf R\omega)+\gamma H_{r-1}(\mathbf R\omega),&\qquad F_2&=H_{r+1}(\omega/\mathbf R)+\gamma H_r(\omega/\mathbf R),
\end{aligned}
\end{gathered}
\end{equation}
with a constant $\mathfrak c>0$.

\begin{theorem}[soliton representation; {\cite[Thm.~10.1, Rem.~10.3]{Gerochcf}}]\label{thm:soliton}
Let $(M,g)$ be as in \S\ref{subsec:rods}. There are an integer $r$, the degree index, a positive measure $\omega=\sum_ka_k\delta_{w_k}$ on the
turning points, a number $\gamma$, with $\gamma=0$ if the end is ALE and $\gamma>0$ otherwise, a constant $\mathfrak c>0$, and a positively oriented
basis $(f_0,f_1)$ of $\mathfrak t$ with the following properties.
\begin{enumerate}[leftmargin=*,label=\textup{(\roman*)}]
\item The Li--Sun degree of the rod structure of $g$ is $n-1-r$ if the end is ALE and $n-r$ otherwise, and that of $g$ with the opposite
orientation is $r$. In particular $0\le r\le n-1$, respectively $0\le r\le n$.
\item In the Weyl coordinates of $(f_0,f_1)$ the metric is \eqref{eq:closed}.
\item In components with respect to $(f_0,f_1)$, and after a change of the common sign in \eqref{eq:det}, the rod vectors are
$v_j=\mathfrak c^{1/2}V_j$, $0\le j\le n$, with $V_j$ the Hankel family \eqref{eq:V} of $(\omega,r,\gamma)$.
\end{enumerate}
\end{theorem}
Part (i) is \cite[Lemma~3.11]{Gerochcf}, whose statement gives the two degrees as an unordered pair. Their order in (i) is the one in
the proof of \cite[Lemma~3.11]{Gerochcf}. Part (iii) combines \cite[Prop.~7.5(a)]{Gerochcf}, by which the \emph{scaled rod vectors}
\cite[Def.~4.17]{LS2025}, which for a smooth metric are the rod vectors up to sign, are $\pm \mathfrak c^{1/2}V_j$, with
$\det(V_{j-1},V_j)<0$ \cite[Prop.~7.5(b)]{Gerochcf}. Since $(f_0,f_1)$ is positively oriented, \eqref{eq:det} then allows only one sign
for all $j$. The basis of
\cite[Thm.~10.1]{Gerochcf} has determinant $1$ relative to the basis of the model of the end. If it is negatively oriented for the
orientation of $\mathfrak t$ fixed in \S\ref{subsec:rods}, replace $(\omega,f_1,z)$ by $(\omega^\flat,-f_1,-z)$, where
$\omega^\flat=\sum_ka_k\delta_{-w_k}$ is the reflected measure. Then $\varpi_k$ becomes $-\bar\varpi_k$, so the moments $h_q$ of $\varpi\omega$ become $(-1)^{q+1}\bar h_q$.
As a result, $\mathcal X$ becomes $(-1)^r\bar{\mathcal X}$ and $\mathcal Y$ becomes $(-1)^{r+1}\bar{\mathcal Y}$, $F_1$ and $F_2$ are unchanged, and
\eqref{eq:closed} is unchanged.

\begin{theorem}[converse; {\cite[Prop.~11.4]{Gerochcf}}]\label{thm:converse}
Let $\omega$ be a positive measure on $n$ points, $1\le r\le n-1$, and $\gamma>0$, and suppose that the Hankel family \eqref{eq:V} satisfies
$V_{j+1}+V_{j-1}=\kappa_jV_j$ for $1\le j\le n-1$ with integers $\kappa_j$. Put $\mathfrak c=|D_0|^{-1}$ with
$D_0=\det(V_0,V_1)$. Then \eqref{eq:closed}, on the torus bundle whose lattice of periods is generated by $\mathfrak c^{1/2}V_0$ and $\mathfrak c^{1/2}V_1$, is a smooth complete simply
connected toric gravitational instanton with the turning points of $\omega$, the rod vectors $v_j=\mathfrak c^{1/2}V_j$, and the chain $\kappa$. Its
end is ALF$^-$, ALF$^+$, or AF according as
\begin{equation}\label{eq:alpha}
\alpha(\omega,\gamma)=\frac{2\bigl(H_{r+1}(\omega)-\gamma^2H_{r-1}(\omega)\bigr)}{\gamma\,H_r(\omega)}
\end{equation}
is positive, negative, or zero, with $\alpha$ as in \cite[(8.36)]{Gerochcf}.
\end{theorem}
We call $\alpha(\omega,\gamma)$ the \emph{asymptotic parameter} of the data. For a rod structure whose outer rod lines are the lines of
$(1,\frac\alpha2)$ and $(1,-\frac\alpha2)$, with $\alpha\ge0$, the number $\alpha$ is the parameter of the model of its end
\cite[(4.7), Defs.~4.4, 4.6]{LS2025}. We call it the asymptotic parameter of the rod structure. For the rod structures of positive
data, written as in Appendix~\ref{app:convex}, the two agree by \eqref{eq:appC-nouter}. Here $D_0<0$ by \cite[Prop.~5.8(i)]{Gerochcf}. In
\cite[Prop.~11.4]{Gerochcf} the relations appear as the hypothesis that in the recursion \cite[(2.12)]{Gerochcf} for the vectors $V_j$ the
coefficient of $V_{j-1}$ is $-1$ and that of $V_j$ an integer $\kappa_j$. In particular, all consecutive determinants $\det(V_{j-1},V_j)$
are equal. The metric there is \eqref{eq:closed} after a change of basis of determinant $1$ \cite[Thm.~7.3(b), (7.1),
Prop.~8.9]{Gerochcf}, with some constant $\mathfrak c'>0$ in place of $\mathfrak c$. Its scaled rod vectors are
$\pm\mathfrak c'^{1/2}V_j$ \cite[Prop.~8.7(i)]{Gerochcf}. For $\lambda>0$, \eqref{eq:closed} with the constant $\lambda^2\mathfrak c$, in
the angles $\lambda\psi$, is $\lambda^2$ times \eqref{eq:closed} with the constant $\mathfrak c$. The lattice generated by
$\lambda\mathfrak c^{1/2}V_0$, $\lambda\mathfrak c^{1/2}V_1$ corresponds to that generated by $\mathfrak c^{1/2}V_0$, $\mathfrak c^{1/2}V_1$.
This gives Theorem~\ref{thm:converse} with $\mathfrak c=|D_0|^{-1}$. Simple connectivity holds because
the orbit space is simply connected and consecutive rod vectors generate the lattice \cite[Rem.~3.9]{Gerochcf}, \cite{OR1970}. Off the axis the metric \eqref{eq:closed} is
positive definite for all positive data: $F_1,F_2>0$ by Heine's formula, and $\operatorname{Im}(\mathcal X\overline{\mathcal Y})=\rho F_1F_2$
\cite[Lemma~7.2]{Gerochcf}.

The following lemma shows that the normalisation of the data used below changes the metric only by a constant factor.

\begin{lemma}\label{lem:normalise}
Let $r=2$, let $\omega=\sum_ka_k\delta_{w_k}$ be a positive measure and $\gamma>0$, let $\lambda,K>0$ and $\mu\in\R$, and put
$\omega'=\sum_kKa_k\delta_{\lambda w_k+\mu}$ and $\gamma'=K\lambda^3\gamma$. Then the Hankel families \eqref{eq:V} satisfy $V'_j=BV_j$ for
$0\le j\le n$, with
\begin{equation}\label{eq:normB}
B=\begin{pmatrix}K^2\lambda^2&0\\-4\mu\gamma'K^2\lambda^2&K^3\lambda^6\end{pmatrix},\qquad \det B=K^5\lambda^8 .
\end{equation}
If the Hankel family of $(\omega,\gamma)$ satisfies the hypotheses of Theorem~\ref{thm:converse}, so does that of $(\omega',\gamma')$. In that case
the metric \eqref{eq:closed} of $(\omega',\gamma')$, in the coordinates $\rho'=\lambda\rho$, $z'=\lambda z+\mu$ and with angles related to
those of $(\omega,\gamma)$ by a constant linear map, is $\lambda$ times the metric \eqref{eq:closed} of $(\omega,\gamma)$.
\end{lemma}

\begin{proof}
At corresponding points the distances become $\lambda R_k$, the numbers $\varpi_k$ do not change, the factors $w_k-w_l$ are multiplied by
$\lambda$, and $w_k+w_l$ becomes $\lambda(w_k+w_l)+2\mu$. By \eqref{eq:PS} and \eqref{eq:heine}, $H_s$ of a measure built from $\omega'$ is
$K^s\lambda^{s(s-1)}$ times $H_s$ of the corresponding measure built from $\omega$, multiplied by $\lambda^{\pm s}$ for the weights $R_ka_k$ and
$a_k/R_k$. The sum $\Sigma_2$ becomes $K^2\lambda^3\Sigma_2+2\mu K^2\lambda^2H_2$. With \eqref{eq:H3gamma} this gives $V'_j=BV_j$, the same relation for
$(\mathcal X,\mathcal Y)$ at corresponding points, $F'_1=K^2\lambda^4F_1$, and $F'_2=K^3\lambda^3F_2$. In particular $D'_0=\det B\,D_0$, so
$\mathfrak c'=\mathfrak c/\det B$, and the relations $V'_{j+1}+V'_{j-1}=\kappa_jV'_j$ hold with the same $\kappa_j$. With $J=\bigl(\begin{smallmatrix}0&1\\-1&0
\end{smallmatrix}\bigr)$, the angles defined by the first relation of
\begin{equation}\label{eq:norm-angles}
d\psi=(\det B)^{-1/2}(JBJ^{-1})^Td\psi',\qquad
\mathcal Y'd\psi'^0-\mathcal X'd\psi'^1=(\det B)^{1/2}\bigl(\mathcal Y\,d\psi^0-\mathcal X\,d\psi^1\bigr)
\end{equation}
satisfy the second. Since
\begin{equation}\label{eq:norm-conf}
F'_1F'_2=\lambda^{-1}\det B\,F_1F_2,\qquad \mathfrak c'F'_1F'_2\,(d\rho'^2+dz'^2)=\lambda\,\mathfrak cF_1F_2\,(d\rho^2+dz^2),
\end{equation}
both terms of \eqref{eq:closed} are multiplied by $\lambda$. The map
$d\psi'\mapsto d\psi$ takes the lattice generated by $\mathfrak c'^{1/2}V'_0,\mathfrak c'^{1/2}V'_1$ to that generated by $\mathfrak c^{1/2}V_0,\mathfrak c^{1/2}V_1$.
\end{proof}

By Lemma~\ref{lem:normalise}, an affine change $z\mapsto\lambda z+\mu$ of the axis, together with a rescaling of $\omega$ and
$\gamma$, changes the metric only by a constant factor and a change of basis of $\mathfrak t$. In particular, soliton data with the degree
index $r=2$, the one that occurs below, can be normalised by $w_1=0$, $w_2=1$, and $a_1=1$. The soliton data are produced in
\cite{Gerochcf} by a surjectivity argument, and whether the soliton data of a metric are determined by it up to these changes is left open
in \cite[Rem.~9.12]{Gerochcf}. For data with given turning points, the lines of the rod vectors determine the soliton data up to a common
factor of the weights and $\gamma$, for ALE, AF, and ALF ends \cite{ADNprep}. For the instanton of
Theorem~\ref{thm:main}(c) the uniqueness follows from the proof in \S\ref{subsec:ALF-proof}, and for those of~(d) it is the third
statement of Remark~\ref{rem:AF-family}. Theorem~\ref{thm:main} classifies metrics, not soliton data, and its proof does not need this
uniqueness.

The entries of $V_j$ are polynomials with integer coefficients in the turning points, the weights $a_k$, and $\gamma$. Smoothness of the
metric is a condition on the stabilisers of consecutive rods (\S\ref{subsec:rods}), and by Theorem~\ref{thm:soliton} it becomes a system
of polynomial equations with integer coefficients for the soliton data. By Theorem~\ref{thm:converse}, every solution with positive
weights, $\gamma>0$, and $1\le r\le n-1$ defines a smooth complete toric gravitational instanton. For each chain, type of end, and degree
index the instantons therefore correspond, up to scale, to the positive real solutions of one polynomial system. For a given number of
turning points there are finitely many such systems, since the hyperk\"ahler instantons are known (\S\ref{subsec:literature}) and the
chains of the others take values in $\{-1,0,1\}^{n-1}$ by Li and Sun's bound (Theorem~\ref{thm:LSbound}).

\subsection{The main result}\label{subsec:main}
Theorem~\ref{thm:mainA} of the introduction is the following theorem without the data of the instantons. Some of these data, such as
the chain and the sign of an ALF end, depend on the orientation of $M$, and the theorem states them for a suitable orientation. By Lemma~\ref{lem:stab}
the notion of a toric instanton used in the introduction, an effective isometric action of $T^2$, agrees for simply connected
instantons with that of \S\ref{subsec:rods}.

\begin{theorem}\label{thm:main}
Let $(M,g)$ be a simply connected toric gravitational instanton whose torus action has four fixed points. Then, up to scale and
possibly after a change of the orientation of $M$, $(M,g)$ is one of the following.
\begin{enumerate}[leftmargin=*,label=\textup{(\alph*)}]
\item A multi-Eguchi--Hanson metric with four collinear centres, on the minimal resolution of $\C^2/\Z_4$. It is hyperk\"ahler, its end
is ALE, and its chain is $(2,2,2)$ up to orientation.
\item A multi-Taub--NUT metric with four collinear centres. It is hyperk\"ahler, its end is ALF, and its chain is $(2,2,2)$ up to
orientation.
\item The ALF instanton $g_{\rm ALF}$ with the chain $(-1,1,-1)$, the degree index $r=2$, and the normalised soliton data
\begin{equation}\label{eq:ALF-data}
\begin{gathered}
w=\Bigl(0,\,1,\,\tfrac12\sqrt{10},\,1+\tfrac12\sqrt{10}\Bigr),\qquad a=(1,A,A,1),\qquad A=\tfrac19\bigl(7+2\sqrt{10}\bigr),\\
\gamma_{\rm ALF}=\tfrac1{82}\bigl(40\sqrt2+36\sqrt5-19\sqrt{10}-65\bigr)\approx0.146143 ,
\end{gathered}
\end{equation}
which lie in $\F=\Q(\sqrt2,\sqrt5)$. It is Hermitian for neither orientation, and its end is ALF$^-$, the quotient of the end of the
Taub--NUT metric by $\Z_3$ acting on the base $\R^3$.
\item An AF instanton on $X_3$ with the chain $(-1,0,1)$ in a suitable orientation of $\mathfrak t$. It is Hermitian for neither
orientation. In this orientation its twist $\beta$, taken in $[0,1)$ as in \S\ref{subsec:rods}, is not $0$. For each $\beta\in(0,1)$
there is, up to scale, exactly one
simply connected toric gravitational instanton with four fixed points and an AF end that has, in one orientation of $\mathfrak t$,
the chain $(-1,0,1)$ and the twist $\beta$. Reversing the orientation of $\mathfrak t$ gives the chain $(1,0,-1)$ and replaces the twist $\beta$ by
$1-\beta$.
\end{enumerate}
Conversely, each of these metrics is a simply connected toric gravitational instanton with four fixed points. In each case it is given by the closed
formula \eqref{eq:closed} with soliton data, which Section~\ref{sec:explicit} lists for (a)--(c).
\end{theorem}
Teo found a two-parameter family of AF instantons on $X_3$, one of the parameters being the scale \cite{Teo2026}. We call it
\emph{Teo's family}, and a smooth complete instanton that is, up to scale, one of its metrics a \emph{member} of it. Each metric of Teo's
family whose parameters are finite and come from soliton data with positive weights is a member, by Theorem~\ref{thm:converse} and the
relation between Teo's metric and the closed formula in \S\ref{subsec:AF-sym}. The instantons of~(d) are the members of Teo's family.
This identification is not part of Theorem~\ref{thm:main}. We outline it in \S\ref{subsec:AF-teo}, using the closed formula, Chen's
construction of Teo's metric \cite{Chen2016,Teo2026}, and a chart of Aksteiner \cite{Aksteiner2026}.

\begin{figure}[ht]
\centering
\begin{minipage}[b]{0.31\textwidth}\centering
\begin{tikzpicture}[scale=0.399,>=stealth]
    \fill[gray!50,opacity=0.3](0,8.405)--(0,0)--(2.5,0)--(4.736,1.118)--(6.816,2.505)--(10.016,4.905)--(10.016,8.405)--cycle;
    \draw (0,8.405)--(0,0)--(2.5,0)--(4.736,1.118)--(6.816,2.505)--(10.016,4.905);
    \draw[->] (0,4.202)--(0.9,4.202) node[anchor=180,inner sep=1.5pt] {$v_{0}$};
    \draw[->] (0.75,0)--(0.75,0.9) node[anchor=270,inner sep=1.5pt] {$v_{1}$};
    \draw[->] (3.618,0.559)--(3.216,1.364) node[anchor=297,inner sep=1.5pt] {$v_{2}$};
    \draw[->] (5.776,1.811)--(5.277,2.56) node[anchor=304,inner sep=1.5pt] {$v_{3}$};
    \draw[->] (8.416,3.705)--(7.876,4.425) node[anchor=307,inner sep=1.5pt] {$v_{4}$};
\end{tikzpicture}

\smallskip\small (a), (b) $\kappa=(2,2,2)$
\end{minipage}\hfill
\begin{minipage}[b]{0.31\textwidth}\centering
\begin{tikzpicture}[scale=0.440,>=stealth]
    \fill[gray!50,opacity=0.3](0,3.5)--(0,0)--(2.5,0)--(3.527,-1.027)--(5.763,-2.145)--(9.092,-4.364)--(9.092,3.5)--cycle;
    \draw (0,3.5)--(0,0)--(2.5,0)--(3.527,-1.027)--(5.763,-2.145)--(9.092,-4.364);
    \draw[->] (0,2.5)--(0.9,2.5) node[anchor=180,inner sep=1.5pt] {$v_{0}$};
    \draw[->] (0.75,0)--(0.75,0.9) node[anchor=270,inner sep=1.5pt] {$v_{1}$};
    \draw[->] (3.014,-0.514)--(2.377,-1.15) node[anchor=45,inner sep=1.5pt] {$v_{2}$};
    \draw[->] (4.645,-1.586)--(4.243,-2.391) node[anchor=63,inner sep=1.5pt] {$v_{3}$};
    \draw[->] (7.427,-3.255)--(7.927,-2.506) node[anchor=236,inner sep=1.5pt] {$v_{4}$};
\end{tikzpicture}

\smallskip\small (c) ALF, $\kappa=(-1,1,-1)$
\end{minipage}\hfill
\begin{minipage}[b]{0.31\textwidth}\centering
\begin{tikzpicture}[scale=0.450,>=stealth]
    \fill[gray!50,opacity=0.3](0,9.268)--(0,0)--(2.5,0)--(4.268,1.768)--(6.768,1.768)--(6.768,5.768)--(6.768,9.268)--cycle;
    \draw (0,9.268)--(0,0)--(2.5,0)--(4.268,1.768)--(6.768,1.768)--(6.768,5.768);
    \draw[->] (0,4.634)--(0.9,4.634) node[anchor=180,inner sep=1.5pt] {$v_{0}$};
    \draw[->] (0.75,0)--(0.75,0.9) node[anchor=270,inner sep=1.5pt] {$v_{1}$};
    \draw[->] (3.384,0.884)--(2.747,1.52) node[anchor=315,inner sep=1.5pt] {$v_{2}$};
    \draw[->] (5.518,1.768)--(5.518,0.868) node[anchor=90,inner sep=1.5pt] {$v_{3}$};
    \draw[->] (6.768,3.768)--(7.668,3.768) node[anchor=180,inner sep=1.5pt] {$v_{4}$};
\end{tikzpicture}

\smallskip\small (d) AF, $\kappa=(1,0,-1)$
\end{minipage}
\caption{The rod structures of Theorem~\ref{thm:main}, drawn as by Li and Sun \cite{LS2025}. The shaded region is the orbit space
$M/T^2$. Each finite side is a rod, drawn to scale and orthogonal to its rod vector, and the arrows are the rod vectors $v_0,\dots,v_4$.
(a), (b) Multi-Eguchi--Hanson and multi-Taub--NUT, with free rod lengths drawn equal. (c) $g_{\rm ALF}$. (d) An AF instanton on
$X_3$.}\label{fig:rods}
\end{figure}
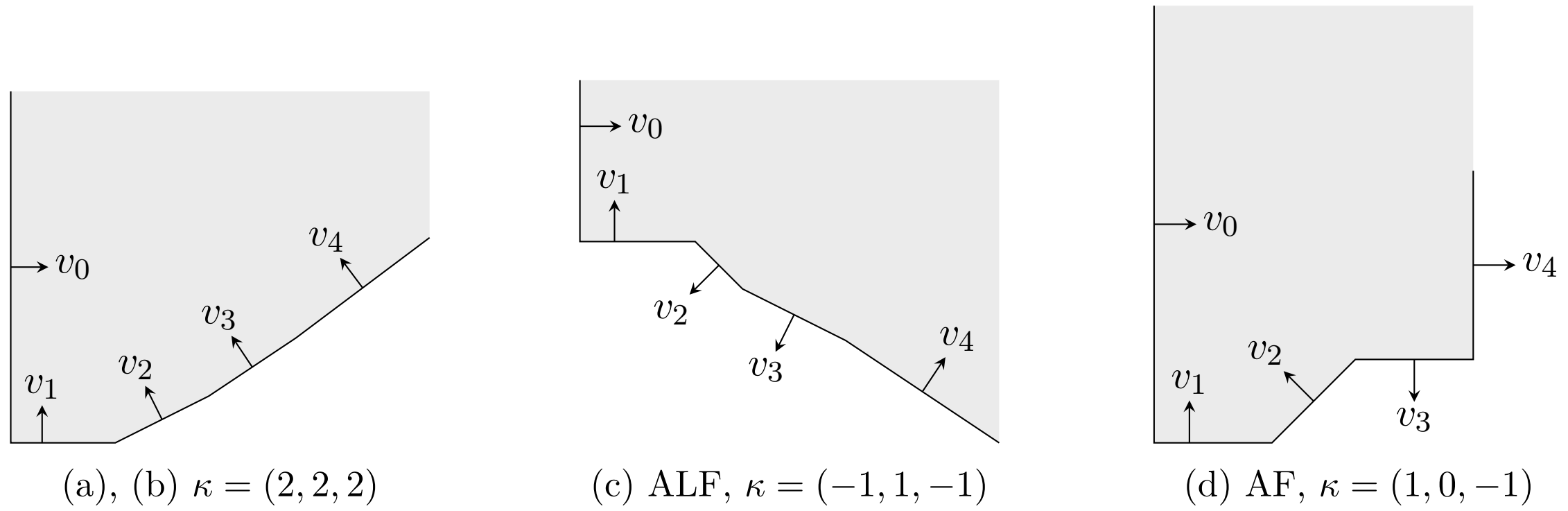

Figure~\ref{fig:rods} shows the rod structures of (a)--(d). It uses the basis $v_0=(1,0)$, $v_1=(0,1)$ of the lattice, which is negatively
oriented by \eqref{eq:det}, and the first rod has length $1$. It shows (a), (b) in the orientation opposite to the hyperk\"ahler one,
with the rod vectors $v_j=(1-j,j)$ and the chain $(2,2,2)$. In (c) the rod vectors are \eqref{eq:ALF-rods} and the rod lengths are
$1$, $\frac12\sqrt{10}-1$, $1$, so that the outer rods are not parallel and $|\det(v_0,v_4)|=3$. In (d) the orientation of $\mathfrak t$
is the one in which the chain is $(1,0,-1)$, the rod vectors are those of \eqref{eq:AF-rods} in reversed order, $v_4=v_0$, and the rods
are drawn with equal lengths. That (a) and (b) are the toric hyperk\"ahler instantons with four turning points is known (\S\ref{subsec:literature}), and the proof
identifies them in the soliton representation (Proposition~\ref{prop:reduction}(i), Remark~\ref{rem:hk}).

The instantons that are hyperk\"ahler for neither orientation are those of (c) and (d). That none of them has an ALE end also follows
from the classification of the toric Hermitian non-K\"ahler ALE instantons by Araneda and Lucietti \cite[Thm.~1.2]{AL2025}
(Remark~\ref{rem:reduction}). The families (a) and (b) are parametrised by the positions of the centres, and in (b) also by the size of
the circle at infinity. The instanton (c) is a single metric, and (d) is, up to scale, a one-parameter family parametrised by the twist,
which takes every value except $0$. Indeed, by \cite[Cor.~11.2]{Gerochcf}, a simply connected toric gravitational instanton with finitely
many fixed points, at least one, and an end of type AF$_0$ has two turning points and is the Euclidean Schwarzschild metric. This proves
Li and Sun's Conjecture~8.1 \cite{LS2025} for toric instantons.

The instanton (c) is the one of \cite[Rem.~11.5(a)]{Gerochcf}, where its data and its existence are given, and
Proposition~\ref{prop:ALF-slice} derives its data. Its rod structure is the entry ALF4-001 of Aksteiner's catalogue of ALF rod structures
\cite[App.~B, Table~30]{Aksteiner2026}. Jones constructed it independently \cite{Jones2026}. His metric is Teo's metric
\cite[(1)--(3)]{Teo2026} with four turning points placed symmetrically about $z=0$ and Teo's parameters (\S\ref{subsec:AF-sym})
satisfying $C_3=1/C_2$ and $C_4=1/C_1$ \cite[(1)--(3), App.~A]{Jones2026}, and Jones determines these parameters in closed form
\cite[(7)]{Jones2026}. He proves that his metric is a toric gravitational instanton with four turning points and an ALF end with the lens
space $L(3,1)$ at infinity. It is simply connected, since consecutive rod vectors in \cite[Fig.~1]{Jones2026} generate the lattice. By
Theorem~\ref{thm:main} it is therefore $g_{\rm ALF}$ up to scale and orientation, since the ends of (a) and (d) are ALE and AF and the
boundary at infinity of (b) is $L(4,1)$ (Table~\ref{tab:topology}).

The proof of Theorem~\ref{thm:main} combines a reduction of the possible rod structures with a convexity principle for the cone
angles. The reduction (Section~\ref{sec:reduction}) uses Li and Sun's bound $|\kappa_j|\le1$ on the chain of an instanton which is
hyperk\"ahler for neither orientation (Theorem~\ref{thm:LSbound}) and an elementary inequality for the rod vectors that follows from
Heine's formula. For ALE ends the inequality bounds both coordinates of the interior rod vectors and leaves no rod structure. For ALF and
AF ends it bounds one coordinate and leaves one class of chains each.

For the ALF chain (Section~\ref{sec:ALF}), the solutions of the smoothness system that are symmetric under the reflection of the axis
are found by hand, and exactly one of them has positive data. That every positive solution is symmetric
(Proposition~\ref{prop:ALF-symmetric}) rests on convexity. For a fixed asymptotic parameter of the end, consider Li and Sun's canonical
harmonic maps of the ALF rod structure and their \emph{cone errors}, the logarithms of the ratios of $2\pi$ to the cone angles along the
finite rods. As for Aksteiner's potential for AF solitons \cite{Aksteiner2026}, the cone errors are the gradient of a strictly convex
function of the rod lengths, and the soliton data of every positive solution give a rod structure at which they vanish
(Appendix~\ref{app:convex}). The cone errors are invariant under a common scaling of the rod lengths and the asymptotic parameter, and
the \emph{balanced} pairs of rod lengths and asymptotic parameter, those at which all cone errors vanish, form a single ray
(Lemma~\ref{lem:appC-scale}, Corollary~\ref{cor:appC-ray}). Every positive solution, with its first rod length normalised to $1$,
therefore has the rod lengths of the symmetric one, and a fold in the first components of the smoothness system
(Appendix~\ref{app:ALF}) forces the symmetry.

For the AF chain (Section~\ref{sec:AF}), the convexity argument of Appendix~\ref{app:convex}, which holds for every constant asymptotic
parameter $\alpha\ge0$, applies at $\alpha=0$ with the twist fixed. It shows that two AF instantons with the same twist agree up to scale
(Proposition~\ref{prop:AF-main}). The existence of an instanton for every twist in $(0,1)$ is Li and Sun's existence theorem
\cite[Thm.~1.1]{LS2025}.

\section{The reduction}\label{sec:reduction}
A simply connected toric gravitational instanton with four turning points that is hyperk\"ahler for neither orientation has no ALE end. If its end is ALF, its chain is $\pm(-1,1,-1)$, and if its end is AF, its chain lies in
the class of $(1,0,-1)$ (Proposition~\ref{prop:reduction}). The hypothesis that the metric is hyperk\"ahler for neither orientation enters twice. By Theorems~\ref{thm:LSdegree} and~\ref{thm:soliton}(i) it excludes the degree
indices $r$ for which one of the two degrees vanishes, and by Theorem~\ref{thm:LSbound} it bounds the chain. The remaining input is a sign pattern in Heine's formula for the Hankel determinants of the measures $\omega_{(j)}$.

Heine's formula \eqref{eq:heine} writes $H_s(\omega_{(j)})$ as the sum of the positive terms of $H_s(\omega)$, each with a sign. The resulting inequality appears in the proof of \cite[Cor.~11.2]{Gerochcf}, see
\cite[(11.5)]{Gerochcf}.

\begin{lemma}\label{lem:strip}
Let $\omega$ be a positive measure on $n\ge1$ points, $\omega_{(j)}$ as in \eqref{eq:V}, and $0\le s\le n$. Then $H_s(\omega)>0$, and the
following hold.
\begin{enumerate}[leftmargin=*,label=\textup{(\roman*)}]
\item $H_s(\omega_{(n)})=(-1)^sH_s(\omega)$, and $H_n(\omega_{(j)})=(-1)^jH_n(\omega)$ for $0\le j\le n$.
\item If $1\le s\le n-1$, then $|H_s(\omega_{(j)})|<H_s(\omega)$ for $0<j<n$.
\end{enumerate}
\end{lemma}

\begin{proof}
By \eqref{eq:PS}, $P_S(\omega_{(j)})=(-1)^{|S\cap\{1,\dots,j\}|}P_S(\omega)$ with $P_S(\omega)>0$. By Heine's formula \eqref{eq:heine}, $H_s(\omega)$
is a sum of positive terms and $H_s(\omega_{(j)})$ is the sum of the same terms with these signs. For $j=n$ every set of $s$ atoms has the
sign $(-1)^s$, and for $s=n$ the only set is $\{1,\dots,n\}$, with the sign $(-1)^j$. This gives (i). For (ii) the triangle inequality gives
$|H_s(\omega_{(j)})|\le H_s(\omega)$, with equality only if all signs agree. Let $0<j<n$ and $1\le s\le n-1$. There is a set
$T\subset\{2,\dots,n-1\}$ of $s-1$ atoms. The sets $T\cup\{1\}$ and $T\cup\{n\}$ carry opposite signs, since the atom $1$ is reversed in
$\omega_{(j)}$ and the atom $n$ is not. Hence the inequality is strict.
\end{proof}

Applied to the vectors $V_j$, Lemma~\ref{lem:strip} constrains the rod vectors. Let $(M,g)$ be as in \S\ref{subsec:rods}, and let $r,\omega,\gamma,\mathfrak c$ be as in
Theorem~\ref{thm:soliton}. Write $V_j=(V^{(1)}_j,V^{(2)}_j)$ and
\begin{equation}\label{eq:x}
x_j=\frac{V^{(1)}_j}{V^{(1)}_0}=\frac{H_r(\omega_{(j)})}{H_r(\omega)},\qquad 0\le j\le n .
\end{equation}
We call $x_j$ the \emph{angular component} of the rod $I_j$.

\begin{corollary}\label{cor:strip}
Let $(M,g)$, $r$, $\omega$, $\gamma$, and $\mathfrak c$ be as in Theorem~\ref{thm:soliton}, and $x_j$ as in \eqref{eq:x}. Then the following
hold.
\begin{enumerate}[leftmargin=*,label=\textup{(\alph*)}]
\item $x_0=1$, $x_{j+1}+x_{j-1}=\kappa_jx_j$ for $1\le j\le n-1$, and $x_n=(-1)^r$.
\item If $1\le r\le n-1$, then $|x_j|<1$ for $0<j<n$.
\item If $r=0$, then $\kappa=(2,\dots,2)$. If $r=n$, or if the end is ALE and $r=n-1$, then $\kappa=(-2,\dots,-2)$.
\item If the end is ALE and $1\le r\le n-2$, then $p\ne0$ and
\begin{equation}\label{eq:parallelogram}
|\det(v_j,v_n)|+|\det(v_0,v_j)|<|p|,\qquad 0<j<n ,
\end{equation}
that is, the interior rod vectors lie in the interior of the parallelogram with the vertices $\pm v_0$, $\pm v_n$.
\end{enumerate}
\end{corollary}

\begin{proof}
By Theorem~\ref{thm:soliton}(iii) the vectors $V_j=\mathfrak c^{-1/2}v_j$ satisfy \eqref{eq:chain}, hence so do their first components.
Lemma~\ref{lem:strip} with $s=r$ gives (a) and (b). If $r=0$, all $V^{(1)}_j=H_0=1$, and \eqref{eq:chain} gives $\kappa_j=2$. If $r=n$, then
$x_j=(-1)^j$ by Lemma~\ref{lem:strip}(i), and \eqref{eq:chain} gives $\kappa_j=-2$. For an ALE end $\gamma=0$ and
$V^{(2)}_j=H_{r+1}(\omega_{(j)})$. If $r=n-1$, Lemma~\ref{lem:strip}(i) gives $V^{(2)}_j=(-1)^jV^{(2)}_0$ and again $\kappa_j=-2$. For (d) write
$H_s=H_s(\omega)>0$. By Lemma~\ref{lem:strip}(i), $V_0=(H_r,H_{r+1})$ and $V_n=\bigl((-1)^rH_r,(-1)^{r+1}H_{r+1}\bigr)$, so
$p=\mathfrak c\det(V_0,V_n)=2(-1)^{r+1}\mathfrak c\,H_rH_{r+1}\ne0$. Write $v_j=m_jv_0+m'_jv_n$, so that $m_j=\det(v_j,v_n)/p$ and $m'_j=\det(v_0,v_j)/p$.
Comparing components gives $H_r(\omega_{(j)})=(m_j+(-1)^rm'_j)H_r$ and $H_{r+1}(\omega_{(j)})=(m_j-(-1)^rm'_j)H_{r+1}$.
Lemma~\ref{lem:strip}(ii) with $s=r$ and $s=r+1$ then gives $|m_j\pm m'_j|<1$, that is, $|m_j|+|m'_j|<1$.
\end{proof}

For four turning points and a chain with entries in $\{-1,0,1\}$, the conditions of Corollary~\ref{cor:strip}(a),(b) on the angular components leave two possibilities.

\begin{lemma}\label{lem:seq}
Let $\kappa\in\{-1,0,1\}^3$, and let $x_0,\dots,x_4$ be real numbers with $x_0=1$, $x_{j+1}+x_{j-1}=\kappa_jx_j$ for $j=1,2,3$,
$|x_1|,|x_2|,|x_3|<1$, and $|x_4|=1$. Then $x_4=1$, and with $\varepsilon\in\{1,-1\}$ one of the following holds.
\begin{enumerate}[leftmargin=*,label=\textup{(\alph*)}]
\item $\kappa=(\varepsilon,0,-\varepsilon)$, $0<\varepsilon x_1<1$, $x_2=\varepsilon x_1-1$, $x_3=-x_1$.
\item $\kappa=(\varepsilon,-\varepsilon,\varepsilon)$ and $(x_1,x_2,x_3)=(\frac\varepsilon3,-\frac23,\frac\varepsilon3)$.
\end{enumerate}
\end{lemma}

\begin{proof}
The relation for $j=1$ gives $x_2=\kappa_1x_1-1\in(-1,0)$, and the one for $j=3$ gives $x_4=\kappa_3x_3-x_2>-1$, so $x_4=1$. Put
$s=1+x_2\in(0,1)$. The relations for $j=1,3$ become $\kappa_1x_1=s=\kappa_3x_3$, so $\kappa_1,\kappa_3\in\{\pm1\}$, $x_1=\kappa_1s$, and
$x_3=\kappa_3s$. The relation for $j=2$ becomes $(\kappa_1+\kappa_3)s=-\kappa_2(1-s)$. If $\kappa_3=-\kappa_1$, then $\kappa_2=0$, which is the
first case with $\varepsilon=\kappa_1$. If $\kappa_3=\kappa_1=\varepsilon$, then $\kappa_2\ne0$, $2s=1-s$, $s=\frac13$, and $\kappa_2=-\varepsilon$, which
is the second case.
\end{proof}

Together with the bound of Theorem~\ref{thm:LSbound}, these results reduce the chains as follows.

\begin{proposition}\label{prop:reduction}
Let $(M,g)$ be a simply connected toric gravitational instanton whose torus action has four fixed points, and let $r$ be as in
Theorem~\ref{thm:soliton}, $\kappa$ as in \eqref{eq:chain}, $p$ as in \eqref{eq:v4}, and $x_j$ as in \eqref{eq:x}. Then the following hold.
\begin{enumerate}[leftmargin=*,label=\textup{(\roman*)}]
\item If $g$ is hyperk\"ahler for one orientation, then its end is ALE or ALF, $g$ is, up to scale, a multi-Eguchi--Hanson or a
multi-Taub--NUT metric with four collinear centres, and $\kappa=\pm(2,2,2)$.
\item If $g$ is hyperk\"ahler for neither orientation, then the end is not ALE, $r=2$, and both Li--Sun degrees are $2$.
\item In case (ii) with an ALF end, $\kappa=\pm(-1,1,-1)$ and $|p|=3$.
\item In case (ii) with an AF end, $\kappa=(\varepsilon,0,-\varepsilon)$ with $\varepsilon\in\{1,-1\}$, and $0<\varepsilon x_1<1$.
\end{enumerate}
\end{proposition}

\begin{proof}
By Theorem~\ref{thm:soliton}(i) and Theorem~\ref{thm:LSdegree}, $g$ is hyperk\"ahler for one orientation if and only if $r\in\{0,3\}$ (ALE) or
$r\in\{0,4\}$ (ALF, AF). Theorem~\ref{thm:soliton}(i), applied to $g$ with both orientations, shows that after a change of the
orientation of $M$ we may assume $r=0$, with $r,\omega,\gamma,\mathfrak c$ given by Theorem~\ref{thm:soliton} for the new orientation.
The hyperk\"ahler orientation is then the opposite one. For $r=0$, $H_0=1$, and
$H^{(\gamma)}_1=H_1$ since the moment of order $2r-1$ does not occur \cite[(5.22)]{Gerochcf}, so \eqref{eq:V} gives
\begin{equation}\label{eq:V-hk}
V_j=\Bigl(1,\ \sum_{k>j}a_k-\sum_{k\le j}a_k\Bigr),\qquad 0\le j\le4 .
\end{equation}
By Theorem~\ref{thm:soliton}(iii) these vectors satisfy \eqref{eq:chain}: for $1\le j\le3$ the first components give $\kappa_j=2$, and the
second ones $a_{j+1}=a_j$. So all weights are equal, to $a$ say. Then $p=\mathfrak c\det(V_0,V_4)=-2\mathfrak c\sum_ka_k\ne0$, so the end is not AF. By
\cite[Prop.~11.3]{Gerochcf}, \eqref{eq:closed} with $r=0$, these weights, and the atom weight $\gamma$ is, up to a constant factor, the Gibbons--Hawking
metric with the centres $w_1,\dots,w_4$ and the potential $\gamma+\sum_ka/R_k$. Hence it is, up to scale, \eqref{eq:GH} with
$U_\infty=\gamma/(2a)$: the
multi-Eguchi--Hanson metric for $\gamma=0$ and the multi-Taub--NUT metric for $\gamma>0$. With the change of orientation,
$\kappa=\pm(2,2,2)$. This is (i).

Otherwise $1\le r\le2$ (ALE) or $1\le r\le3$ (ALF, AF), and $\kappa\in\{-1,0,1\}^3$ by Theorem~\ref{thm:LSbound}, so $|p|\le3$ by
\eqref{eq:v4}. Suppose the end is ALE. By \eqref{eq:det}, \eqref{eq:chain}, and \eqref{eq:v4}, $\det(v_0,v_1)=\det(v_3,v_4)=-1$,
$\det(v_1,v_4)=1-\kappa_2\kappa_3$, and $\det(v_0,v_3)=1-\kappa_1\kappa_2$, so \eqref{eq:parallelogram} for $j=1$ and $j=3$ gives
\begin{equation}\label{eq:ALEineq}
|p|\ge2+|1-\kappa_2\kappa_3|,\qquad |p|\ge2+|1-\kappa_1\kappa_2| .
\end{equation}
If $|p|=2$, then $\kappa_2\kappa_3=\kappa_1\kappa_2=1$, so $\kappa_1=\kappa_2=\kappa_3=\varepsilon$ and $p=2\varepsilon-\varepsilon^3=\varepsilon$, a
contradiction. If $|p|=3$, the three terms of $p$ in \eqref{eq:v4} equal $\varepsilon=p/3$, so $\kappa_1=\kappa_3=\varepsilon=-\kappa_2$.
Then $\kappa_2\kappa_3=-1$, and \eqref{eq:ALEineq} gives $|p|\ge4$, again a contradiction. Since $|p|\le1$ contradicts \eqref{eq:ALEineq} directly, the end is not ALE.

For ALF and AF ends, Corollary~\ref{cor:strip}(a),(b) show that $x_0,\dots,x_4$ satisfy the hypotheses of Lemma~\ref{lem:seq}. Hence $x_4=1$,
and $x_4=(-1)^r$ gives $r=2$. Both degrees are $2$ by Theorem~\ref{thm:soliton}(i). In the first case of Lemma~\ref{lem:seq},
\eqref{eq:v4} gives $v_4=v_0$ and $p=0$. In the second it gives $v_4=2v_0-3\varepsilon v_1$ and $p=3\varepsilon$. Since ALF ends have $p\ne0$ and
AF ends $p=0$, this proves (iii) and (iv).
\end{proof}

\begin{remark}\label{rem:reduction}
For ALE ends the reduction also follows from other results. By \eqref{eq:degid}, $\min(d,\bar d)\le1$, and degree $1$ means
Hermitian non-K\"ahler. By the theorem of Araneda and Lucietti, which needs no condition on the structure group at infinity, the
only smooth toric ALE instanton of this kind is the Eguchi--Hanson metric, with two turning points \cite[Thm.~1.2]{AL2025}. The
difference between ALE and ALF ends lies in the second component of $V_j$: for $\gamma>0$ it has no sign, and only the first component obeys
Lemma~\ref{lem:strip}. The ALF chain $(-1,1,-1)$ satisfies the inequality of Corollary~\ref{cor:strip}(b) while lying on the boundary of
the parallelogram \eqref{eq:parallelogram}, with equality for $j=1,3$.
\end{remark}

\begin{remark}[the hyperk\"ahler instantons]\label{rem:hk}
The proof of Proposition~\ref{prop:reduction}(i) applies to any number $n\ge2$ of turning points and describes the toric hyperk\"ahler
instantons in the soliton representation. They are those with the degree index $r=0$ for one orientation of $M$, the hyperk\"ahler
orientation being the opposite one. For $r=0$ all Hankel determinants in \eqref{eq:closed} have order at most one: $\mathcal X=F_1=1$,
$F_2=\gamma+\sum_ka_k/R_k$ is the Gibbons--Hawking potential, and the closed formula is the Gibbons--Hawking form
\cite[Prop.~11.3]{Gerochcf}. By \eqref{eq:V-hk}, with $n$ in place of $4$, the smoothness conditions \eqref{eq:chain} are linear in the
weights and involve neither the turning points nor $\gamma$. They say that $\kappa=(2,\dots,2)$ and that all weights are equal, and then
$p=-2\mathfrak c\sum_ka_k\ne0$, so that the end is ALE or ALF. This is the observation of Li and Sun that the scaled rod vectors of a
Gibbons--Hawking metric, and hence its cone angles, do not depend on the rod lengths \cite[proof of Thm.~7.7, Rem.~7.8]{LS2025}.
Every choice of turning points, of a common weight $a>0$, and of $\gamma\ge0$ therefore gives a smooth toric hyperk\"ahler instanton, whose
metric depends only on the turning points and on $\gamma/(2a)$. Up to scale these metrics form a family of dimension $n-2$ for ALE and
$n-1$ for ALF ends, the multi-Eguchi--Hanson and multi-Taub--NUT metrics with $n$ collinear centres. They are the toric ones among
Kronheimer's hyperk\"ahler ALE instantons of type $A_{n-1}$ \cite{Kronheimer1989a,Kronheimer1989b} and among the multi-Taub--NUT
metrics. By Minerbe, the multi-Taub--NUT metrics are the non-flat hyperk\"ahler ALF instantons with a circle fibration over the end of $\R^3$
\cite{Minerbe2011}. In the hyperk\"ahler orientation the degree index is $r=n-1$ (ALE), respectively $r=n$ (ALF), and
$\kappa=(-2,\dots,-2)$ by Corollary~\ref{cor:strip}(c), so that the intersection form $\mathrm{tridiag}(-2;1)$ of \eqref{eq:selfint} is minus
the Cartan matrix of $A_{n-1}$, as in Kronheimer's description. In this orientation the smoothness conditions are of higher degree in the
weights, but linear in the dual weights $a_k^\vee=1/(a_kq'(w_k)^2)$, $q(T)=\prod_k(T-w_k)$, of \eqref{eq:duality}. The identity
$c_\omega P_S(\omega^\vee_{(j)})=(-1)^jP_{S^c}(\omega_{(j)})$ of the proof of Lemma~\ref{lem:duality}, with
$c_\omega=\prod_ka_k\prod_{k<l}(w_k-w_l)^2$, holds for any number $n$ of atoms. For ALE ends $r=n-1$ and $\gamma=0$. For ALF ends $r=n$,
$H_{n+1}=0$, and $\Sigma_n(\omega_{(j)})=(-1)^jc_\omega\mathfrak s$, with $\Sigma_n$ as in \eqref{eq:H3gamma} and $\mathfrak s=\sum_kw_k$
\cite[(5.20)]{Gerochcf}. With \cite[(5.22)]{Gerochcf} for ALF ends, the identity therefore gives
\begin{equation}\label{eq:hk-dual}
V_j=(-1)^jc_\omega\bigl(H_1(\omega^\vee_{(j)}),\,1\bigr)\quad\text{(ALE)},\qquad
V_j=(-1)^jc_\omega\bigl(1,\,-2\gamma\mathfrak s-\gamma^2H_1(\omega^\vee_{(j)})\bigr)\quad\text{(ALF)}.
\end{equation}
As for \eqref{eq:V-hk}, the smoothness conditions then say that $\kappa=(-2,\dots,-2)$ and that $a_kq'(w_k)^2$ does not depend on $k$. For $1\le r\le n-2$ (ALE), respectively $1\le r\le n-1$ (ALF, AF), the Hankel determinants have higher order and the smoothness
conditions are not linear. For $n=4$, Proposition~\ref{prop:reduction} leaves only $r=2$, with one positive solution for ALF ends (Section~\ref{sec:ALF})
and a one-parameter family of them for AF ends (Section~\ref{sec:AF}).
\end{remark}

\section{The ALF instanton: proof of Theorem~\ref{thm:main}(c)}\label{sec:ALF}
Up to normalisation, the smoothness system for the chain $(-1,1,-1)$ has exactly one positive solution, and this solution is symmetric under the
reflection of the axis. It gives part (c) of Theorem~\ref{thm:main}, and with the opposite orientation of $M$ it gives the soliton data
for the chain $(1,-1,1)$.

\subsection{The smoothness system}\label{subsec:ALF-system}
The smoothness conditions for the chain $(-1,1,-1)$ are polynomial equations. By Proposition~\ref{prop:reduction}(iii),
after a change of the orientation of $M$ the chain is $\kappa=(-1,1,-1)$ and $r=2$. By Theorem~\ref{thm:soliton}(iii) the soliton data
$(\omega,\gamma)$, with $\omega=\sum_{k=1}^4a_k\delta_{w_k}$ positive and $\gamma>0$, satisfy
\begin{equation}\label{eq:ALF-lin}
V_{j+1}+V_{j-1}=\kappa_jV_j\quad(j=1,2,3)\iff V_2=-V_0-V_1,\quad V_3=-V_0-2V_1,\quad V_4=2V_0+3V_1,
\end{equation}
with $V_j=\bigl(H_2(\omega_{(j)}),\,H_3(\omega_{(j)})-2\gamma\Sigma_2(\omega_{(j)})-\gamma^2H_1(\omega_{(j)})\bigr)$ by \eqref{eq:H3gamma}.
The chain forces the end to be ALF$^-$: $p=\det(v_0,v_4)=-3$ by \eqref{eq:v4}, while $p=\mathfrak c\det(V_0,V_4)=-\mathfrak c\gamma\,\alpha(\omega,\gamma)H_2(\omega)^2$
by \cite[(5.29)]{Gerochcf} and \eqref{eq:alpha}. Hence $\alpha(\omega,\gamma)>0$ (Theorem~\ref{thm:converse}). By
Lemma~\ref{lem:normalise}, an affine change of the axis and a rescaling of the data map solutions of \eqref{eq:ALF-lin} to solutions and
change the metric only by a constant factor. Exactly one such change gives $w_1=0$, $w_2=1$, and $a_1=1$. In these \emph{normalised} data
the unknowns are
\begin{equation}\label{eq:ALF-unknowns}
\mathbf d=(L_2,L_3,a_2,a_3,a_4,\gamma),\qquad w=(0,\,1,\,1+L_2,\,1+L_2+L_3),\qquad a=(1,a_2,a_3,a_4),
\end{equation}
and the data are positive if and only if all six coordinates of $\mathbf d$ are. In these unknowns \eqref{eq:ALF-lin} is a system of six polynomial equations
$f_{j,s}(\mathbf d)=0$, $j=1,2,3$, $s=1,2$, with integer coefficients: $f_{j,s}$ is the $s$-th component of $V_{j+1}+V_{j-1}-\kappa_jV_j$.

\subsection{The reflection and the symmetric solutions}\label{subsec:ALF-sym}
The reflection $z\mapsto-z$ of the axis
reverses the order of the turning points: the reflected measure $\omega^\flat=\sum_ka_{5-k}\delta_{-w_{5-k}}$, with its atoms
numbered in increasing order, has $\omega^\flat_{(j)}$ equal to the reflection of $-\omega_{(4-j)}$. By \eqref{eq:PS}, $P_S$ involves the turning
points only through squared differences and satisfies $P_S(-\tilde\omega)=(-1)^{|S|}P_S(\tilde\omega)$, while $w_k+w_l$ changes sign under the
reflection. Hence $H_2(\omega^\flat_{(j)})=H_2(\omega_{(4-j)})$, while $H_1$, $H_3$, and $\Sigma_2$ change sign. Since $\kappa$ is a palindrome,
the reflected data satisfy \eqref{eq:ALF-lin} again, with the vectors $\mathrm{diag}(1,-1)V_{4-j}$. After the normalisation of
\S\ref{subsec:ALF-system}, with $\lambda=L_3^{-1}$, $\mu=(1+L_2+L_3)/L_3$, and $K=a_4^{-1}$, the reflection acts on the unknowns by
\begin{equation}\label{eq:sigma}
\sigma(L_2,L_3,a_2,a_3,a_4,\gamma)=\Bigl(\frac{L_2}{L_3},\ \frac1{L_3},\ \frac{a_3}{a_4},\ \frac{a_2}{a_4},\ \frac1{a_4},\
\frac{\gamma}{a_4L_3^3}\Bigr).
\end{equation}
It is an involution of the points with $L_3a_4\ne0$, and it maps positive points to positive points and solutions to solutions, also
complex ones. Indeed, the identities for the reflection above and the relation $V'_j=BV_j$ with \eqref{eq:normB} in the proof of
Lemma~\ref{lem:normalise} are polynomial identities in the data and in $\lambda,\mu,K$. So they hold for complex data with $\lambda K\ne0$,
although the lemma is stated for positive data. We call the fixed points of $\sigma$ the \emph{symmetric solutions}. A point with non-zero coordinates is fixed if and only if $L_3=1$, $a_4=1$, and $a_2=a_3$. Indeed, a fixed point has
$L_3^2=1$ and $a_4^2=1$. Since $L_3=-1$ would give $L_2=-L_2$, we have $L_3=1$. Then $\gamma=\gamma/a_4$ gives $a_4=1$, and hence $a_2=a_3$. We call the set
\begin{equation}\label{eq:locus}
L_3=1,\qquad a_4=1,\qquad a_2=a_3
\end{equation}
the \emph{symmetric locus}. Its points are the normalised data whose turning points and weights are symmetric under the reflection.

\begin{proposition}\label{prop:ALF-slice}
On the symmetric locus write $u=L_2$ and $A=a_2=a_3$. The solutions of \eqref{eq:ALF-lin} on it with non-zero $A$ and pairwise
distinct turning points are the four points with
\begin{equation}\label{eq:slice-sol}
2u^2+4u-3=0,\qquad A=\frac{u+2}{3u}=\frac1{2u^2},\qquad (2A+1)\gamma^2+A(2u+7)\gamma-\tfrac54A=0 .
\end{equation}
Exactly one of them is positive: $u=\frac12\sqrt{10}-1$, $A=\frac19(7+2\sqrt{10})$, and $\gamma=\gamma_{\rm ALF}$ of \eqref{eq:ALF-data}. The four points
are conjugate under the Galois group of $\F=\Q(\sqrt2,\sqrt5)$, and their values of $\gamma$ are the roots of
\begin{equation}\label{eq:E4}
\mathcal E_4(T)=656T^4+2080T^3-120T^2-200T+25 .
\end{equation}
\end{proposition}

\begin{proof}
On the locus the turning points are $0,1,1+u,2+u$ and the weights $1,A,A,1$, so by \eqref{eq:PS}
\begin{equation}\label{eq:slice-PS}
\begin{gathered}
P_{12}=P_{34}=A,\qquad P_{13}=P_{24}=A(u+1)^2,\qquad P_{14}=(u+2)^2,\qquad P_{23}=A^2u^2,\\
P_{123}=P_{234}=A^2u^2(u+1)^2,\qquad P_{124}=P_{134}=A(u+1)^2(u+2)^2,
\end{gathered}
\end{equation}
where $P_{kl}=P_{\{k,l\}}(\omega)$ and $P_{klm}=P_{\{k,l,m\}}(\omega)$. The sums $w_k+w_l$ are $1$, $u+1$, $u+2$, $u+2$, $u+3$, $2u+3$ for
$kl=12,13,14,23,24,34$. Inserting these into \eqref{eq:ALF-lin} with the signs $P_S(\omega_{(j)})=(-1)^{|S\cap\{1,\dots,j\}|}P_S(\omega)$ gives
\begin{equation}\label{eq:slice-f}
\begin{aligned}
f_{1,1}&=f_{3,1}=A^2u^2+4A-(u+2)^2,\qquad f_{2,1}=(3Au-u-2)(Au+u+2),\\
f_{1,2}&=2A^2u^2(u+1)^2-2\gamma(u+2)\bigl(f_{1,1}+2A(u+1)\bigr)-2(2A+1)\gamma^2,\\
f_{2,2}&=-2\gamma(u+2)f_{2,1},\qquad f_{3,2}=-f_{1,2}-4\gamma(u+2)f_{1,1} .
\end{aligned}
\end{equation}
So the six equations on the locus are equivalent to $f_{1,1}=f_{2,1}=0$ and
\begin{equation}\label{eq:slice-g0}
A^2u^2(u+1)^2-2A(u+1)(u+2)\gamma-(2A+1)\gamma^2=0 .
\end{equation}
If $Au+u+2=0$, then $f_{1,1}=4A\ne0$. Hence $3Au=u+2$, and then $f_{1,1}=-4(u+2)(2u^2+4u-3)/(9u)$ vanishes if and only if $2u^2+4u-3=0$. Then
$u^2+2u=\frac32$, $Au^2=\frac12$, $(u+1)^2=\frac52$, and $(u+1)(u+2)=u+\frac72$, and \eqref{eq:slice-g0} becomes the last equation of
\eqref{eq:slice-sol}. Conversely, these points satisfy all six equations. A positive point has $u>0$, so $u=\frac12\sqrt{10}-1$ and
$A=1/(7-2\sqrt{10})=\frac19(7+2\sqrt{10})$. The product of the two roots in $\gamma$ is $-5A/(4(2A+1))<0$, so exactly one is positive.
Solving the quadratic gives $\gamma_{\rm ALF}$. The conjugation $\sqrt{10}\mapsto-\sqrt{10}$ gives the other value of $u$. The four values of $\gamma$
are the roots of $\mathcal E_4$, permuted transitively by the Galois group of $\F$.
\end{proof}

\subsection{Positive solutions are symmetric}\label{subsec:ALF-forcing}
The positive symmetric point of Proposition~\ref{prop:ALF-slice} is the only positive solution of
\eqref{eq:ALF-lin} (Proposition~\ref{prop:ALF-symmetric}). The proof compares a positive solution with the symmetric one through Li and Sun's
canonical harmonic maps and then uses a fold in the first components of \eqref{eq:ALF-lin}.

Let $v_0,\dots,v_4$ be the vectors $(1,0)$, $(0,1)$, $(-1,-1)$, $(-1,-2)$, $(2,3)$ of \eqref{eq:ALF-rods}, which satisfy
$v_{j+1}+v_{j-1}=\kappa_jv_j$. For $\alpha>0$ let $F_\alpha$ be the linear map of $\R^2$ with $F_\alpha(v_0)=(1,\frac\alpha2)$ and
$F_\alpha(v_4)=(1,-\frac\alpha2)$, and $\ell_j=F_\alpha(v_j)$ \eqref{eq:appC-ell}. These vectors match the outer rods of positive data. In the scaling of
\cite[(8.42)]{Gerochcf}, the vectors $\mathbf n_j=A_{\omega,\gamma}V_j$ of \eqref{eq:appC-n}, the outer rods of any positive data have
$\mathbf n_0=(1,\frac\alpha2)$ and $\mathbf n_4=(1,-\frac\alpha2)$ with $\alpha=\alpha(\omega,\gamma)$ \eqref{eq:appC-nouter}. For a positive
solution of \eqref{eq:ALF-lin} with this $\alpha$ the vectors $\mathbf n_j$ satisfy the same relations as the $v_j$, so $\mathbf n_j=\ell_j$
for all $j$. For $L=(L_1,L_2,L_3)\in\R^3_{>0}$ let $R_{L,\alpha}$ be the
rod structure of type ALF$^-$ \cite[Def.~3.1]{LS2025} with the turning points $0$, $L_1$, $L_1+L_2$, $L_1+L_2+L_3$, the rod lines
$[\ell_0],\dots,[\ell_4]$, and the model end of type ALF$^-$ with the parameter $\alpha$ \cite[(4.7)]{LS2025}. By
\cite[Thm.~4.24, Prop.~4.26]{LS2025} there is a unique harmonic map strongly tamed by $R_{L,\alpha}$ with the end-normalised augmentation,
the \emph{canonical harmonic map} of $R_{L,\alpha}$. Along each finite rod $I_j$ its scaled rod vector \cite[Def.~4.17]{LS2025} is
$\pm e^{\mathfrak b_j(L,\alpha)}\ell_j$. For the lattice $F_\alpha(\Z^2)$, in which any two consecutive vectors $\ell_{j-1},\ell_j$ form a
basis, the cone angle along $I_j$ \cite[Def.~4.19]{LS2025} is $2\pi e^{-\mathfrak b_j(L,\alpha)}$. We call
$\mathfrak b(L,\alpha)=(\mathfrak b_1,\mathfrak b_2,\mathfrak b_3)(L,\alpha)$ the \emph{cone errors} of $R_{L,\alpha}$, and the pair $(L,\alpha)$
\emph{balanced} if $\mathfrak b(L,\alpha)=0$, that is, if all cone angles along the axis equal $2\pi$. For $c>0$,
$\mathfrak b(cL,c\alpha)=\mathfrak b(L,\alpha)$ (Lemma~\ref{lem:appC-scale}).

Appendix~\ref{app:convex} adapts Aksteiner's convex potential for AF solitons \cite{Aksteiner2026}
to rod structures with a constant asymptotic parameter $\alpha\ge0$, among them the rod structures $R_{L,\alpha}$. For each $\alpha>0$ the cone errors $\mathfrak b(\cdot,\alpha)$ are the
gradient of a strictly convex real-analytic function on $\R^3_{>0}$ (Proposition~\ref{prop:appC-convex}), so for each $\alpha$ at most one $L$ is balanced.
Together with Lemma~\ref{lem:appC-scale} this shows that any two balanced pairs are positive multiples of each other (Corollary~\ref{cor:appC-ray}).
Every positive solution of \eqref{eq:ALF-lin} gives a balanced pair, formed by its rod lengths $L=(1,L_2,L_3)$ and the number
$\alpha=\alpha(\omega,\gamma)>0$ of \eqref{eq:alpha} (Proposition~\ref{prop:appC-balanced}).
The first components of \eqref{eq:ALF-lin} do not involve $\gamma$. In coordinates in which three turning points are $0$, $1$, and $\infty$,
Appendix~\ref{app:ALF} solves them for the weights. On their solutions with positive weights
the cross-ratio $t$ of the turning points is at least $\frac52$, and it equals $\frac52$ only where the ratio $\varrho$ of the two middle weights in
these coordinates, multiplied by $t$, is $1$, the value of the symmetric solutions. This is the fold of Lemma~\ref{lem:appALF-fold}.

\begin{proposition}\label{prop:ALF-symmetric}
Let $\mathbf d$ be a solution of \eqref{eq:ALF-lin} with positive coordinates \eqref{eq:ALF-unknowns}. Then $\mathbf d$ is the positive point
of Proposition~\ref{prop:ALF-slice}, with the data \eqref{eq:ALF-data}. In particular, it lies on the symmetric locus \eqref{eq:locus}.
\end{proposition}

\begin{proof}
Let $(\omega,\gamma)$ be the normalised data of $\mathbf d$, with the rod lengths $L=(1,L_2,L_3)$. Let $\mathbf d^*$ be the positive point
of Proposition~\ref{prop:ALF-slice}, with the data $(\omega^*,\gamma_{\rm ALF})$ of \eqref{eq:ALF-data} and the rod lengths $L^*=(1,u,1)$,
$u=\frac12\sqrt{10}-1$. Both are positive solutions of \eqref{eq:ALF-lin}, $\mathbf d^*$ by Proposition~\ref{prop:ALF-slice}. Hence, by
Proposition~\ref{prop:appC-balanced}, the pairs $(L,\alpha)$ and
$(L^*,\alpha^*)$, with $\alpha=\alpha(\omega,\gamma)>0$ and $\alpha^*=\alpha(\omega^*,\gamma_{\rm ALF})>0$, are balanced. By
Corollary~\ref{cor:appC-ray} they lie on one ray, $(L,\alpha)=c(L^*,\alpha^*)$ with $c>0$, and the first rod lengths $L_1=L^*_1=1$ give
$c=1$. So $L_2=u$ and $L_3=1$, and the numbers $t$ and $y$ of \eqref{eq:appALF-ty} are
\begin{equation}\label{eq:ALF-ty-star}
t=\frac{(1+L_2)(L_2+L_3)}{L_3}=(1+u)^2=\tfrac52,\qquad y=L_2+L_3=1+u,\qquad y^2=t .
\end{equation}
Let $b_2,b_3,b_\infty$ be the gauge weights of \eqref{eq:appALF-ty} and $\varrho=tb_3/b_2$. By Corollary~\ref{cor:appALF-reduced}, $\varrho>0$, the
gauge weights satisfy \eqref{eq:appALF-b}, and $\mathcal C(\varrho,t)=0$. Since $t=\frac52$, the fold, Lemma~\ref{lem:appALF-fold}, gives
$\varrho=1$. For $\varrho=1$, \eqref{eq:appALF-b} gives $b_\infty=t$ and $b_2=tb_3$, and \eqref{eq:appALF-ty} with $y^2=t$ gives
\begin{equation}\label{eq:ALF-sym-a}
a_4=\frac{b_\infty}{y^2}=1,\qquad a_2=\frac{b_2(1+y)^2}{y^2}=b_3(1+y)^2=\frac{b_3(t+y)^2}{y^2}=a_3 .
\end{equation}
Together with $L_3=1$ this says that $\mathbf d$ lies on the symmetric locus \eqref{eq:locus}, and by Proposition~\ref{prop:ALF-slice} the
only positive solution of \eqref{eq:ALF-lin} on it is $\mathbf d^*$.
\end{proof}

The symmetry of the system alone does not give Proposition~\ref{prop:ALF-symmetric}, since the solutions off the symmetric locus come in pairs
$\{\mathbf d,\sigma(\mathbf d)\}$. The order of the turning points also matters. An exact computation of all complex solutions of \eqref{eq:ALF-lin} gives a real pair with $\gamma>0$, positive weights, and
$L_3>0$ whose only negative coordinate is the middle rod length $L_2$. In the proof the order of the
turning points enters through the rod structures $R_{L,\alpha}$, which require $L\in\R^3_{>0}$, and through the conditions $t>1$ and
$y>0$ of the gauge of Appendix~\ref{app:ALF}.

\subsection{Proof of Theorem~\ref{thm:main}(c)}\label{subsec:ALF-proof}
We write $g_{\rm ALF}$ for the metric \eqref{eq:closed} of the data \eqref{eq:ALF-data}, with $r=2$ and $\mathfrak c=|D_0|^{-1}$, on the
torus bundle whose lattice is generated by $\mathfrak c^{1/2}V_0$ and $\mathfrak c^{1/2}V_1$. By the converse below it is the instanton of
Theorem~\ref{thm:main}(c). Let $(M,g)$ be a simply connected toric gravitational instanton whose torus
action has four fixed points and whose end is ALF, and suppose that $g$ is hyperk\"ahler for neither orientation. By
Proposition~\ref{prop:reduction}(iii) and Theorem~\ref{thm:soliton}, after a change of orientation $\kappa=(-1,1,-1)$, $r=2$, and $g$ has
positive soliton data satisfying \eqref{eq:ALF-lin}. By Theorem~\ref{thm:soliton}, $g$ is \eqref{eq:closed} with these data and some
constant $\mathfrak c>0$, on the lattice generated by $v_0=\mathfrak c^{1/2}V_0$ and $v_1=\mathfrak c^{1/2}V_1$. Changing $\mathfrak c$ only rescales the metric (after
Theorem~\ref{thm:converse}), so up to scale we may assume $\mathfrak c|D_0|=1$. The normalisation of these soliton data is a positive solution, which is the point \eqref{eq:ALF-data}
by Proposition~\ref{prop:ALF-symmetric}. By Lemma~\ref{lem:normalise} the
normalisation changes the metric \eqref{eq:closed} by a constant factor and a change of the angles, so $g$ agrees with $g_{\rm ALF}$ up to scale.

Conversely, the data \eqref{eq:ALF-data} are positive and satisfy \eqref{eq:ALF-lin} by Proposition~\ref{prop:ALF-slice}, so by
Theorem~\ref{thm:converse} they define a smooth complete simply connected toric gravitational instanton. Its chain is $(-1,1,-1)$, which gives
$p=\det(v_0,v_4)=3\det(v_0,v_1)\ne0$ and $\alpha(\omega,\gamma_{\rm ALF})>0$ (\S\ref{subsec:ALF-system}), so its end is ALF$^-$.
Since the chain is not $\pm(2,2,2)$, the instanton is
hyperk\"ahler for neither orientation by Proposition~\ref{prop:reduction}(i). Both its Li--Sun degrees are therefore $2$ by
Proposition~\ref{prop:reduction}(ii), and it is Hermitian for neither orientation by Theorem~\ref{thm:LSdegree}.
The description of the end is proved in Section~\ref{sec:explicit}. \qed

The asymptotic parameter of the data \eqref{eq:ALF-data} is, by \eqref{eq:alpha},
\begin{equation}\label{eq:ALF-alpha}
\alpha(\omega,\gamma_{\rm ALF})=\tfrac1{123}\bigl(1300+840\sqrt2+592\sqrt5+380\sqrt{10}\bigr),
\end{equation}
with $H_1(\omega)=\frac19(32+4\sqrt{10})$, $H_2(\omega)=\frac13(28+8\sqrt{10})$, and
$H_3(\omega)=\frac13(80+25\sqrt{10})$.

\subsection{The other orientation}\label{subsec:ALF-orient}
A duality of Hankel determinants gives the soliton data of $g_{\rm ALF}$ with the opposite orientation of $M$, in which the chain is
$(1,-1,1)$. For a measure $\omega=\sum_{k=1}^4a_k\delta_{w_k}$ on four distinct points with non-zero weights and
$\gamma\ne0$ put
\begin{equation}\label{eq:duality}
\begin{gathered}
(\omega,\gamma)\longmapsto(\omega^\vee,\gamma^{-1}),\qquad \omega^\vee=\sum_{k=1}^4\frac1{a_k\,q'(w_k)^2}\,\delta_{w_k},\qquad
q(T)=\prod_{k=1}^4(T-w_k),\\
c_\omega=a_1a_2a_3a_4\prod_{k<l}(w_k-w_l)^2,\qquad \mathfrak s=w_1+w_2+w_3+w_4 .
\end{gathered}
\end{equation}

\begin{lemma}\label{lem:duality}
Let $\omega$, $\gamma$, $\omega^\vee$, $c_\omega$, and $\mathfrak s$ be as in \eqref{eq:duality}. Then the Hankel families \eqref{eq:V} with $r=2$ satisfy
\begin{equation}\label{eq:dualV}
V_j(\omega^\vee,\gamma^{-1})=(-1)^j\,\Theta\,V_j(\omega,\gamma),\qquad
\Theta=\frac1{c_\omega\gamma^2}\begin{pmatrix}\gamma^2&0\\-2\gamma \mathfrak s&-1\end{pmatrix},\qquad 0\le j\le4 .
\end{equation}
In particular, $(\omega,\gamma)$ satisfies $V_{j+1}+V_{j-1}=\kappa_jV_j$ for $j=1,2,3$ if and only if $(\omega^\vee,\gamma^{-1})$ satisfies
these relations with $-\kappa$ in place of $\kappa$.
\end{lemma}

\begin{proof}
Let $S\subset\{1,2,3,4\}$ with complement $S^c$. Since $q'(w_k)^2=\prod_{l\ne k}(w_k-w_l)^2$,
\begin{equation}\label{eq:dualQ}
\prod_{k\in S}q'(w_k)^2=\prod_{k<l\in S}(w_k-w_l)^4\prod_{k\in S,\,l\in S^c}(w_k-w_l)^2 ,
\end{equation}
so by \eqref{eq:PS} $c_\omega P_S(\omega^\vee)=P_{S^c}(\omega)$: the weights $a_k$ with $k\in S$ cancel, and of the squared differences only
those within $S^c$ remain. The measures $\omega^\vee_{(j)}$ and $\omega_{(j)}$ reverse the signs of the same atoms, and the signs
$(-1)^{|S\cap\{1,\dots,j\}|}$ and $(-1)^{|S^c\cap\{1,\dots,j\}|}$ have the product $(-1)^j$. Hence
$c_\omega P_S(\omega^\vee_{(j)})=(-1)^jP_{S^c}(\omega_{(j)})$. Summing over the sets with $s$ atoms and multiplying, for two atoms, by
$w_k+w_l=\mathfrak s-\sum_{m\in S^c}w_m$ before summing, Heine's formula \eqref{eq:heine} and \eqref{eq:H3gamma} give
\begin{equation}\label{eq:dualH}
c_\omega H_s(\omega^\vee_{(j)})=(-1)^jH_{4-s}(\omega_{(j)})\quad(0\le s\le4),\qquad
c_\omega\Sigma_2(\omega^\vee_{(j)})=(-1)^j\bigl(\mathfrak sH_2(\omega_{(j)})-\Sigma_2(\omega_{(j)})\bigr).
\end{equation}
With $H_s=H_s(\omega_{(j)})$ and $\Sigma_2=\Sigma_2(\omega_{(j)})$, the first component of $V_j(\omega^\vee,\gamma^{-1})$ is therefore
$(-1)^jH_2/c_\omega$. By \eqref{eq:H3gamma} the second is
\begin{equation}\label{eq:dualproof}
\frac{(-1)^j}{c_\omega}\Bigl(H_1-\frac2\gamma\bigl(\mathfrak sH_2-\Sigma_2\bigr)-\frac{H_3}{\gamma^2}\Bigr)
=\frac{(-1)^j}{c_\omega\gamma^2}\Bigl(-\bigl(H_3-2\gamma\Sigma_2-\gamma^2H_1\bigr)-2\gamma \mathfrak sH_2\Bigr),
\end{equation}
which is \eqref{eq:dualV}. Finally, with $V^\vee_j=(-1)^j\Theta V_j$ we have
$V^\vee_{j+1}+V^\vee_{j-1}+\kappa_jV^\vee_j=(-1)^{j+1}\Theta\,(V_{j+1}+V_{j-1}-\kappa_jV_j)$, and $\Theta$ is invertible.
\end{proof}

\begin{corollary}\label{cor:ALF-other}
Let $g$ be the instanton $g_{\rm ALF}$ of Theorem~\ref{thm:main}\textup{(c)} with the opposite orientation of $M$. Then its chain is
$(1,-1,1)$, and its normalised soliton data are the turning points \eqref{eq:ALF-data}, the weights $(1,9A,9A,1)$, and the atom weight
$45A/(4\gamma_{\rm ALF})$.
\end{corollary}

\begin{proof}
The chain is $(1,-1,1)$ by \S\ref{subsec:rods}, and the degree index is again $r=2$, since both Li--Sun degrees are $2$
(Theorem~\ref{thm:soliton}(i)). The duality \eqref{eq:duality} maps positive data to positive data and is an involution. It is compatible
with the changes of Lemma~\ref{lem:normalise}. Let $\omega'=\sum_kKa_k\delta_{\lambda w_k+\mu}$ and $\gamma'=K\lambda^3\gamma$. The polynomial of
the points $\lambda w_k+\mu$ has the derivative $\lambda^3q'(w_k)$ at $\lambda w_k+\mu$, so $(\omega'^\vee,\gamma'^{-1})$ is obtained from
$(\omega^\vee,\gamma^{-1})$ by the change with the same $\lambda,\mu$ and with $K^{-1}\lambda^{-6}$ in place of $K$. Let $(\omega,\gamma)$ be
soliton data of $g$ (Theorem~\ref{thm:soliton}). They are positive and satisfy the relations for $(1,-1,1)$. By
Lemma~\ref{lem:duality}, $(\omega^\vee,\gamma^{-1})$ is a positive solution of \eqref{eq:ALF-lin}, so by
Proposition~\ref{prop:ALF-symmetric} its normalisation is \eqref{eq:ALF-data}. Hence the normalisation of
$(\omega,\gamma)$ is the normalisation of the dual of \eqref{eq:ALF-data}. For \eqref{eq:ALF-data}, with $u=\frac12\sqrt{10}-1$ as in
Proposition~\ref{prop:ALF-slice}, $q'(w_1)^2=q'(w_4)^2=(u+1)^2(u+2)^2=\frac{45}4A$ and $q'(w_2)^2=q'(w_3)^2=u^2(u+1)^2$, by the symmetry of
the turning points. Since $Au^2=\frac12$ by \eqref{eq:slice-sol}, the dual weights multiplied by $q'(w_1)^2$, which makes the first
one equal to $1$, are $1$, $(u+2)^2/(Au^2)=2(u+2)^2=9A$, $9A$, and $1$. The same factor, the change with $\lambda=1$, $\mu=0$, and
$K=q'(w_1)^2$, takes $\gamma_{\rm ALF}^{-1}$ to $45A/(4\gamma_{\rm ALF})$.
\end{proof}

\section{The AF instantons: proof of Theorem~\ref{thm:main}(d)}\label{sec:AF}
By the convexity of Appendix~\ref{app:convex} (Proposition~\ref{prop:AF-main}), an AF instanton with four turning points is
determined up to scale by its twist, once the smoothness conditions are written down (\S\ref{subsec:AF-sym}) and the twist is
identified with the number $-x_1$ of \eqref{eq:x} (\S\ref{subsec:AF-twist}). The existence for every twist in $(0,1)$ is Li and
Sun's theorem \cite{LS2025}. Sections~\ref{subsec:AF-teo} and~\ref{subsec:AF-family} relate these instantons to Teo's
family.

\subsection{Teo's parameters and Teo's locus}\label{subsec:AF-sym}
By
Proposition~\ref{prop:reduction}(iv) the chain lies in the class of $(1,0,-1)$, the chain of $X_3$. We fix the
orientation of $\mathfrak t$ in which, with \eqref{eq:det}, the chain is $(-1,0,1)$, so that $\varepsilon=-1$ in
Proposition~\ref{prop:reduction}(iv). The smoothness conditions are
\begin{equation}\label{eq:AF-lin}
V_2=-V_1-V_0,\qquad V_3=-V_1,\qquad V_4=V_0 .
\end{equation}

Teo constructs his metrics \cite[before (1)]{Teo2026} as four-soliton solutions on flat space in the sense of Chen
\cite[\S3.1]{Chen2016}. They have poles at four turning points $z_1<\dots<z_4$ and Belinski--Zakharov parameters
$C_1,\dots,C_4\in\R\cup\{\infty\}$. We normalise the turning points to $w_1=0$, $w_2=1$, $1<w_3<w_4$, and let
$q(T)=\prod_{l=1}^4(T-w_l)$. For positive soliton data with $r=2$ on these points (Theorem~\ref{thm:soliton}), the numbers
$c_k=-a_kq'(w_k)$ together with $\gamma$ are the data of a far-node configuration with four finite nodes in the sense of
\cite[\S8]{Gerochcf}. The Gram matrix of the Killing fields in \eqref{eq:closed} is then, up to a constant factor and a constant change
of frame, that of the Belinski--Zakharov soliton on $\R^3\times S^1$ with the poles $w_k$ and the parameters
$b_k=(c_k+\gamma)/(c_k-\gamma)$ \cite[Thms.~7.3(b), 8.3, 8.5]{Gerochcf}. The computation below, which also treats the conformal factor,
shows that Teo's parameters are $C_k=1/b_k$, that is,
\begin{equation}\label{eq:dictionary}
C_k=\frac{c_k-\gamma}{c_k+\gamma},\qquad c_k=-a_k\,q'(w_k),\qquad q(T)=\prod_{l=1}^4(T-w_l).
\end{equation}
Since $q'(w_k)$ has the signs $-,+,-,+$, the weights are positive if and only if
\begin{equation}\label{eq:AF-positive}
|C_1|<1<|C_2|,\qquad |C_3|<1<|C_4| ,
\end{equation}
where $C_2$ and $C_4$ may be infinite ($c_k=-\gamma$) and $C_1$ and $C_3$ may vanish ($c_k=\gamma$). Conversely,
\eqref{eq:dictionary} maps $|c_k|/\gamma\in(0,\infty)$ bijectively onto the values of $C_k$ allowed by \eqref{eq:AF-positive}. So every
$\mathbf C=(C_1,\dots,C_4)$ with \eqref{eq:AF-positive} comes from positive data on the given turning points, unique up to a common
positive factor of the weights and $\gamma$.

Teo's rod vectors $\ell^{\rm T}_1,\dots,\ell^{\rm T}_5$ \cite[(4a)--(4e)]{Teo2026} are explicit rational functions of his parameters and
the turning points. Where $C_2$ and $C_4$ are finite, a direct computation with Heine's formula \eqref{eq:heine}, which expands both
sides in the Heine terms of the data, shows that
\begin{equation}\label{eq:AF-Teorods}
\ell^{\rm T}_{j+1}=s_j\,\mathsf M\,V_j,\qquad 0\le j\le4,\qquad (s_0,\dots,s_4)=(1,1,-1,-1,1),
\end{equation}
with an invertible matrix $\mathsf M$ that depends on the data but not on $j$. That the two sides span the same lines already
follows from \cite[Thm.~7.3(b), Prop.~8.7(i), Thm.~8.5]{Gerochcf}, which describe the rod lines of the soliton through the Hankel
family, and the computation determines the signs $s_j$. By \eqref{eq:AF-Teorods}, Teo's conditions $\ell^{\rm T}_1=\ell^{\rm T}_5$
and $\ell^{\rm T}_2=\ell^{\rm T}_4$ \cite[before (5)]{Teo2026} become $V_4=V_0$ and $V_3=-V_1$, and his condition
$\ell^{\rm T}_1+\ell^{\rm T}_2=\ell^{\rm T}_3$ \cite[(6)]{Teo2026} becomes $V_2=-V_1-V_0$. Thus Teo's smoothness conditions are
\eqref{eq:AF-lin}.

Teo's metric \cite[(1)--(3)]{Teo2026} is, up to a constant factor, the closed formula in other coordinates. For the Gram matrix of
the Killing fields, Chen's dressing formula \cite[(2.16)--(2.20), (3.1)--(3.6)]{Chen2016} is the dressing formula of
\cite[(4.3), (4.13)]{Gerochcf} with the parameters $b_k=1/C_k$. Chen's conformal factor is an explicit function of $\rho$ and the
$\mu_k$ times the determinant of his matrix $\Gamma$, and it is computed in \cite[Lemmas~8.2, 8.4(i) and the proof of
Thm.~8.5]{Gerochcf}. Together with that theorem this shows the following. For positive data with finite parameters
\eqref{eq:dictionary}, Teo's metric with the turning points $z_k=w_k$, these parameters, and any value of his constant $k_0>0$ is a
constant multiple of \eqref{eq:closed}, in angles that are a constant linear function of Teo's angles. A translation of the turning
points does not change Teo's metric, since $z$ and the $z_k$ enter it only through Chen's pole functions $\mu_k$ \cite[(2)]{Teo2026}.
The same computations show that Teo's rod vectors are, up to sign, the scaled rod vectors of his metric. Since a constant factor does
not change scaled rod vectors, the linear map of the angles takes the lattice of Teo's identifications \cite[(9)]{Teo2026}, generated by
$\ell^{\rm T}_1$ and $\ell^{\rm T}_2$, to the lattice of Theorem~\ref{thm:converse}, generated by $\mathfrak c^{1/2}V_0$ and
$\mathfrak c^{1/2}V_1$. Consequently, for positive data with finite parameters that satisfy \eqref{eq:AF-lin}, Teo's metric with his
identifications is, up to scale, the instanton of Theorem~\ref{thm:converse}. In particular it is smooth and complete, and it is
positive definite, which Teo had checked numerically \cite[after (11)]{Teo2026}. For the ALF instanton of Theorem~\ref{thm:main}(c),
which is also of Teo's form, Jones proves these properties directly \cite{Jones2026} (\S\ref{subsec:main}).

Teo solves his first two conditions by a symmetric ansatz \cite[(5a), (5b)]{Teo2026}. In our normalisation it is \emph{Teo's locus}
\begin{equation}\label{eq:Teo-locus}
w_4-w_3=w_2-w_1=1,\qquad C_3=-1/C_2,\qquad C_4=-1/C_1 .
\end{equation}
On Teo's locus the relations $V_4=V_0$ and $V_3=-V_1$ hold identically. Where $C_2$ and $C_4$ are finite, one checks this by inserting \eqref{eq:Teo-locus} into Teo's
rod vectors \cite[(4a)--(4e)]{Teo2026} and using \eqref{eq:AF-Teorods}. Elsewhere it follows by the continuity of the
$V_j$ in the weights.

\subsection{The twist}\label{subsec:AF-twist}
The twist of the end is the number $\vartheta=-x_1$ (Lemma~\ref{lem:twist}). Let
$\xi$ be the Killing field that stays bounded at infinity, scaled to the period of the model, so that $\Lambda$ is generated by $v_0$ and
$\beta v_0+\xi$ \cite[before Cor.~11.2]{Gerochcf}.

\begin{lemma}\label{lem:twist}
Let $(M,g)$ be a simply connected toric gravitational instanton whose torus action has four fixed points, with an AF end and the chain
$(-1,0,1)$ in the orientation of $\mathfrak t$ of \S\ref{subsec:AF-sym}. Let $\beta$ be its twist, $x_1$ as in \eqref{eq:x}, and
$\vartheta=-x_1$. Then
\begin{equation}\label{eq:twist}
\beta\equiv\vartheta\pmod\Z ,
\end{equation}
that is, $\beta=\vartheta$.
\end{lemma}

\begin{proof}
By the definition of the twist (\S\ref{subsec:rods}), $\Lambda$ is generated by $v_0$ and $\beta v_0+\xi$, where $\xi$ is taken
with $\det(v_0,\xi)=1$. Changing the common sign of the rod vectors changes the signs of $v_0$ and $\xi$ and not $\beta$. Write
$v_1=p_1v_0+q_1(\beta v_0+\xi)$ with $p_1,q_1\in\Z$. By \eqref{eq:det}, $-1=\det(v_0,v_1)=q_1\det(v_0,\xi)=q_1$. In the basis $(f_0,f_1)$ of Theorem~\ref{thm:soliton}, the Killing field $f_1$ has the length $|\mathcal X|(F_1F_2)^{-1/2}$ by
\eqref{eq:closed}. At infinity all $\varpi_k$ tend to the same unimodular number, so that
$|\mathcal X|\to H_2(\omega)$ and $F_1F_2\to\gamma H_2(\omega)^2$. Hence the length of $f_1$ tends to $\gamma^{-1/2}$ at infinity. The Killing field $f_0$ has the length $|\mathcal Y|(F_1F_2)^{-1/2}$, which grows like $\rho$
\cite[Prop.~8.8 and the proof of Cor.~11.2]{Gerochcf}. Here $r=2$ by Proposition~\ref{prop:reduction}(ii). So the Killing fields that stay bounded at infinity are the multiples of $f_1$, and
$\xi\in\R f_1$. By Theorem~\ref{thm:soliton}(iii) and \eqref{eq:x}, $x_j=\mathsf x(v_j)$ for the linear form $\mathsf x$ on $\mathfrak t$ that vanishes on $f_1$ and
takes the value $1$ on $v_0$. Hence $x_1=\mathsf x(v_1)=p_1+q_1(\beta+\mathsf x(\xi))=p_1-\beta$, which is \eqref{eq:twist}. Since $\beta\in[0,1)$ and
$\vartheta\in(0,1)$ by Proposition~\ref{prop:reduction}(iv), it follows that $\beta=\vartheta$.
\end{proof}

The proof uses only that the metric is \eqref{eq:closed} in the positively oriented basis $(f_0,f_1)$, that $v_j=\mathfrak c^{1/2}V_j$
with $\det(v_0,v_1)=-1$, that $r=2$, and that $0<\vartheta<1$. So it applies, with $x_1$ computed from the data and in the orientation of
$\mathfrak t$ for which $(f_0,f_1)$ is positively oriented, to the instanton that
Theorem~\ref{thm:converse} gives from positive data with $r=2$ that satisfy \eqref{eq:AF-lin} and $0<-x_1<1$. In this orientation $\det(v_0,v_1)=\mathfrak cD_0=-1$, since $D_0<0$. This is how
Lemma~\ref{lem:twist} is used in the proofs of Proposition~\ref{prop:AF-main} and of Theorem~\ref{thm:main}(d).

\subsection{Convexity and the twist}\label{subsec:AF-main}
The twist fixes rod lines to which the convexity of Appendix~\ref{app:convex} applies, at the asymptotic parameter $\alpha=0$, and
this gives the uniqueness of the instanton with a given twist (Proposition~\ref{prop:AF-main}).

For $\vartheta\in(0,1)$ let
\begin{equation}\label{eq:AF-ell}
\ell_0(\vartheta)=(1,0),\quad \ell_1(\vartheta)=(-\vartheta,-1),\quad \ell_2(\vartheta)=(\vartheta-1,1),\quad \ell_3(\vartheta)=(\vartheta,1),
\quad \ell_4(\vartheta)=(1,0),
\end{equation}
and let $\Lambda_\vartheta=\Z\ell_0(\vartheta)+\Z\ell_1(\vartheta)$. These vectors satisfy $\ell_{j+1}+\ell_{j-1}=\kappa_j\ell_j$ for the
chain $\kappa=(-1,0,1)$ and $\det(\ell_{j-1},\ell_j)=-1$ for $1\le j\le4$, so every pair $(\ell_{j-1},\ell_j)$ is a basis of
$\Lambda_\vartheta$. Their first components $-\vartheta$ and $\vartheta-1$ for $j=1,2$ are negative, so they satisfy \eqref{eq:appC-adm}
with $\alpha=0$. For $L\in\R^3_{>0}$ let $R_{L,\vartheta}$ be the rod structure $R_L$ of \S\ref{subsec:appC-rods} formed with these
vectors. It has the turning points $w(L)$ of \eqref{eq:appC-wL}, the rod lines $[\ell_j(\vartheta)]$, and an end of type AF$_0$, and the
twist is recorded by the lattice. Let $\mathfrak b(L,\vartheta)\in\R^3$ be its cone errors \eqref{eq:appC-b}, and call $(L,\vartheta)$
\emph{balanced} if $\mathfrak b(L,\vartheta)=0$. For the lattice $\Lambda_\vartheta$ the cone angle along the rod $I_j$ is
$2\pi e^{-\mathfrak b_j(L,\vartheta)}$ \cite[Def.~4.19]{LS2025}. By Proposition~\ref{prop:appC-convex}, for fixed $\vartheta$ the map
$\mathfrak b(\cdot,\vartheta)$ is the gradient of a strictly convex function on $\R^3_{>0}$, so at most one $L$ is balanced. The
balanced $L$, if it exists, is therefore unique, not only unique up to scale. This is consistent with the scaling of the turning points
because the twist fixes the scale. Indeed, by the third identity of \eqref{eq:appC-nscale}, a dilation of the turning points by $c\ne1$
replaces the lines $[\ell_j(\vartheta)]$ by $[\mathrm{diag}(1,c)\,\ell_j(\vartheta)]$. For $j=1,2$, written with the second components
$-1$ and $1$, the vectors spanning these lines have first components adding up to $-1/c$, whereas those of $\ell_1(\vartheta')$ and
$\ell_2(\vartheta')$ add up to $-1$ for every $\vartheta'$.

\begin{proposition}\label{prop:AF-main}
Let $\vartheta\in(0,1)$. Up to scale there is exactly one simply connected toric gravitational instanton whose torus action has four
fixed points, with an AF end, the chain $(-1,0,1)$ in the orientation of $\mathfrak t$ of \S\ref{subsec:AF-sym}, and the twist
$\vartheta$. Up to scale its rod structure is $R_{L,\vartheta}$ with the lattice $\Lambda_\vartheta$, for the balanced $L$.
\end{proposition}

\begin{proof}
Let $(M,g)$ be an instanton as in the proposition. By Proposition~\ref{prop:reduction}(i) it is hyperk\"ahler for neither orientation.
So by Theorem~\ref{thm:soliton} and Proposition~\ref{prop:reduction}(ii),(iv) it has positive soliton data $(\omega,\gamma)$ with $r=2$
that satisfy \eqref{eq:AF-lin}. By Lemma~\ref{lem:normalise} we may take them on points $w=(0,1,w_3,w_4)$, and we put
$L=(1,w_3-1,w_4-w_3)$. By \eqref{eq:appC-A}, the vectors $\mathbf n_j=\mathbf n_j(\omega,\gamma)$ of \eqref{eq:appC-n} are the images of
the $V_j$ under one invertible matrix, so they satisfy $\mathbf n_{j+1}+\mathbf n_{j-1}=\kappa_j\mathbf n_j$ for $j=1,2,3$. In particular
$\mathbf n_4=\mathbf n_0$, and \eqref{eq:appC-nouter} gives $\alpha(\omega,\gamma)=0$ and $\mathbf n_0=(1,0)$. The first component of
$\mathbf n_1$ is the number $x_1$ of \eqref{eq:x}, which is $-\vartheta$ by Lemma~\ref{lem:twist}. Its second component
$y_1=\det(\mathbf n_0,\mathbf n_1)$ is negative (\S\ref{subsec:appC-rods}). With $c=-1/y_1$, the matrix $\mathrm{diag}(1,c)$ maps
$\mathbf n_0$ and $\mathbf n_1$ to $\ell_0(\vartheta)$ and $\ell_1(\vartheta)$. Both families satisfy the same recursion, so
$\mathrm{diag}(1,c)\,\mathbf n_j=\ell_j(\vartheta)$ for $0\le j\le4$. By the third identity of \eqref{eq:appC-nscale}, the dilated data
$\tilde\omega=\sum_ka_k\delta_{cw_k}$ and $\tilde\gamma=c^3\gamma$ on the points $w(cL)$ have $\mathbf n_j(\tilde\omega,\tilde\gamma)
=\ell_j(\vartheta)$ for $0\le j\le4$, and $\alpha(\tilde\omega,\tilde\gamma)=0$. So they satisfy the hypotheses of
Proposition~\ref{prop:appC-ident}(ii) for $\tilde L=cL$ and the vectors $\ell_j(\vartheta)$. Then \eqref{eq:appC-nb} gives
$\mathfrak b(\tilde L,\vartheta)=0$. Hence $\tilde L$ is the balanced $L$ for $\vartheta$.

By Lemma~\ref{lem:normalise} with $\lambda=c$, $\mu=0$, and $K=1$, the metric $g$ agrees up to scale with the metric \eqref{eq:closed} of
$(\tilde\omega,\tilde\gamma)$, with some constant $\mathfrak c>0$. The basis $(a_{\tilde\gamma},b_{\tilde\gamma})$ of
\cite[Prop.~8.9]{Gerochcf} differs from the basis of \eqref{eq:closed} by a change of determinant $1$ (after
Theorem~\ref{thm:converse}). Read in this basis, the metric is the toric metric of the far-node member with the data $(\tilde\omega,\tilde\gamma)$ and the
constant $\mathfrak c$ (\S\ref{subsec:appC-rods}). By Proposition~\ref{prop:appC-ident}(ii), its harmonic map is the canonical harmonic
map $\Phi_{\tilde L}$ of $R_{\tilde L,\vartheta}$. By the same proposition, its augmentation is $\nu_{\tilde L}+\log\varsigma$, where $\nu_{\tilde L}$ is the
end-normalised augmentation and $\varsigma=(\mathfrak c\tilde\gamma H_2(\tilde\omega)^2)^{1/2}$. By \eqref{eq:appC-farnorm} and
\eqref{eq:appC-nb}, its scaled rod vectors are $\pm\varsigma\,\ell_j(\vartheta)$. Since $g$ is smooth, these are its rod vectors
\cite[Def.~4.19]{LS2025}. Consecutive rod vectors generate the lattice \cite{OR1970}, which is therefore $\varsigma\Lambda_\vartheta$.
In the angles divided by $\varsigma$, in which the lattice is $\Lambda_\vartheta$, the metric is $\varsigma^2$ times the toric metric
$\rho\,\Phi_{\tilde L}(d\phi,d\phi)+e^{2\nu_{\tilde L}}(d\rho^2+dz^2)$ \cite[(2.1)]{Gerochcf}. So, up to scale, the rod structure of $g$
is $R_{\tilde L,\vartheta}$ with the lattice $\Lambda_\vartheta$. On the union of the principal orbits, $g$ is the toric metric of
$\Phi_{\tilde L}$ and $\nu_{\tilde L}$. Both depend only on $\vartheta$. A second instanton $g'$ as in the proposition therefore has the
same rod structure as $g$ up to scale. By the uniqueness theorem of Kunduri and Lucietti \cite[Thm.~1.1]{KL2021}, a toric AF instanton
is determined by its rod structure. Its hypotheses hold for $g$ and $g'$. The action has a fixed point and, since consecutive rod
vectors generate the lattice, no discrete isotropy groups \cite[Def.~2.1]{KL2021}. The decay conditions of \cite[Def.~2.2]{KL2021} hold
since both metrics are given by \eqref{eq:closed}. The metric functions of this formula are real-analytic in $(1/R,\theta)$ near infinity, where
$(\rho,z)=(R\sin\theta,R\cos\theta)$ \cite[Prop.~8.8]{Gerochcf}. So $g'$ agrees with $g$ up to scale.

The existence comes from Li and Sun's work on $X_3$. Consider their rod structure $\mathcal R_{\mathbf l}$ on $X_3$ with the
lattice $\Z^2$, the rod lengths $\mathbf l$, and the type AF$_\beta$. In their lattice basis
$(\partial_{\phi_1},\beta\partial_{\phi_1}+\partial_{\phi_2})$ of \cite[(3.8)]{LS2025} (\S\ref{subsec:rods}), the rod lines are
$[(1,0)],[(0,1)],[(-1,1)],[(0,1)],[(1,0)]$ \cite[\S6.2]{LS2025}. With the signs of \eqref{eq:det}, the rod vectors are
$(1,0),(0,-1),(-1,1),(0,1),(1,0)$, and the chain is $(-1,0,1)$. So Li and Sun's orientation of $\mathfrak t$ is that of
\S\ref{subsec:AF-sym}. In the basis $(\partial_{\phi_1},\partial_{\phi_2})$ of the model, in which the end is AF$_0$, the
coordinates $(a,b)$ in the lattice basis become $(a+\beta b,b)$. So these rod vectors become the vectors $\ell_j(\beta)$ of \eqref{eq:AF-ell}, the lattice $\Z^2$
becomes $\Lambda_\beta$, and $\mathcal R_{\mathbf l}$ is $R_{\mathbf l,\beta}$ with the lattice $\Lambda_\beta$. For every $\beta\in(0,1)$ they find a balanced $\mathbf l$ \cite[Thm.~6.4]{LS2025} and a smooth complete toric
gravitational instanton on $X_3$ with this rod structure and the twist $\beta$ \cite[Thm.~1.1]{LS2025}. Its torus action has four fixed
points, and $X_3$ is simply connected. With $\beta=\vartheta$ this is an instanton as in the proposition.
\end{proof}

Li and Sun's existence theorem rests on the behaviour of the canonical harmonic maps at the boundary of the space of rod lengths. The
uniqueness argument above, by contrast, works at fixed rod lines and needs neither this behaviour nor the existence of balanced rod
lengths. The existence for every twist also follows from Aksteiner's existence theorem \cite{Aksteiner2026}, whose first family contains
the rod structure of $X_3$, once its hypotheses on the collisions of turning points are removed, as is done in \cite{ADNprep}. The
solutions on Teo's locus in \S\ref{subsec:AF-teo} give a third proof of the existence.

\subsection{Proof of Theorem~\ref{thm:main}(d)}\label{subsec:AF-proof}
Part~(d) of Theorem~\ref{thm:main} follows from Proposition~\ref{prop:AF-main}, and with it the proof of Theorem~\ref{thm:main} is complete.

\begin{proof}[Proof of Theorem~\ref{thm:main}(d)]
Let $(M,g)$ be a simply connected toric gravitational instanton whose torus action has four fixed points and whose end is AF. It is not
hyperk\"ahler for either orientation by Proposition~\ref{prop:reduction}(i), so both its Li--Sun degrees are $2$ by
Proposition~\ref{prop:reduction}(ii), and it is Hermitian for neither orientation by Theorem~\ref{thm:LSdegree}. Its chain lies in the
class of $(1,0,-1)$ by Proposition~\ref{prop:reduction}(iv). In the orientation of $\mathfrak t$ of \S\ref{subsec:AF-sym} it is
$(-1,0,1)$, and its twist is $\beta=-x_1\in(0,1)$ by Lemma~\ref{lem:twist} and Proposition~\ref{prop:reduction}(iv). By
Proposition~\ref{prop:AF-main}, $g$ is, up to scale, the only such instanton with the twist $\beta$, and every $\beta\in(0,1)$ is the
twist of such an instanton. The manifold is $X_3$ by Section~\ref{sec:explicit} (Table~\ref{tab:topology}). Reversing the orientation of
$\mathfrak t$ gives the chain $(1,0,-1)$ and replaces $\beta$ by $-\beta\equiv1-\beta$.
\end{proof}

Theorem~\ref{thm:main}(d) adds to Li and Sun's existence theorem on $X_3$ \cite[Thm.~1.1]{LS2025} that the instanton with the chain
$(-1,0,1)$ and a given twist is unique up to scale. By the paragraph before Proposition~\ref{prop:AF-main}, its balanced rod lengths
are unique for each twist, not only up to scale. The cone angles of $R_{L,\vartheta}$ with the lattice $\Lambda_\vartheta$ are all
$2\pi$ if and only if $L$ is balanced, and by the proof of Proposition~\ref{prop:AF-main}, Li and Sun's rod structure
$\mathcal R_{\mathbf l}$ on $X_3$ with the twist $\beta$ is $R_{\mathbf l,\beta}$ with this lattice. This proves Li and Sun's
Conjecture~8.5 \cite{LS2025} for $X_3$. For all $X_m$ the uniqueness of the balanced rod lengths is proved in
\cite{ADNprep}.

That no instanton with four turning points has the twist $\beta=0$, which follows from \cite[Cor.~11.2]{Gerochcf}
(\S\ref{subsec:main}), can also be seen directly from the strip inequality of
Corollary~\ref{cor:strip}(b), as in the proof of that corollary. Indeed, Proposition~\ref{prop:reduction}(iv) gives
$0<\varepsilon x_1<1$, so $x_1\notin\Z$, while $x_1\equiv-\beta$ modulo $\Z$ by Lemma~\ref{lem:twist} and its proof.

\begin{proof}[Proof of Theorem~\ref{thm:main}]
Let $(M,g)$ be as in the theorem. If $g$ is hyperk\"ahler for one orientation,
Proposition~\ref{prop:reduction}(i) gives (a) or (b). Otherwise the end is not ALE by Proposition~\ref{prop:reduction}(ii), and
\S\S\ref{subsec:ALF-proof} and~\ref{subsec:AF-proof} give (c) and (d). The converse statements are the classical ones for (a) and (b)
\cite{GH1978}, the converse part of \S\ref{subsec:ALF-proof} for (c), and for (d) the instantons of Proposition~\ref{prop:AF-main}. By
Theorem~\ref{thm:soliton}, each of them is given by \eqref{eq:closed} with soliton data.
\end{proof}

\subsection{Teo's metrics}\label{subsec:AF-teo}
The instantons of Theorem~\ref{thm:main}(d) are, up to scale, the members of Teo's family, and the twists of the members fill the
interval $(0,1)$. This rests on the relation between Teo's metric and the closed formula in \S\ref{subsec:AF-sym} and on an
intermediate value argument in a chart of Aksteiner.

On Teo's locus \eqref{eq:Teo-locus} we write the turning points as $w=(0,1,1+L_2,2+L_2)$, with $L_2>0$, and $\mathbf C=(C_1,C_2,-1/C_2,-1/C_1)$. There $V_4=V_0$
and $V_3=-V_1$ hold identically. By \eqref{eq:AF-positive}, the positive data are the points with $|C_1|<1<|C_2|$. For finite
$C_2$, inserting the locus into Teo's rod vectors and using \eqref{eq:AF-Teorods} shows that the remaining relation $V_2=-V_1-V_0$
becomes the two equations $Z_1=Z_2=0$ for the polynomials
\begin{equation}\label{eq:AF-Z}
\begin{aligned}
Z_1&=(L_2+2)(1-C_1)\,C_2^2+\bigl(C_1^2L_2-2C_1L_2-L_2-4\bigr)C_2+L_2+2+C_1(L_2-2),\\
Z_2&=\mathfrak q\,(C_2^2+1)+\bigl((C_1^2+1)L_2(L_2+2)+4C_1\bigr)C_2,\qquad \mathfrak q=C_1^2-\bigl((L_2+1)^2+1\bigr)C_1+1 ,
\end{aligned}
\end{equation}
and that at the common zeros the number $-x_1$ of \eqref{eq:x} is
\begin{equation}\label{eq:AF-theta}
-x_1=\frac{(1+C_1^2)(1+C_2^2)+4C_1C_2}{(L_2+1)^2(C_2-C_1)^2+(1+C_1C_2)^2} .
\end{equation}
The denominator is positive, since $C_1=C_2$ and $C_1C_2=-1$ are incompatible for real numbers.

In the rest of this subsection $u,v,t,\eta$ are real variables, unrelated to the rod vectors $v_j$. Aksteiner's chart of Teo's locus
\cite[App.~A.4, (140) and the chart before (143)]{Aksteiner2026}, with his $r=L_2/(L_2+2)=(1-t)/(1+t)$, is
\begin{equation}\label{eq:AF-chart}
C_1=\frac{1-u}{1+u},\qquad C_2=\frac{1+v}{1-v},\qquad L_2=\frac{1-t}t,\qquad 0<u,v<1,\quad 0<t\le1 .
\end{equation}
It gives $0<C_1<1<C_2<\infty$. With
\begin{equation}\label{eq:AF-G}
\mathcal G_0(u,v)=u^2(1+v)+2uv^2+v^2-v,\qquad \mathcal G_1(u,v)=3u^2v^2-u^2v+2uv^2+v-1 ,
\end{equation}
the polynomials \eqref{eq:AF-Z} and the right side of \eqref{eq:AF-theta} become, at the point \eqref{eq:AF-chart},
\begin{equation}\label{eq:AF-Zchart}
\begin{gathered}
t^2Z_2=\frac{4\bigl(u^2-v^2+t^2(1+3u^2v^2)\bigr)}{(1+u)^2(1-v)^2},\qquad tZ_1=\frac{4\bigl(\mathcal G_0(u,v)+t\,\mathcal G_1(u,v)\bigr)}{(1+u)^2(1-v)^2},\\
\frac{(1+C_1^2)(1+C_2^2)+4C_1C_2}{(L_2+1)^2(C_2-C_1)^2+(1+C_1C_2)^2}=\frac{2t^2(1+u^2v^2)}{(u+v)^2+t^2(1-uv)^2} .
\end{gathered}
\end{equation}
So $Z_2=0$ is the graph $t^2=(v^2-u^2)/(1+3u^2v^2)$, which is the first equation of \cite[(143)]{Aksteiner2026}. The equation $Z_1=0$ is
$\mathcal G_0+t\mathcal G_1=0$, the second. Following \cite[after (143)]{Aksteiner2026} we put $u=\eta v$, and for
$(v,\eta)\in(0,1]\times[0,1]$ we let
\begin{equation}\label{eq:AF-rect}
\begin{gathered}
t(v,\eta)=v\Bigl(\frac{1-\eta^2}{1+3\eta^2v^4}\Bigr)^{1/2},\qquad
\mathcal G(v,\eta)=\mathcal G_0(\eta v,v)+t(v,\eta)\,\mathcal G_1(\eta v,v),\\
\hat\vartheta(v,\eta)=\frac{(1-\eta)(1+\eta^2v^4)}{1-\eta(1-\eta)v^2+\eta^2(2+\eta)v^4} .
\end{gathered}
\end{equation}
These functions are continuous, and the denominator of $\hat\vartheta$ is at least $1-\frac14v^2>0$.

For $0<v<1$ and $0<\eta<1$, the point \eqref{eq:AF-chart} with $u=\eta v$ and $t=t(v,\eta)$ has $0<t<1$, so $L_2>0$. It satisfies
$Z_2=0$ by \eqref{eq:AF-Zchart}. If moreover $\mathcal G(v,\eta)=0$, then $Z_1=0$. Inserting $t^2$ into the last identity of
\eqref{eq:AF-Zchart} shows that $-x_1=\hat\vartheta(v,\eta)$ there. The derivative of $\hat\vartheta$ in $\eta$ is negative for
$0<\eta\le1$, so $\hat\vartheta(v,\cdot)$ decreases strictly from $1$ at $\eta=0$ to $0$ at $\eta=1$. For $0<v\le\frac13$ and
$0\le\eta<1$ both $\mathcal G_0(\eta v,v)$ and $\mathcal G_1(\eta v,v)$ are negative, so $\mathcal G(v,\eta)<0$. At $v=1$,
$\mathcal G_0(\eta,1)=\mathcal G_1(\eta,1)=2\eta(1+\eta)$, so $\mathcal G(1,\eta)>0$ for $\eta>0$. Let $\beta\in(0,1)$. For
$v\in[\frac13,1]$ let $\eta_\beta(v)\in(0,1)$ be the solution of $\hat\vartheta(v,\eta)=\beta$. It depends continuously on $v$, and the function
$\mathcal G(v,\eta_\beta(v))$ is negative at $v=\frac13$ and positive at $v=1$. So the intermediate value theorem gives $v_*\in(\frac13,1)$
at which it vanishes. The corresponding point is a positive solution of \eqref{eq:AF-lin} on Teo's locus with $-x_1=\beta$,
$0<C_1<1$, and $C_2=(1+v_*)/(1-v_*)>2$.

Since this solution satisfies $V_4=V_0$, it has $\alpha(\omega,\gamma)=0$ by \cite[Prop.~5.8(ii)]{Gerochcf}, and
Theorem~\ref{thm:converse} gives a simply connected toric gravitational instanton with four turning points, an AF end, and the chain
$(-1,0,1)$. Its twist is $\beta$ by Lemma~\ref{lem:twist} and the remark after it, so by Proposition~\ref{prop:AF-main} it is, up to
scale, the instanton of Theorem~\ref{thm:main}(d) with the twist $\beta$. Since its parameters are finite, \S\ref{subsec:AF-sym} shows
that it is, up to scale, Teo's metric with these parameters and his identifications, on the branch $C_2>2$ that Teo considers
\cite[after (8)]{Teo2026}. Hence every instanton of Theorem~\ref{thm:main}(d) is, up to scale, one of Teo's metrics and thus a member of
Teo's family on $X_3$. Conversely, every member is an instanton as in Theorem~\ref{thm:main}(d), so its twist lies in $(0,1)$. With
Proposition~\ref{prop:AF-main} it follows that the set of twists of the members of Teo's family is the interval $(0,1)$.

\subsection{Teo's family: soliton data}\label{subsec:AF-family}
The convexity of Appendix~\ref{app:convex}, together with the bijectivity of the map from the weights to the rod lines, also yields three statements on the soliton data of Teo's family.

\begin{remark}\label{rem:AF-family}
The following three statements are not used in the proofs of this paper, and we only sketch their proofs. First, every solution of $V_3=-V_1$, $V_4=V_0$ with positive data \eqref{eq:AF-positive}
and $1<w_3<w_4$ lies on Teo's locus. Second, the positive solutions of \eqref{eq:AF-lin} on Teo's locus are the solutions of
$Z_1=Z_2=0$ of \eqref{eq:AF-Z} (\S\ref{subsec:AF-teo}) with $L_2>0$ and $|C_1|<1<|C_2|$. They form one real-analytic arc, a graph
$L_2\mapsto(C_1(L_2),C_2(L_2))$ over $(0,\infty)$. Along it the number $\vartheta=-x_1$, the twist, decreases strictly from $1$ as
$L_2\to0$ to $0$ as $L_2\to\infty$, with $d\vartheta/dL_2\ne0$. Third, the soliton data of each instanton of Theorem~\ref{thm:main}(d) are
unique up to the changes of Lemma~\ref{lem:normalise}.

All three rest on the fact that positive data with $\alpha=0$ on given turning points are determined, up to a factor, by the lines of their
vectors $\mathbf n_j$ (\S\ref{subsec:appC-rods}). Indeed, the pencil transport of Lemma~\ref{lem:appD-pencil} shows that the map from the weights
to these lines has an invertible differential at every point (Lemma~\ref{lem:appC-local}). By \cite[Thm.~9.10]{Gerochcf} this map is proper
of degree $\pm1$. For the first statement, the reflection of the axis composed with the duality \eqref{eq:duality} maps a positive solution
of $V_3=-V_1$, $V_4=V_0$ to positive data with $\alpha=0$, the rod lengths in the reverse order, and the same lines. The convexity of the
cone errors at these fixed lines, proved as Proposition~\ref{prop:appC-convex}, gives $w_4-w_3=w_2-w_1$. The bijectivity then gives
the conditions on the $C_k$ in \eqref{eq:Teo-locus}. The relation $V_2=-V_1-V_0$ is not needed. For the second, a comparison of two
graphs in Aksteiner's chart \eqref{eq:AF-chart} shows that exactly one positive solution on Teo's locus has a given middle rod length.
The convexity with the bijectivity shows that each value of
$\vartheta$ in $(0,1)$ is taken by a unique normalised positive solution of \eqref{eq:AF-lin}. The two parametrisations are inverse to each other. The third follows from the convexity and the bijectivity, with
the balance of positive solutions in the proof of Proposition~\ref{prop:AF-main} and Lemma~\ref{lem:twist}. The first two statements have also been confirmed by exact computations with polynomials with integer coefficients.
Along the arc, $(C_1,C_2)$ tends to $(1,\infty)$ as $L_2\to0$ and to $(\frac12,2)$ as $L_2\to\infty$, in agreement with the limits that
Teo describes (Section~\ref{sec:explicit}).
\end{remark}

\section{Explicit solutions and metrics}\label{sec:explicit}
The soliton data, the metrics, and the topology of the instantons of Theorem~\ref{thm:main} are as follows.

The instantons (a) and (b) are the Gibbons--Hawking metrics \eqref{eq:GH} with four centres. In the soliton representation they are
\eqref{eq:closed} with $r=0$, $\omega=\frac12\sum_k\delta_{w_k}$, $\gamma=U_\infty$, and $\mathfrak c=|\det(V_0,V_1)|^{-1}=1$, with the lattice of periods
generated by $V_0=(1,2)$ and $V_1=(1,1)$ \cite[Prop.~11.3]{Gerochcf}
(Remark~\ref{rem:hk}). Here $w_1<\dots<w_4$ are the positions of the centres on the line. The rod lengths are the distances between consecutive centres, and the chain is
$(2,2,2)$, or $(-2,-2,-2)$ in the hyperk\"ahler orientation.

The ALF instanton (c) is \eqref{eq:closed} with the data \eqref{eq:ALF-data} and $r=2$. In the notation of \eqref{eq:closed} and \eqref{eq:PS}, its functions $\mathcal X,\mathcal Y,F_1,F_2$ are
\begin{equation}\label{eq:ALF-metric}
\begin{aligned}
\mathcal X&=\sum_{k<l}\varpi_k\varpi_la_ka_l(w_k-w_l)^2,\qquad
F_1=\sum_{k<l}R_kR_la_ka_l(w_k-w_l)^2+\gamma_{\rm ALF}\sum_kR_ka_k,\\
\mathcal Y&=\sum_{k<l<m}\varpi_k\varpi_l\varpi_m\,P_{klm}-2\gamma_{\rm ALF}\sum_{k<l}\varpi_k\varpi_l\,P_{kl}\,(w_k+w_l)
-\gamma_{\rm ALF}^2\sum_k\varpi_ka_k-i\gamma_{\rm ALF}\rho\,\mathcal X,\\
F_2&=\sum_{k<l<m}\frac{P_{klm}}{R_kR_lR_m}+\gamma_{\rm ALF}\sum_{k<l}\frac{P_{kl}}{R_kR_l},
\end{aligned}
\end{equation}
with $P_{kl}=a_ka_l(w_k-w_l)^2$ and $P_{klm}=a_ka_la_m(w_k-w_l)^2(w_k-w_m)^2(w_l-w_m)^2$. The constant in \eqref{eq:closed} is $\mathfrak c=|D_0|^{-1}$ with
\begin{equation}\label{eq:ALF-D0}
D_0=\det(V_0,V_1)=-\frac{160\bigl(139197+90142\sqrt2+57028\sqrt5+43998\sqrt{10}\bigr)}{3^4\cdot41^2}\approx-626.704,
\end{equation}
and the lattice of periods is generated by $\mathfrak c^{1/2}V_0$ and $\mathfrak c^{1/2}V_1$. By \eqref{eq:ALF-lin}, in the basis $(v_0,v_1)$ of the lattice
the rod vectors are
\begin{equation}\label{eq:ALF-rods}
v_0=(1,0),\qquad v_1=(0,1),\qquad v_2=(-1,-1),\qquad v_3=(-1,-2),\qquad v_4=(2,3) .
\end{equation}
The rod lengths are $1,\frac12\sqrt{10}-1,1$. Since $\det(v_0,v_1)=-1$ by \eqref{eq:det}, the basis $(v_0,v_1)$ is negatively oriented. The
coordinates of rod vectors in \eqref{eq:ALF-rods}, \eqref{eq:AF-rods}, \eqref{eq:AF-ell}, Figure~\ref{fig:rods}, and Appendix~\ref{app:convex}
are taken in such a basis, so the counter-clockwise rotations in the definition of the Li--Sun degree
(\S\ref{subsec:rods}) appear clockwise in them.

As in the proof of Lemma~\ref{lem:twist}, in the basis $(f_0,f_1)$ of Theorem~\ref{thm:soliton} the Killing field $f_1$ has bounded
length at infinity, while the length of $f_0$ grows like $\rho$. The first components of the rod vectors in this basis are
$\mathfrak c^{1/2}H_2(\omega_{(j)})$, proportional to the numbers $x_j$ of \eqref{eq:x}, which are $1,-\frac13,-\frac23,-\frac13,1$ by
Lemma~\ref{lem:seq}(b) with $\varepsilon=-1$. Hence $\xi=v_4-v_0=(1,3)$ lies on the line $\R f_1$. It is primitive, and the circle it
generates is the circle at infinity. Its isotropy group along the rod $I_j$ is its intersection with the stabiliser of $I_j$, whose
order is the index $|\det(v_j,\xi)|$ of $\Z v_j+\Z\xi$ in the lattice. By \eqref{eq:ALF-rods}, the orders along $I_0,\dots,I_4$ are
$3,1,2,1,3$.

Since $|\det(v_0,v_4)|=3$, the lattice $\Lambda'=\Z v_0+\Z v_4$ has index $3$ in $\Lambda$. The corresponding threefold cover of the end
has outer rod vectors that form a basis of its lattice, and the end is the quotient of this cover by $\Lambda/\Lambda'\cong\Z_3$. As
$\xi\in\Lambda'$ is primitive in $\Lambda$, we have $\R\xi\cap\Lambda=\Z\xi\subset\Lambda'$. Thus no non-trivial element of
$\Lambda/\Lambda'$ lies in the circle generated by $\xi$, and $\Z_3$ acts effectively on the base $\R^3$, by rotations about the axis.
The end is therefore the quotient of the end of the Taub--NUT metric by a group $\Z_3$ in the torus which acts on the base $\R^3$ by
rotations, as stated in \cite[Rem.~11.5(a)]{Gerochcf}. In particular it is not the end of a multi-Taub--NUT metric
\cite[Rem.~11.6]{Gerochcf}. By \cite[Prop.~8.8(ii), (8.37)]{Gerochcf}, the total NUT charge of soliton
data $(\omega,\gamma)$ is $-\frac14\gamma\,\alpha(\omega,\gamma)$. For the normalised data \eqref{eq:ALF-data}, by
\eqref{eq:ALF-alpha}, it is
\begin{equation}\label{eq:ALF-NUT}
-\tfrac14\alpha(\omega,\gamma_{\rm ALF})\,\gamma_{\rm ALF}=-\tfrac1{10086}\bigl(4265+2390\sqrt2+1700\sqrt5+1130\sqrt{10}\bigr)\approx-1.48916 .
\end{equation}

For the AF instantons (d), let $L_2>0$, and let $C_1,C_2$ with $|C_1|<1<|C_2|$ be such that the data on the turning points
$w=(0,1,1+L_2,2+L_2)$ and the parameters $\mathbf C=(C_1,C_2,-1/C_2,-1/C_1)$ satisfy \eqref{eq:AF-lin}. This is a positive solution
on Teo's locus, that is, a solution of the equations $Z_1=Z_2=0$ of \eqref{eq:AF-Z}. By \S\ref{subsec:AF-teo} there is such a solution, with $0<C_1<1<C_2$, for every twist $\vartheta=-x_1$ in
$(0,1)$. According to the second statement of Remark~\ref{rem:AF-family}, whose proof is sketched there, there is exactly one for
every $L_2>0$, namely $(C_1,C_2)=(C_1(L_2),C_2(L_2))$. The metric is \eqref{eq:closed} with $r=2$, the weights $a_k=-c_k/q'(w_k)$ with
$c_k=\gamma(1+C_k)/(1-C_k)$ from \eqref{eq:dictionary}, and any atom weight $\gamma>0$, since replacing $\gamma$ by $K\gamma$ replaces
$\omega$ by $K\omega$, which changes the metric only by a change of the angles (Lemma~\ref{lem:normalise} with $\lambda=1$, $\mu=0$).
The AF condition $\alpha(\omega,\gamma)=0$, that is, $\gamma^2=H_3(\omega)/H_1(\omega)$ by \eqref{eq:alpha} \cite[(8.36)]{Gerochcf},
holds on Teo's locus. The constant in \eqref{eq:closed} is $\mathfrak c=|\det(V_0,V_1)|^{-1}$, and by \eqref{eq:AF-lin}, in the basis
$(v_0,v_1)$ of the lattice, the rod vectors are
\begin{equation}\label{eq:AF-rods}
v_0=(1,0),\qquad v_1=(0,1),\qquad v_2=(-1,-1),\qquad v_3=(0,-1),\qquad v_4=(1,0) .
\end{equation}
If the parameters are finite, as for the solutions of \S\ref{subsec:AF-teo}, then by \S\ref{subsec:AF-sym} the metric is, up to
scale and a constant linear change of the angles, Teo's metric with the parameters $C_k$ \cite[(1)--(3)]{Teo2026}, which is therefore a
smooth complete instanton. According to Teo, the metrics tend to flat space as $L_2\to0$ and contain two Taub--bolt metrics with
opposite NUT charges as limits as $L_2\to\infty$ \cite{Teo2026}.

The invariants in Table~\ref{tab:topology} follow from \S\ref{subsec:rods}. Reversing the orientation of $\mathfrak t$
reverses $\kappa$ and changes none of the other invariants. Reversing the orientation of $M$ replaces $\kappa$ by $-(\kappa_3,\kappa_2,\kappa_1)$
and $Q_M$ by $-Q_M$, up to the order of the spheres. So Figure~\ref{fig:rods} shows (a), (b), and (d) with representatives of their
chains that differ from those in the table. The Euler characteristic
is the number of fixed points. The spheres $S_1,S_2,S_3$ form a basis of $H_2(M;\Z)$ with the intersection matrix
$Q_M=\mathrm{tridiag}(\kappa_j;1)$ by \eqref{eq:selfint}. From $Q_M$ the signs of the eigenvalues, the signature, $\det Q_M=-p$, and the parity
(even if and only if all $\kappa_j$ are even) are read off. Here $\Pl(\kappa_1,\kappa_2,\kappa_3)$ is the linear plumbing with Euler numbers $\kappa_j$ (\S\ref{subsec:rods}), $\Pl(-2,-2)$ the $A_2$
manifold, and $\Pl(3)$ the disc bundle of Euler number $3$ over $S^2$. The boundary of the plumbing $\Pl(\kappa_1,\kappa_2,\kappa_3)$, which up
to orientation is the boundary at infinity of $M$, is, up to
orientation, the lens space $L(|p|,\det(v_1,v_4))$ if $p\ne0$ \cite[Lemma~3]{Nilsson2024} and $S^2\times S^1$ if $p=0$. The diffeomorphism types follow by blowing down \cite{GS1999}. In a
linear plumbing, a sphere with Euler number $-1$, respectively $+1$, can be blown down. This splits off a summand $\CPb^2$, respectively
$\CP^2$, and raises, respectively lowers, the Euler numbers of its neighbours by one. The neighbours then meet in a point. Blowing down
the two outer spheres of $\Pl(-1,1,-1)$, respectively its middle sphere, gives the second, respectively the first, description in the
table. Blowing down the first and then
the last sphere of $\Pl(-1,0,1)$ gives $\Pl(0)\#\CP^2\#\CPb^2$. Its interior is $(\CP^2\#\CPb^2)\setminus S^1$, since the interior
$S^2\times\R^2$ of $\Pl(0)=S^2\times D^2$ is $S^4\setminus S^1$. By the
Chern--Gauss--Bonnet formula, $\frac1{32\pi^2}\int_M|\mathrm{Rm}|^2=\chi(M)-1/|\Gamma|$ for an ALE end with group $\Gamma$, and $=\chi(M)$ for ALF and
AF ends, including the twisted AF ends of (d) and the $\Z_3$-quotient of a Taub--NUT end of (c) \cite{EGH1980}, \cite[Lemma~4.40]{LS2025}. So
this energy is $\frac{15}4$ in case (a) and $4$ in cases (b), (c), and (d).

\begin{table}[ht]
\centering\small
\setlength{\tabcolsep}{3.5pt}
\begin{tabular}{>{\raggedright\arraybackslash}p{0.17\textwidth}lccccc>{\raggedright\arraybackslash}p{0.22\textwidth}}
\toprule
case & $\kappa$ & $\chi$ & $(n_+,n_-,n_0)$ & sign. & $\det Q_M$ & $\partial\Pl$ & manifold \\
\midrule
(a), (b) & $(-2,-2,-2)$ & $4$ & $(0,3,0)$ & $-3$ & $-4$, even & $L(4,1)$ & minimal resolution of $\C^2/\Z_4$ \\
(c) ALF & $(-1,1,-1)$ & $4$ & $(1,2,0)$ & $-1$ & $3$, odd & $L(3,1)$ & $\Pl(-2,-2)\#\CP^2\cong\Pl(3)\#2\,\CPb^2$ \\
(d) AF & $(-1,0,1)$ & $4$ & $(1,1,1)$ & $0$ & $0$, odd & $S^2\times S^1$ & $(\CP^2\#\CPb^2)\setminus S^1$ \\
\bottomrule
\end{tabular}
\caption{Topology of the instantons with four fixed points, in the orientation in which the chain is $\kappa$ (for (a), (b) the
hyperk\"ahler one). The columns give the Euler characteristic, the numbers of positive, negative, and zero eigenvalues of the intersection
form $Q_M$, the signature, $\det Q_M$ with its parity, the boundary of the plumbing $\Pl(\kappa)$, and the diffeomorphism type.}\label{tab:topology}
\end{table}

\section{The way forward}\label{sec:outlook}
Of the tools of the proof, the fold of Lemma~\ref{lem:appALF-fold} and Aksteiner's chart of Teo's locus are specific to four turning
points, whereas the convexity argument extends to any number of them \cite{ADNprep}.

Beyond four turning points there are infinitely many toric instantons that are Hermitian for neither orientation, with AF
ends as well as with ALF ends. For AF ends, Li and Sun's manifolds $X_m$ carry an AF$_\beta$ instanton for every twist $\beta\in(0,1)$
\cite[Thm.~1.1]{LS2025}, which is Hermitian for neither orientation when $m\ge3$ (Section~\ref{sec:intro}). Aksteiner \cite{Aksteiner2026}
gives five infinite families of AF rod structures. He also gives existence for the $356$ rod structures of his catalogue with at most
twelve turning points, under hypotheses on the collisions of turning points that are not proved there. These results are proved without
the hypotheses in \cite{ADNprep}. Every member of the five families and each of the $356$ rod structures carries an
AF$_\beta$ instanton for every admissible twist $\beta$, in the sense explained below, and for each twist its \emph{balanced} rod lengths,
those at which all cone errors vanish, are unique. For Li and Sun's $X_m$ this is the uniqueness of the rod lengths conjectured in
\cite[Conj.~8.5]{LS2025}, proved in \cite{ADNprep}.

Aksteiner introduced the admissible twists of an AF rod structure \cite[\S5.4, App.~A]{Aksteiner2026}. They depend only on its chain
$W$. In the sign convention of \S\ref{subsec:rods}, which is minus Aksteiner's, they form an open interval $I_W$ with rational
endpoints. It is cut out by the conditions $|x_j|<1$ on the angular components \eqref{eq:x} of the finite rods, which depend affinely on the
twist. Every AF instanton with the chain $W$ has its twist in $I_W$ \cite{ADNprep}. Aksteiner found that $I_W$ can be a proper subinterval of
$(0,1)$. His catalogue contains such rod structures from eight turning points on, the first with admissible twists forming an interval
of length $\frac23$. A computer enumeration of the catalogue shows that $159$ of its $356$ rod structures, all with eight to twelve
turning points, are of this kind. For Li and Sun's $X_m$, for Teo's family, and for every rod structure of the catalogue with at most
seven turning points, $I_W=(0,1)$.

For a chain $W$ of the five families or of the catalogue, the rod lengths of every AF instanton with the chain $W$ are balanced. By the
argument of the proof of Proposition~\ref{prop:AF-main} \cite{ADNprep}, such an instanton is determined up to scale by its twist and its
rod lengths. It follows that the AF instantons with the chain $W$ form, up to scale, a one-parameter family
over $I_W$: for each $\beta\in I_W$ there is exactly one, and for $\beta\notin I_W$ there is none. The AF instantons with these chains
are thus classified by pairs of a chain and a twist in its interval, and the uniqueness for a given twist, as in
Theorem~\ref{thm:main}(d) and \cite[Conj.~8.5]{LS2025}, is a statement on $I_W$. How the instantons degenerate as $\beta$ tends to an
endpoint of $I_W$, where the angular component of a finite rod tends to $\pm1$, is not known, and Aksteiner \cite{Aksteiner2026}
asserts no bound there.

For ALF ends, it is proved in \cite{ADNprep} that for all $a,b\ge0$ the chain
\begin{equation}\label{eq:kappaab}
\kappa_{a,b}=(-1,\underbrace{0,\dots,0}_{2a},1,\underbrace{0,\dots,0}_{2b},-1)
\end{equation}
is, up to scale, the chain of exactly one simply connected toric instanton with an ALF end, for the orientation in which the end is
ALF$^-$. This instanton has $4+2(a+b)$ turning points and is Hermitian for neither orientation, and $\kappa_{0,0}$ is the chain of
$g_{\rm ALF}$. There are therefore infinitely many toric ALF instantons that are Hermitian for neither orientation, besides the
infinitely many multi-Taub--NUT metrics, which are hyperk\"ahler. Building on Aksteiner's census of ALF rod structures with at most
eight turning points \cite{Aksteiner2026}, the work \cite{ADNprep} shows that a simply connected toric ALF instanton with four to eight
turning points that is hyperk\"ahler for neither orientation is, up to scale, one of $52$ instantons, each the only one with its chain.
For six of them and for the exclusion of sixteen chains, \cite{ADNprep} uses interval arithmetic.

With five turning points and an ALF end there are, up to scale, exactly two simply connected toric instantons that are
hyperk\"ahler for neither orientation \cite{ADNprep}. In the orientation in which the end is ALF$^-$ their chains are $(-1,1,0,1)$ and
$(1,-1,-1,1)$. These rod structures were found independently in \cite[Rem.~11.5, Table~4]{Gerochcf} and \cite[\S7.3]{DMT2026}. The
existence and the uniqueness of each instanton are proved by hand in \cite{ADNprep}, and both instantons are Hermitian for neither
orientation. For AF ends with five turning points there is, up to orientation and reversal, one chain that the necessary conditions
allow, $(1,0,0,1)$, that of Li and Sun's $X_4$. It carries an AF$_\beta$ instanton for every twist $\beta\in(0,1)$ \cite{LS2025}, whose
balanced rod lengths are unique for each twist. By the argument of the proof of Proposition~\ref{prop:AF-main}
\cite{ADNprep}, as above, this instanton is, up to scale, the only one with its chain and twist.

For ALE ends, the toric Bando--Kasue--Nakajima conjecture says that a simply connected toric Ricci-flat ALE instanton is hyperk\"ahler
for one orientation. This is the toric case of the conjecture of Bando, Kasue, and Nakajima that every Ricci-flat ALE $4$-manifold is
self-dual or anti-self-dual \cite{BKN1989}. Five turning points are the first case of the toric conjecture in which both degrees can be
$2$, by \eqref{eq:degid}, and there are only the hyperk\"ahler metrics. Indeed, if both Li--Sun degrees of a toric instanton are
positive, its rod vectors satisfy the \emph{block determinant bound}
\begin{equation}\label{eq:blockbound}
|\det(v_a,v_b)|<b-a,\qquad 0\le a<b\le n,\quad b-a\ge2
\end{equation}
\cite{ADNprep}. For five turning points this bound and the parallelogram condition \eqref{eq:parallelogram} hold for none of the
chains in $\{-1,0,1\}^4$ with $p\ne0$. For every number $n$ of turning points, the chains of the instantons that are hyperk\"ahler for
neither orientation lie in $\{-1,0,1\}^{n-1}$ by Theorem~\ref{thm:LSbound}. They are therefore finitely many, and the analogue of
Proposition~\ref{prop:reduction} is a finite enumeration. For five turning points a computer enumeration shows that the necessary conditions reduce the $81$ chains to the
cases above, and a proof by hand will appear in \cite{ADNprep}.

From six turning points on, the following new difficulties appear.
\begin{enumerate}[leftmargin=*,label=\textup{(\roman*)}]
\item For ALF ends with three to five turning points no instanton that is hyperk\"ahler for neither orientation has a rod parallel to the
asymptotic circle, that is, a finite rod whose angular component $x_j$ in \eqref{eq:x} vanishes. With six turning points such rod
structures occur, as they do for the Taub--bolt metric, which has two turning points. For them the transport of
Lemma~\ref{lem:appC-transport} need not exist, since it needs leading minors that do not vanish, and one of them vanishes when a rod is
parallel to the asymptotic circle. The pencil transport of Lemma~\ref{lem:appD-pencil} needs no such hypothesis. With it the
convexity argument of Appendix~\ref{app:convex}, extended to any number of turning points in \cite{ADNprep}, gives the strict
convexity for these rod structures as well, and with it the uniqueness of the balanced rod lengths for a given asymptotic
parameter. In \cite{ADNprep} this is used for the
ALF$^-$ rod structures in which the angular components of all finite rods lie in $(-1,1)$, the value $0$ included. Existence is
proved there by minimising the strictly convex function whose gradient is the vector of cone errors. Under conditions on the rod
vectors, it is shown there that near every face of the boundary of the space of rod lengths, where turning points collide, this function
decreases in a direction pointing into the interior, and that it increases at infinity. This uses Li and Sun's analysis of the
canonical harmonic maps when turning points collide and when the span of the rods tends to infinity \cite{LS2025}. The function then
attains its minimum in the interior, where all cone errors vanish. For the sixteen instantons with at most eight turning points that
have a rod parallel to the asymptotic circle, \cite{ADNprep} verifies these conditions with exact, computer-assisted data on the rod
vectors, without interval arithmetic.
\item For ALF ends, from six turning points on there are chains that satisfy the necessary conditions but carry no instanton. Their
exclusion uses explicit trial rod lengths at which all cone errors are positive, and these are so far enclosed by interval arithmetic.
From eight turning points on there are chains $\kappa$ such that $\kappa$ and $-(\kappa_{n-1},\dots,\kappa_1)$ both carry an ALF
instanton, for the orientation in which the end is ALF$^-$. The existence of these instantons is so far known by interval arithmetic.
\item For ALE ends with six turning points, one rod structure, with the chain $(1,-1,1,0,1)$, satisfies the block determinant bound
\eqref{eq:blockbound} and the parallelogram condition \eqref{eq:parallelogram}, but carries no instanton. An instanton would be a point
at which all cone errors take the common value $0$. The cone errors of an ALE rod structure do not change when all rod lengths are
multiplied by a common factor, and they are the gradient of a convex function of the rod lengths that is homogeneous of degree one and
strictly convex in the directions transversal to the rays \cite{ADNprep}. The rod lengths at which all cone errors are equal are therefore
unique up to scale. For this rod structure the common value there is $-\log c_*$ with $c_*\approx1.094$, so that all cone angles are
$2\pi c_*\ne2\pi$. The value of $c_*$ is known by interval arithmetic, not by hand. Numerical computations give rod structures
with about twenty turning points whose common cone error has absolute value as small as about $6\cdot10^{-10}$, so the toric
Bando--Kasue--Nakajima conjecture becomes an arithmetic problem.
\item Exact computer enumerations show that the number of chains that the necessary conditions leave, up to orientation and reversal,
grows. For AF ends it is $3$, $3$, $9$, $15$, and $40$ with six to ten turning points, and for ALF ends $6$, $16$, and $40$ with six to
eight.
\end{enumerate}

The smoothness conditions give the classification an arithmetic structure. Its discrete data are the chain, the type of the end,
the degree index, and, for ALF ends, the sign of $v_n$ relative to the asymptotic circle. For given discrete data, the smoothness
conditions are polynomial equations with integer coefficients in the normalised soliton data. Since the turning points are among the
unknowns, the rod lengths are not part of the data. The discrete data therefore define a scheme over $\Q$, the solutions of these
equations, and the instantons with these data are, up to scale, its real points with positive coordinates. Its dimension at an instanton
depends on the instanton. At the hyperk\"ahler instantons it is
positive, $n-2$ for ALE and $n-1$ for ALF ends (Remark~\ref{rem:hk}). At the data \eqref{eq:ALF-data} of $g_{\rm ALF}$ it is $0$: the
Jacobian determinant of the six equations \eqref{eq:ALF-lin} there is a non-zero element of $\F$, so these data are an isolated and
simple solution. For Teo's rod structures the scheme contains, near the instantons, the real-analytic arc of the second
statement of Remark~\ref{rem:AF-family}, along which the twist varies. We expect the same at the other instantons that are hyperk\"ahler for neither
orientation: dimension $0$ for ALF ends, and a curve along which the twist varies for AF ends. The number fields of the isolated
points, or for AF ends the function field of the curve, are the arithmetic of this scheme. For ALE ends the analogous object is the
scheme of the equations saying that all cone errors are equal. On it the toric Bando--Kasue--Nakajima conjecture asks whether the common cone error can vanish at
a point with positive data.

Where the normalised soliton data of an ALF instanton are an isolated solution of the smoothness system, a square polynomial system
with rational coefficients, they are algebraic. They generate a number field, that of the irreducible factor of an eliminant which vanishes
at them. For AF ends the same holds at algebraic twists. Nothing known predicts the degree of this field, but the examples suggest that
symmetry under the reflection of the axis lowers it. For $g_{\rm ALF}$ the field is $\F=\Q(\sqrt2,\sqrt5)$, of degree $4$. With five turning points,
the normalised soliton data computed for the instanton with the chain $(1,-1,-1,1)$, a palindrome, generate a field of degree $16$ with
Galois group $\mathfrak S_{16}$. Those computed for the instanton with the chain $(-1,1,0,1)$ generate a field of degree $412$ with Galois
group $\mathfrak S_{412}$. These fields and their Galois groups were determined by exact computation. Closed formulas like
\eqref{eq:ALF-data} are therefore not expected beyond four turning points.

\section*{Computation and the use of an AI assistant}
No proof in this paper depends on a computer calculation. The identification of the instantons of Theorem~\ref{thm:main}(d) with
Teo's metrics rests on computations that \S\S\ref{subsec:AF-sym} and~\ref{subsec:AF-teo} describe but do not reproduce. Computer
calculations also enter statements that no proof uses: the exact and decimal values in Sections~\ref{sec:ALF} and~\ref{sec:explicit},
the solutions off the symmetric locus in \S\ref{subsec:ALF-forcing}, the checks in Remark~\ref{rem:AF-family}, and the numerical
evidence, enumerations, interval arithmetic, and arithmetic data of Section~\ref{sec:outlook}, some of them taken from the cited works. The programs and their outputs are available from the authors on request.

The project was carried out with the help of an AI assistant, Claude (Anthropic), used through the Claude Code environment. The
assistant took part in developing the arguments and in writing the text, and it wrote and ran the verification computations. We also
used it to review the literature and to read the drafts critically. The authors take responsibility for the content.

\subsection*{Acknowledgements}
We are grateful to Steffen Aksteiner and Alfie
Jones for sharing their work before publication, and for helpful discussions.

\appendix
\section{The first components for ALF ends}\label{app:ALF}
Once the rod lengths of a positive solution of \eqref{eq:ALF-lin} are known, the proof of Proposition~\ref{prop:ALF-symmetric} is
completed by an analysis of the first components of the system. We use the notation of
Section~\ref{sec:ALF}, with $\kappa=(-1,1,-1)$ and $r=2$. For data $(\omega,\gamma)$ that are not necessarily normalised we write
$f_{j,s}(\omega,\gamma)$ for the $s$-th component of $V_{j+1}+V_{j-1}-\kappa_jV_j$, so that $f_{j,s}(\mathbf d)$ is its value at the
normalised data \eqref{eq:ALF-unknowns}. The first components $f_{j,1}$ do not involve $\gamma$. A M\"obius transformation of the axis
moves three turning points to $0$, $1$, and $\infty$ (\S\ref{subsec:appALF-gauge}), and in these coordinates the three equations
$f_{j,1}=0$ are solved explicitly. On their solutions with positive weights the cross-ratio of the turning points is at least $\frac52$,
with equality at exactly one point of the resulting plane curve, the fold (\S\ref{subsec:appALF-first}).

\subsection{Heine terms and a M\"obius gauge}\label{subsec:appALF-gauge}
The first components are sums of Heine terms, and a M\"obius gauge moves three turning points to $0$, $1$, and
$\infty$.

Write $P_S=P_S(\omega)$ and $P_{kl}=P_{\{k,l\}}$. As in the proof of Lemma~\ref{lem:strip}, $P_S(\omega_{(j)})=\varepsilon_j(S)P_S$
with $\varepsilon_j(S)=(-1)^{|S\cap\{1,\dots,j\}|}$. Put $\hat\varepsilon_j(S)=\varepsilon_{j+1}(S)+\varepsilon_{j-1}(S)-\kappa_j\varepsilon_j(S)$
for $j=1,2,3$. By \eqref{eq:heine},
\begin{equation}\label{eq:appALF-f}
f_{j,1}=\sum_{|S|=2}\hat\varepsilon_j(S)\,P_S ,
\end{equation}
and explicitly
\begin{equation}\label{eq:appALF-f1}
\begin{aligned}
f_{1,1}&=P_{12}-P_{13}-P_{14}+P_{23}+P_{24}+3P_{34},\qquad f_{2,1}=-P_{12}+P_{13}-P_{14}+3P_{23}+P_{24}-P_{34},\\
f_{3,1}&=3P_{12}+P_{13}-P_{14}+P_{23}-P_{24}+P_{34}.
\end{aligned}
\end{equation}

\begin{lemma}\label{lem:appALF-first}
Let $\omega$ be a positive measure on points $w_1<\dots<w_4$ and $f_{j,1}$ as in \eqref{eq:appALF-f}. Then $f_{1,1}=f_{2,1}=f_{3,1}=0$ if
and only if
\begin{equation}\label{eq:appALF-Psolved}
P_{14}=P_{12}+P_{13}+2P_{23},\qquad P_{24}+P_{34}=P_{12}+P_{13},\qquad P_{24}-P_{34}=2P_{12}-P_{23}.
\end{equation}
\end{lemma}

\begin{proof}
By \eqref{eq:appALF-f1}, $f_{2,1}+f_{3,1}$, $f_{1,1}-f_{3,1}$, and $f_{3,1}-f_{2,1}$ are $-2$, $2$, and $-2$ times the differences of the
two sides of the three equations \eqref{eq:appALF-Psolved}. The matrix of these combinations, with the rows $(0,1,1)$,
$(1,0,-1)$, and $(0,-1,1)$, has determinant $-2$.
\end{proof}

We call numbers $t>1$, $y>0$, and $b_1,b_2,b_3,b_\infty>0$ \emph{gauge data}. For such data put $(\zeta_1,\zeta_2,\zeta_3)=(0,1,t)$ and
$\hat\zeta_k=\zeta_k+y$, and let $\omega_{t,y,b}$ be the measure with
\begin{equation}\label{eq:appALF-gauge}
w_k=\frac1y-\frac1{\hat\zeta_k},\quad a_k=b_k\,\hat\zeta_k^2\quad(k=1,2,3),\qquad w_4=\frac1y,\quad a_4=b_\infty .
\end{equation}
It is the image of $\sum_kb_k\delta_{\zeta_k}+b_\infty\delta_\infty$ under the M\"obius map $\varphi(\zeta)=\zeta/(y(\zeta+y))$, which
sends $-y$ to $\infty$ and $\infty$ to $1/y$, with the weights divided by $\varphi'(\zeta_k)=\hat\zeta_k^{-2}$. Since
$(\varphi(z)-\varphi(z'))^2=\varphi'(z)\varphi'(z')(z-z')^2$, this change of weights preserves the Heine terms of pairs, which
therefore do not involve $y$:

\begin{lemma}\label{lem:appALF-gauge}
Let $(t,y,b)$ be gauge data and $\omega_{t,y,b}$ the measure \eqref{eq:appALF-gauge}. Then the following hold.
\begin{enumerate}[leftmargin=*,label=\textup{(\alph*)}]
\item $0=w_1<w_2<w_3<w_4$, and $P_{12}=b_1b_2$, $P_{13}=b_1b_3t^2$, $P_{23}=b_2b_3(t-1)^2$, and $P_{k4}=b_kb_\infty$ for $k\le3$.
\item For $L_2,L_3,a_2,a_3,a_4>0$ the numbers
\begin{equation}\label{eq:appALF-ty}
\begin{gathered}
y=L_2+L_3,\qquad t=\frac{(1+L_2)(L_2+L_3)}{L_3},\\
(b_1,b_2,b_3,b_\infty)=\Bigl(1,\,\frac{a_2y^2}{(1+y)^2},\,\frac{a_3y^2}{(t+y)^2},\,a_4y^2\Bigr)
\end{gathered}
\end{equation}
are gauge data, and for them the measure with the turning points $(0,1,1+L_2,1+L_2+L_3)$ and the weights $(1,a_2,a_3,a_4)$ is the
image of $\omega_{t,y,b}$ under $w\mapsto y(1+y)w$, $a\mapsto y^{-2}a$.
\end{enumerate}
\end{lemma}

\begin{proof}
To prove (a), note that since $0<\hat\zeta_1<\hat\zeta_2<\hat\zeta_3$, the numbers $w_k=1/\hat\zeta_1-1/\hat\zeta_k$ increase from $0$ and stay below $w_4=1/\hat\zeta_1$. Moreover
$w_k-w_l=(\zeta_k-\zeta_l)/(\hat\zeta_k\hat\zeta_l)$ for $k,l\le3$ and $w_k-w_4=-1/\hat\zeta_k$, which with \eqref{eq:appALF-gauge} and \eqref{eq:PS}
gives the pair terms. To prove (b), note that $t-1=L_2(1+L_2+L_3)/L_3>0$, so that the numbers \eqref{eq:appALF-ty} are gauge data. For $b_1=1$ the change maps $\omega_{t,y,b}$ to the turning points $(0,1,t(1+y)/(t+y),1+y)$ and the
weights $(1,b_2(1+y)^2/y^2,b_3(t+y)^2/y^2,b_\infty/y^2)$. For the values \eqref{eq:appALF-ty} these are the given ones.
\end{proof}

By \eqref{eq:appALF-ty}, $t$ is the cross-ratio $(w_3-w_1)(w_4-w_2)/\bigl((w_2-w_1)(w_4-w_3)\bigr)$ of the turning points.

\subsection{The first components and the fold}\label{subsec:appALF-first}
In the gauge the first components are solved explicitly (Proposition~\ref{prop:appALF-curve}). On the solutions with positive weights the
cross-ratio is at least $\frac52$, with equality at a single point (Lemma~\ref{lem:appALF-fold}), and the same holds for the positive
solutions of \eqref{eq:ALF-lin} (Corollary~\ref{cor:appALF-reduced}). Throughout,
\begin{equation}\label{eq:appALF-N}
N=\varrho^2t-\varrho t^2+3\varrho+t .
\end{equation}

\begin{proposition}\label{prop:appALF-curve}
Let $(t,y,b)$ be gauge data with $b_1=1$, and put $\varrho=tb_3/b_2$. Then $\omega_{t,y,b}$ satisfies $f_{1,1}=f_{2,1}=f_{3,1}=0$ if
and only if
\begin{equation}\label{eq:appALF-b}
b_\infty=\frac{t(\varrho t+1)}{\varrho+t},\qquad b_2=\frac{tN}{\varrho(\varrho+t)(t-1)^2},\qquad b_3=\frac{N}{(\varrho+t)(t-1)^2},
\end{equation}
and
\begin{equation}\label{eq:appALF-C}
\begin{aligned}
0=\mathcal C(\varrho,t)&:=2N^2+(\varrho+t)(\varrho t+1)\bigl(N-\varrho(t-1)^2\bigr)\\
&\phantom{:}=3\varrho^4t^2-5\varrho^3t^3+2\varrho^3t^2+15\varrho^3t+2\varrho^2t^3-6\varrho^2t^2+2\varrho^2t+20\varrho^2\\
&\phantom{:=}-5\varrho t^3+2\varrho t^2+15\varrho t+3t^2 .
\end{aligned}
\end{equation}
\end{proposition}

\begin{proof}
By Lemma~\ref{lem:appALF-gauge}(a) with $b_1=1$ and $b_3=\varrho b_2/t$, the equations \eqref{eq:appALF-Psolved}, the last two divided
by $b_2$, are
\begin{equation}\label{eq:appALF-Xgauge}
\begin{gathered}
b_\infty=(1+\varrho t)\,b_2+\frac{2\varrho(t-1)^2}{t}\,b_2^2,\\
\Bigl(1+\frac\varrho t\Bigr)b_\infty=1+\varrho t,\qquad \Bigl(1-\frac\varrho t\Bigr)b_\infty=2-\frac{\varrho(t-1)^2}{t}\,b_2 .
\end{gathered}
\end{equation}
The second is the formula for $b_\infty$. Then the third reads $\varrho(t-1)^2b_2/t=2-(t-\varrho)(\varrho t+1)/(\varrho+t)=N/(\varrho+t)$, which is the
formula for $b_2$. The first, multiplied by $\varrho(\varrho+t)^2(t-1)^2/t$, becomes $\varrho(\varrho+t)(\varrho t+1)(t-1)^2=(\varrho+t)(\varrho t+1)N+2N^2$,
that is, $\mathcal C=0$. The steps are reversible.
\end{proof}

\begin{lemma}[the fold]\label{lem:appALF-fold}
Let $\varrho,t>0$ with $\mathcal C(\varrho,t)=0$. Then $t\ge\frac52$, with equality if and only if $\varrho=1$, and $N(\varrho,t)>0$.
\end{lemma}

\begin{proof}
Let $s=\varrho+1/\varrho\ge2$, with $s=2$ only for $\varrho=1$. Dividing $N$, $(\varrho+t)(\varrho t+1)$, and $N-\varrho(t-1)^2$ by $\varrho$ gives
$\mathcal C(\varrho,t)=\varrho^2\bar{\mathcal C}(s,t)$ with
\begin{equation}\label{eq:appALF-Cbar}
\bar{\mathcal C}(s,t)=2(ts-t^2+3)^2+(ts+t^2+1)(ts-2t^2+2t+2)=3t^2s^2-(5t^3-2t^2-15t)\,s+2t^3-12t^2+2t+20 .
\end{equation}
Here $\bar{\mathcal C}(2,t)=4(5-2t)(t+1)^2$ and $\partial_s\bar{\mathcal C}=t(6ts-5t^2+2t+15)$, which for $s\ge2$ and $0<t\le\frac52$
is at least $t\bigl(15+t(14-5t)\bigr)>0$. So for $0<t\le\frac52$ and $s\ge2$, $\bar{\mathcal C}(s,t)\ge\bar{\mathcal C}(2,t)\ge0$, with
equality only for $(s,t)=(2,\frac52)$. Hence $t\ge\frac52$, $t=\frac52$ forces $\varrho=1$, and $\varrho=1$ forces $\bar{\mathcal C}(2,t)=0$,
that is, $t=\frac52$. At $(\varrho,t)=(1,\frac52)$, $N=\frac74$. For $t>\frac52$ the quadratic polynomial $\bar{\mathcal C}(\cdot,t)$ has the
leading coefficient $3t^2>0$ and is negative at $s=2$ and at $s_0=t-3/t$. At $s_0$ the first expression in \eqref{eq:appALF-Cbar}
gives $\bar{\mathcal C}(s_0,t)=-2(t+1)(t-1)^3$. So both lie between its roots, its root $s\ge2$ satisfies $s>s_0$, and
$N/\varrho=ts-t^2+3=t(s-s_0)>0$.
\end{proof}

\begin{corollary}\label{cor:appALF-reduced}
Let $\mathbf d$ be a solution of \eqref{eq:ALF-lin} with positive coordinates, let $(t,y,b)$ be given by \eqref{eq:appALF-ty}, and put
$\varrho=tb_3/b_2$. Then $(t,y,b)$ are gauge data with $b_1=1$, $\varrho>0$, the weights satisfy \eqref{eq:appALF-b}, and $\mathcal C(\varrho,t)=0$.
In particular $t\ge\frac52$, with equality if and only if $\varrho=1$.
\end{corollary}

\begin{proof}
By Lemma~\ref{lem:appALF-gauge}(b), $(t,y,b)$ are gauge data with $b_1=1$, so $\varrho>0$, and the normalised measure $\omega$ of $\mathbf d$
is mapped to $\omega_{t,y,b}$ by $w\mapsto w/(y(1+y))$, $a\mapsto y^2a$. This is the change of Lemma~\ref{lem:normalise} with
$\lambda=1/(y(1+y))$, $\mu=0$, and $K=y^2$, under which the first components of the Hankel family are multiplied by $K^2\lambda^2$ by
\eqref{eq:normB}. So $f_{j,1}(\omega_{t,y,b})=K^2\lambda^2f_{j,1}(\omega)=0$ for $j=1,2,3$. Proposition~\ref{prop:appALF-curve} gives
\eqref{eq:appALF-b} and $\mathcal C(\varrho,t)=0$, and Lemma~\ref{lem:appALF-fold} gives the last statement.
\end{proof}

At a positive point of the symmetric locus \eqref{eq:locus}, where $L_3=1$ and $a_2=a_3$, \eqref{eq:appALF-ty} gives
$t=(1+L_2)^2=y^2$, $b_2/b_3=(t+y)^2/(1+y)^2=y^2=t$, and $b_\infty=a_4y^2=t$, that is, $\varrho=1$, $y^2=t$, and $b_\infty=t$. Conversely,
$y^2=t$ gives $L_3=1$, and then $\varrho=1$ gives $a_2=a_3$ and $b_\infty=t$ gives $a_4=1$. So in the gauge the positive part of the symmetric
locus is the set $\{\varrho=1,\ y^2=t,\ b_\infty=t\}$. On the solutions of \eqref{eq:appALF-b} the last condition follows from $\varrho=1$. By Lemma~\ref{lem:appALF-fold} the points of this set that satisfy $f_{1,1}=f_{2,1}=f_{3,1}=0$ lie over the
single point $(\varrho,t)=(1,\frac52)$ of the curve $\mathcal C=0$, $\varrho>0$, the fold.
\section{Uniqueness by convexity}\label{app:convex}
For fixed rod lines, Li and Sun's canonical harmonic map assigns to every vector $L\in\R^3_{>0}$ of rod lengths a vector
$\mathfrak b(L)$ of cone errors. This vector vanishes if and only if the scaled rod vectors along the finite rods are the given rod vectors
(\S\ref{subsec:appC-rods}). For rod structures with four turning points, an end of type ALF$^-$ or AF$_0$, and a constant asymptotic
parameter $\alpha\ge0$, we show that $\mathfrak b$ is the gradient of a strictly convex function (Proposition~\ref{prop:appC-convex}), so
that it has at most one zero. This adapts Aksteiner's convexity argument for AF rod structures \cite{Aksteiner2026}. For the rod
structures of type ALF$^-$ of the chain $(-1,1,-1)$, with the parameter $\alpha>0$ and the cone errors $\mathfrak b(\cdot,\alpha)$, we show
that the zeros for all $\alpha$ lie on one ray (Corollary~\ref{cor:appC-ray}), since $\mathfrak b(cL,c\alpha)=\mathfrak b(L,\alpha)$ for
$c>0$ (Lemma~\ref{lem:appC-scale}), and that every positive solution of \eqref{eq:ALF-lin} gives a zero
(Proposition~\ref{prop:appC-balanced}). The proof of Proposition~\ref{prop:ALF-symmetric} in \S\ref{subsec:ALF-forcing} rests on these
three results, and the uniqueness for AF ends in \S\ref{subsec:AF-main} on the convexity. Throughout, the
degree index is $r=2$ except in the paragraph before Lemma~\ref{lem:appC-sign} and in the proof of Proposition~\ref{prop:appC-ident}.

\subsection{The rod structures and their canonical maps}\label{subsec:appC-rods}
We identify the canonical maps of the rod structures $R_L$ defined below with far-node members of the Hankel family of
\cite{Gerochcf}, whose scaled rod vectors are explicit Hankel quantities (Proposition~\ref{prop:appC-ident}).

Let $\alpha\ge0$, and let $\ell_0,\dots,\ell_4\in\R^2$ be vectors with
\begin{equation}\label{eq:appC-adm}
\ell_0=\bigl(1,\tfrac\alpha2\bigr),\qquad \ell_4=\bigl(1,-\tfrac\alpha2\bigr),\qquad \det(\ell_{j-1},\ell_j)<0\quad(1\le j\le4),\qquad
\ell^{(1)}_1<0,\quad \ell^{(1)}_2<0,
\end{equation}
where $\ell^{(1)}_j$ is the first component of $\ell_j$. The vectors \eqref{eq:appC-ell} of \S\ref{subsec:appC-ALF}, with $\alpha>0$, and the
vectors \eqref{eq:AF-ell} of \S\ref{subsec:AF-main}, with $\alpha=0$, are of this kind. We identify the Lie algebra of the torus
with $\R^2$ by a basis, in which the Weyl coordinates, the harmonic maps, and the models \cite[(2.13)]{Gerochcf} are computed. In this basis the
model $\Phi^\sharp_{1,\alpha}$ of \cite[(2.13)]{Gerochcf} is Li and Sun's model $\mathrm{TN}^-_\alpha$ \cite[Def.~4.4]{LS2025} if $\alpha>0$,
and their model $\Phi^{\mathrm{AF}_0}$ \cite[Def.~4.6]{LS2025} if $\alpha=0$. For $L\in\R^3_{>0}$ put
\begin{equation}\label{eq:appC-wL}
w(L)=\bigl(0,\,L_1,\,L_1+L_2,\,L_1+L_2+L_3\bigr),
\end{equation}
and let $R_L$ be the rod structure \cite[Def.~3.1]{LS2025} with the turning points $w(L)$, the rod lines $[\ell_0],\dots,[\ell_4]$, and the
asymptotic type ALF$^-$ if $\alpha>0$, AF$_0$ if $\alpha=0$. If $\alpha>0$, then, since $[\ell_0]$ and $[\ell_4]$ are the lines of
$(1,\pm\frac\alpha2)$, the rotation $P_{[\ell_0],[\ell_4]}$ of \cite[(1.10)]{LS2025} is the identity and $[\ell_0]\wedge[\ell_4]=\alpha$
\cite[(1.11)]{LS2025}. So the model of $R_L$ \cite[(4.7)]{LS2025} is $\Phi^\sharp_{1,\alpha}$. If $\alpha=0$, then $[\ell_0]=[\ell_4]$ is the
line of $(1,0)$, the model of $R_L$ is $\Phi^\sharp_{1,0}$, with the model augmentation $0$ \cite[Def.~4.6]{LS2025}. The twist of the end is then
recorded by the lattice, as in \cite[Def.~2.11]{Gerochcf}. In both cases the model augmentation is $\nu^\sharp=\frac12\log(1+\alpha/(2R))$,
$R=(\rho^2+z^2)^{1/2}$. Let $\Phi_L$ be the unique harmonic map strongly tamed by $R_L$ \cite[Def.~4.10, Thm.~4.24]{LS2025}, and $\nu_L$ its
end-normalised augmentation, the unique augmentation with $\nu_L-\nu^\sharp=O(R^{-1})$ \cite[Prop.~4.26, Def.~4.27]{LS2025}. Let
$\tilde v_j(L)\in[\ell_j]$ be its scaled rod vectors \cite[Def.~4.17]{LS2025}, which are defined up to sign. The \emph{cone errors}
$\mathfrak b_j(L)$ are defined by
\begin{equation}\label{eq:appC-b}
\tilde v_j(L)=\pm e^{\mathfrak b_j(L)}\,\ell_j,\qquad 0\le j\le4,
\end{equation}
and we write $\mathfrak b(L)=(\mathfrak b_1,\mathfrak b_2,\mathfrak b_3)(L)$. We call $L$ \emph{balanced} if $\mathfrak b(L)=0$. These
objects depend on $\ell=(\ell_0,\dots,\ell_4)$. Where this matters we write $R_{L,\ell}$, $\Phi_{L,\ell}$, $\nu_{L,\ell}$, and
$\mathfrak b(L,\ell)$.

For a positive measure $\omega=\sum_{k=1}^4a_k\delta_{w_k}$ on points $w_1<\dots<w_4$ and $\gamma>0$, write $H_s=H_s(\omega)$,
$\Sigma_2=\Sigma_2(\omega)$, and $\alpha=\alpha(\omega,\gamma)$ as in \eqref{eq:H3gamma} and \eqref{eq:alpha}, and put, for $0\le j\le4$,
\begin{equation}\label{eq:appC-n}
\mathbf n_j(\omega,\gamma)=\Bigl(x_j,\ \frac{\hat Y_j}{\gamma H_2}\Bigr),\qquad x_j=\frac{H_2(\omega_{(j)})}{H_2},\qquad
\hat Y_j=H^{(\gamma)}_3(\omega_{(j)})+2\gamma\,\Sigma_2\,x_j ,
\end{equation}
with $\omega_{(j)}$ as in \eqref{eq:V}. These are the vectors $T_0\hat V_j$ of \cite[(8.42)]{Gerochcf} for $r=2$: $\mathbf n_j=A_{\omega,\gamma}V_j$ for the Hankel vectors $V_j$ of \eqref{eq:V} and the matrix
\begin{equation}\label{eq:appC-A}
A_{\omega,\gamma}=\begin{pmatrix}H_2^{-1}&0\\0&(\gamma H_2)^{-1}\end{pmatrix}\begin{pmatrix}1&0\\2\gamma\Sigma_2/H_2&1\end{pmatrix},
\qquad \det A_{\omega,\gamma}=\frac1{\gamma H_2^2}>0 .
\end{equation}
Since $\omega_{(4)}=-\omega$, $H_s(-\omega)=(-1)^sH_s(\omega)$, and $\Sigma_2(-\omega)=\Sigma_2(\omega)$, \eqref{eq:H3gamma} and \eqref{eq:alpha} give
\begin{equation}\label{eq:appC-nouter}
\mathbf n_0(\omega,\gamma)=\bigl(1,\tfrac\alpha2\bigr),\qquad \mathbf n_4(\omega,\gamma)=\bigl(1,-\tfrac\alpha2\bigr),\qquad \alpha=\alpha(\omega,\gamma),
\end{equation}
and $\det(\mathbf n_{j-1},\mathbf n_j)=\det A_{\omega,\gamma}\det(V_{j-1},V_j)<0$ for $1\le j\le4$ \cite[Prop.~5.8(i)]{Gerochcf}. For soliton data of an instanton, $x_j$
is the number of \eqref{eq:x}. By Lemma~\ref{lem:normalise}, with $(\lambda,\mu,K)=(1,\mu,1)$, respectively
$(1,0,c)$, respectively $(c,0,1)$, the translated measure $\sum_ka_k\delta_{w_k+\mu}$ with the same $\gamma$, the data $(c\omega,c\gamma)$,
and the dilated measure $\sum_ka_k\delta_{cw_k}$ with $c^3\gamma$ satisfy
\begin{equation}\label{eq:appC-nscale}
\begin{gathered}
\mathbf n_j\Bigl(\sum_ka_k\delta_{w_k+\mu},\gamma\Bigr)=\mathbf n_j(\omega,\gamma),\qquad
\mathbf n_j(c\omega,c\gamma)=\mathbf n_j(\omega,\gamma),\\
\mathbf n_j\Bigl(\sum_ka_k\delta_{cw_k},c^3\gamma\Bigr)=\begin{pmatrix}1&0\\0&c\end{pmatrix}\mathbf n_j(\omega,\gamma),
\end{gathered}
\end{equation}
for $\mu\in\R$ and $c>0$. The values of $\alpha$ for these data are $\alpha$, $\alpha$, and $c\alpha$. Indeed, $H_s$ and $\Sigma_2$ of the translated
measure are $H_s$ and $\Sigma_2+2\mu H_2$, and $V_j$ is replaced by $\bigl(H_2x_j,V^{(2)}_j-4\mu\gamma H_2x_j\bigr)$, which leaves $\hat Y_j$
unchanged. The data $(c\omega,c\gamma)$ leave $x_j$ unchanged and multiply $H_2$, $\hat Y_j$, and $\gamma H_2$ by $c^2$, $c^3$, and $c^3$. The dilation
multiplies $H_s$, $\Sigma_2$, and $\hat Y_j$ by $c^{s(s-1)}$, $c^3$, and $c^6$.

For $\omega$ and $\gamma$ as above there is a far-node configuration whose four finite nodes are $w_1,\dots,w_4$, none of them heavy or
light, that is, with all parameters $c_k$ of the finite nodes finite and non-zero \cite[\S8, before (8.1)]{Gerochcf}. Its far-node
limit, as the fifth node $\Omega$ tends to infinity, has the family data $(\omega,2,\gamma,\mathfrak s)$ with $\mathfrak s=w_1+\dots+w_4$
\cite[Thms.~7.3(b), 8.3]{Gerochcf}. For a constant $\mathfrak c>0$ let $\nu_{\mathfrak c}$ be the augmentation given by $e^{2\nu_{\mathfrak c}}=\mathfrak c\,F_1F_2$, with $F_1,F_2$
as in \eqref{eq:closed} \cite[Thm.~8.3(ii), (8.12)]{Gerochcf}. We call the harmonic map of such a limit, read in the basis
$(a_\gamma,b_\gamma)$ of \cite[Prop.~8.9]{Gerochcf}, together with $\nu_{\mathfrak c}$, a \emph{far-node member} with data $(\omega,\gamma)$.
By \cite[Prop.~8.9 and the display before it]{Gerochcf}, a far-node member is locally tamed by the rod structure with the turning points
$w_1,\dots,w_4$ and the rod lines $[\mathbf n_0],\dots,[\mathbf n_4]$, and its distance from $\Phi^\sharp_{1,\alpha(\omega,\gamma)}$ is
$O(R^{-1})$, uniformly in the angle. The proof of \cite[Prop.~8.9(i)]{Gerochcf} writes the
rod vector \cite[(8.35)]{Gerochcf}, given there in the frame $(e_\gamma,e_1)$, as $\gamma^{1/2}H_2\mathbf n_j$ in the basis
$(a_\gamma,b_\gamma)$. By this and \cite[Prop.~8.7(i)]{Gerochcf}, the scaled rod vectors of a far-node member are
\begin{equation}\label{eq:appC-farnorm}
\pm \mathfrak c^{1/2}\gamma^{1/2}H_2\,\mathbf n_j(\omega,\gamma),\qquad 0\le j\le4 .
\end{equation}
Here \cite[Prop.~8.7(i)]{Gerochcf} is stated for the constant $\mathfrak c_\nu$ of \cite[Thm.~7.3]{Gerochcf}. By the definition of the scaled rod vectors in
\cite[\S2.2]{Gerochcf}, which agrees with that of \cite[Def.~4.17]{LS2025}, replacing it by $\mathfrak c$ multiplies the scaled rod vectors by $(\mathfrak c/\mathfrak c_\nu)^{1/2}$. No other statement used here
depends on the choice of the constant \cite[Rem.~2.3]{Gerochcf}. By \cite[Prop.~8.8(iii)]{Gerochcf}, $e^{2\nu_{\mathfrak c}}$ is a real-analytic
function of $(R^{-1},\theta)$ near $R^{-1}=0$, where $z+i\rho=Re^{i\theta}$, with the limit $\mathfrak c\gamma H_2^2$ as $R\to\infty$.

For an integer $0\le r\le4$ let $\mathbf n^{[r]}_j$, $0\le j\le4$, be the vectors $T_0\hat V_j$ of \cite[(8.42)]{Gerochcf} formed with the
degree index $r$, so that $\mathbf n^{[2]}_j=\mathbf n_j(\omega,\gamma)$. Then $\mathbf n^{[r]}_0=(1,\frac\alpha2)$ and
$\mathbf n^{[r]}_4=(-1)^r(1,-\frac\alpha2)$ \cite[display after (8.42)]{Gerochcf}. Moreover $\det(\mathbf n^{[r]}_{j-1},\mathbf n^{[r]}_j)<0$
\cite[Prop.~5.8(i)]{Gerochcf}, since they are obtained from the $V_j$ by one matrix of positive determinant. The
following sign property is used for these vectors.

\begin{lemma}[the sign property]\label{lem:appC-sign}
Let $\alpha\in\R$, and let $\mathbf N_0,\dots,\mathbf N_4$ and $\mathbf N'_0,\dots,\mathbf N'_4$ be vectors in $\R^2$ with
$\mathbf N_0=\mathbf N'_0=(1,\frac\alpha2)$, $\mathbf N_4=(1,-\frac\alpha2)$, $\mathbf N'_4=\pm(1,-\frac\alpha2)$, and
$\det(\mathbf N_{j-1},\mathbf N_j)<0$, $\det(\mathbf N'_{j-1},\mathbf N'_j)<0$ for $1\le j\le4$. If $\mathbf N'_j=\lambda_j\mathbf N_j$ with
$\lambda_j\in\R$ for $0\le j\le4$, then $\lambda_j>0$ for $0\le j\le4$, and $\mathbf N'_4=\mathbf N_4$.
\end{lemma}

\begin{proof}
We have $\lambda_0=1$, and $\lambda_{j-1}\lambda_j\det(\mathbf N_{j-1},\mathbf N_j)=\det(\mathbf N'_{j-1},\mathbf N'_j)<0$ gives
$\lambda_{j-1}\lambda_j>0$ for $1\le j\le4$. Hence $\lambda_j>0$ for $0\le j\le4$. The first component of $\mathbf N'_4=\lambda_4\mathbf N_4$
is $\pm1=\lambda_4>0$, so $\mathbf N'_4=\mathbf N_4$.
\end{proof}

\begin{proposition}\label{prop:appC-ident}\label{prop:appD-ident}
Let $\alpha\ge0$, let $\ell_0,\dots,\ell_4$ be as in \eqref{eq:appC-adm}, $L\in\R^3_{>0}$, and let $\Phi_L$, $\nu_L$, and $\mathfrak b(L)$ be as
in \eqref{eq:appC-b}. Then the following hold.
\begin{enumerate}[leftmargin=*,label=\textup{(\roman*)}]
\item There are a positive measure $\omega$ on the points $w(L)$ of \eqref{eq:appC-wL} and $\gamma>0$ with $\alpha(\omega,\gamma)=\alpha$
such that $[\mathbf n_j(\omega,\gamma)]=[\ell_j]$ for $1\le j\le3$.
\item Let $\omega$ be a positive measure on the points $w(L)$ and $\gamma>0$, with $\alpha(\omega,\gamma)=\alpha$ and
$[\mathbf n_j(\omega,\gamma)]=[\ell_j]$ for $1\le j\le3$. Then every far-node member with data $(\omega,\gamma)$ has the harmonic map
$\Phi_L$, it has the augmentation $\nu_L$ if and only if $\mathfrak c=1/(\gamma H_2(\omega)^2)$, and
\begin{equation}\label{eq:appC-nb}
\mathbf n_j(\omega,\gamma)=e^{\mathfrak b_j(L)}\,\ell_j,\qquad 0\le j\le4 .
\end{equation}
In particular $\mathfrak b_0(L)=\mathfrak b_4(L)=0$.
\end{enumerate}
\end{proposition}

\begin{proof}
To prove (i), we note that the vectors $\ell_0,\dots,\ell_4$ satisfy the hypotheses of \cite[Lemma~3.10]{Gerochcf}, by
\eqref{eq:appC-adm}. Let $r$ be its integer. By \cite[Prop.~9.11]{Gerochcf} there are a positive measure $\omega$ on the
points $w(L)$, $\gamma>0$ with $\alpha(\omega,\gamma)=\alpha$, and a far-node limit with the family data $(\omega,r,\gamma,\mathfrak s)$, for some
$\mathfrak s$, whose rod lines in the basis $(a_\gamma,b_\gamma)$ are $[\ell_0],\dots,[\ell_4]$. That $\gamma$ can be chosen with $\alpha(\omega,\gamma)=\alpha$ follows from
\cite[Rem.~8.10]{Gerochcf}, from part (b) of that remark if $\alpha>0$ and from
part (a), with $\gamma=\gamma_*(\omega)$, if $\alpha=0$. The condition on the rod lines says that $[\mathbf n^{[r]}_j]=[\ell_j]$ for
$0\le j\le4$ \cite[Prop.~8.9(i)]{Gerochcf}. Write $\mathbf n^{[r]}_j=\lambda_j\ell_j$. By Lemma~\ref{lem:appC-sign}, with
$\mathbf N_j=\ell_j$ and $\mathbf N'_j=\mathbf n^{[r]}_j$, $\lambda_j>0$ for $0\le j\le4$ and $\mathbf n^{[r]}_4=(1,-\frac\alpha2)$, so $r$
is even. The first component of $\mathbf n^{[r]}_j$ is $H_r(\omega_{(j)})/H_r(\omega)$, and since $\lambda_j>0$ it has the sign of
$\ell^{(1)}_j$, which is negative for $j=1,2$ by \eqref{eq:appC-adm}. This excludes $r=0$, where it is $1$, and $r=4$, where it is
$(-1)^j$ by Lemma~\ref{lem:strip}(i). So $r=2$, $\mathbf n^{[2]}_j=\mathbf n_j(\omega,\gamma)$, and (i) holds.

Under the hypotheses of (ii), a far-node member with data $(\omega,\gamma)$ is locally tamed by the rod structure with the turning points
$w(L)$ and the rod lines $[\mathbf n_j]$, which are the lines $[\ell_j]$, for $j=0,4$ by \eqref{eq:appC-nouter}. Its harmonic map is at
distance $O(R^{-1})$ from $\Phi^\sharp_{1,\alpha}$. The map $\Phi_L$ has the same two properties, since it is strongly tamed by $R_L$. By
\cite[Lemma~3.1]{Gerochcf} the two harmonic maps agree. Two augmentations of the same harmonic map differ by a constant
\cite[Rem.~2.3]{Gerochcf}, so $\nu_L=\nu_{\mathfrak c}$ for exactly one $\mathfrak c>0$. Since $e^{2\nu_{\mathfrak c}}$ is real-analytic in $(R^{-1},\theta)$
with the limit $\mathfrak c\gamma H_2^2$, and $\nu^\sharp=O(R^{-1})$, the condition $\nu_{\mathfrak c}-\nu^\sharp=O(R^{-1})$ holds if and only if $\mathfrak c\gamma H_2^2=1$.
For this constant \eqref{eq:appC-farnorm} gives $\tilde v_j(L)=\pm\mathbf n_j(\omega,\gamma)$. Write $\mathbf n_j=\lambda_j\ell_j$. By
Lemma~\ref{lem:appC-sign}, with $\mathbf N_j=\ell_j$ and $\mathbf N'_j=\mathbf n_j$, which satisfy its hypotheses by \eqref{eq:appC-adm},
\eqref{eq:appC-nouter}, and the inequality after it, $\lambda_j>0$. So $e^{\mathfrak b_j(L)}=\lambda_j$ by \eqref{eq:appC-b}, which is \eqref{eq:appC-nb}. By
\eqref{eq:appC-nouter}, $\lambda_0=\lambda_4=1$.
\end{proof}

We call data $(\omega,\gamma)$ on the points $w(L)$ \emph{direction matched} for $L$, or for $(L,\ell)$, if they satisfy the hypotheses of
Proposition~\ref{prop:appC-ident}(ii). By that proposition their vectors $\mathbf n_j$ determine $\mathfrak b(L)$, whatever data are chosen.

\subsection{Aksteiner's matrices and the dictionary}\label{subsec:appC-dict}
The vectors $\mathbf n_j$ can be expressed through the matrices of Aksteiner's real matrix reconstruction of solitons
\cite{Aksteiner2026} (Proposition~\ref{prop:appC-dict}). Aksteiner uses these matrices for AF data, at the weight scale at which
$\mathbf u^T\mathbf v=0$ below. The ALF case needs them at every scale, and we write out the proofs.

Let $\omega=\sum_ka_k\delta_{w_k}$ be a positive measure on points $w_1<\dots<w_4$ and $\gamma>0$. For a real measure $\tilde\omega$ on
these points with moments $h_q$ let $M_s(\tilde\omega)=(h_{i+l})_{i,l=0}^{s-1}$, so that $H_s(\tilde\omega)=\det M_s(\tilde\omega)$, and
write $M_s=M_s(\omega)$. Put $\mathsf w_b=(a_k^{1/2}w_k^b)_{k=1}^4\in\R^4$ for $b=0,1,2$, $\mathsf W=[\mathsf w_0,\mathsf w_1]$, and
$\mathsf W'=[\mathsf w_0,\mathsf w_1,\mathsf w_2]$, so that $\mathsf W^T\mathsf W=M_2$ and $\mathsf W'^T\mathsf W'=M_3$. Let $\mathsf Z$ be the diagonal
matrix with the entries $w_1,\dots,w_4$, and let $\mathbf e_0,\mathbf e_1$ and $\mathbf e_0,\mathbf e_1,\mathbf e_2$ be the unit vectors of $\R^2$ and
$\R^3$. Let $E_j$ be the diagonal $4\times4$ matrix whose first $j$ diagonal entries are $1$ and whose others are $0$. Define
\begin{equation}\label{eq:appC-TUV}
\begin{gathered}
\Pi=\mathsf WM_2^{-1}\mathsf W^T,\qquad T=2\Pi-I,\qquad \mathsf a=\mathsf WM_2^{-1}\mathbf e_1,\qquad \mathsf c=(I-\Pi)\,\mathsf w_2,\\
\mathbf u=\gamma^{1/2}\mathsf a+\gamma^{-1/2}\mathsf c,\qquad \mathbf v=\gamma^{1/2}\mathsf a-\gamma^{-1/2}\mathsf c .
\end{gathered}
\end{equation}
These are the matrices $P$, $T$ and the vectors $a$, $c$, $u$, $v$ of \cite{Aksteiner2026} for the degree index $2$ and the weights
$\omega/\gamma$. For $0\le j\le4$ let $T_j$ be the leading $j\times j$ block of $T$ and $\mathbf v_{[j]}$ the vector of the first $j$ entries
of $\mathbf v$, with $\det T_0=1$ and $\mathbf v_{[0]}^T\operatorname{adj}(T_0)\mathbf v_{[0]}=0$. Define Aksteiner's rods
\begin{equation}\label{eq:appC-rods}
\mathbf q_j=(-1)^j\bigl(\det T_j,\ 2\,\mathbf v_{[j]}^T\operatorname{adj}(T_j)\,\mathbf v_{[j]}\bigr),\qquad 0\le j\le4,
\end{equation}
which are those of \cite{Aksteiner2026} with the two components exchanged.

\begin{lemma}\label{lem:appC-interp}
Let $\Pi$, $T$, $\mathsf a$, $\mathsf c$, $\mathbf u$, $\mathbf v$ be as in \eqref{eq:appC-TUV}, and let $\mathcal W$ be the span of $\mathsf w_0$ and
$\mathsf w_1$. Then the following hold.
\begin{enumerate}[leftmargin=*,label=\textup{(\alph*)}]
\item $\Pi$ is the orthogonal projection onto $\mathcal W$, which has dimension $2$. Moreover $\mathsf a\in\mathcal W$, $\mathsf w_0^T\mathsf a=0$,
$\mathsf w_1^T\mathsf a=1$, $\mathsf c\perp\mathcal W$, and
\begin{equation}\label{eq:appC-xichi}
\|\mathsf a\|^2=\frac{H_1}{H_2},\qquad \mathsf w_2^T\mathsf a=\frac{\Sigma_2}{H_2},\qquad \mathsf w_2^T\mathsf c=\|\mathsf c\|^2=\frac{H_3}{H_2},\qquad
\|\mathbf v\|^2=\frac{H_3+\gamma^2H_1}{\gamma H_2} .
\end{equation}
\item $T$ is symmetric, $T^2=I$, $\det T=1$, $T\mathbf v=\mathbf u$, and
\begin{equation}\label{eq:appC-interp}
[\mathsf Z,T]=2(\mathsf c\mathsf a^T-\mathsf a\mathsf c^T),\qquad [\mathsf Z+\mathbf v\mathbf v^T,T]=0 .
\end{equation}
\item $\Pi$, $T$, $\mathbf u$, and $\mathbf v$ do not change if the points $w_k$ are translated by a common $\mu\in\R$, or if $(\omega,\gamma)$ is
replaced by $(c\omega,c\gamma)$ with $c>0$.
\end{enumerate}
\end{lemma}

\begin{proof}
To prove (a), note that since $\det M_2=H_2>0$, $\mathsf W$ has rank $2$, and $\Pi$ is symmetric with $\Pi^2=\Pi$ and range $\mathcal W$. As
$\mathsf W^T\mathsf a=\mathbf e_1$, $\mathsf a$ has the stated products with $\mathsf w_0,\mathsf w_1$, and $\|\mathsf a\|^2=\mathbf e_1^TM_2^{-1}\mathbf e_1=h_0/H_2=H_1/H_2$.
Next, $\mathsf w_2^T\mathsf a=(h_2,h_3)M_2^{-1}\mathbf e_1=(h_0h_3-h_1h_2)/H_2$, and
\begin{equation}\label{eq:appC-tau}
h_0h_3-h_1h_2=\frac12\sum_{k,l}a_ka_l\bigl(w_k^3+w_l^3-w_kw_l(w_k+w_l)\bigr)=\sum_{k<l}a_ka_l(w_k-w_l)^2(w_k+w_l)=\Sigma_2 .
\end{equation}
The vector $\mathsf c$ is orthogonal to $\mathcal W$, and $\mathsf w_2^T\mathsf c=\|\mathsf c\|^2=h_4-(h_2,h_3)M_2^{-1}(h_2,h_3)^T$ is the Schur complement of
$M_2$ in $M_3$, which is $H_3/H_2$. Since $\mathsf a\perp\mathsf c$, $\|\mathbf v\|^2=\gamma\|\mathsf a\|^2+\gamma^{-1}\|\mathsf c\|^2$.

To prove (b), note that $T^2=4\Pi^2-4\Pi+I=I$, and $T$ has the eigenvalues $1$ and $-1$, each twice, so $\det T=1$. Since $T\mathsf a=\mathsf a$ and $T\mathsf c=-\mathsf c$, $T\mathbf v=\mathbf u$.
For $x=\mathsf Wy\in\mathcal W$, $\mathsf Zx=y_0\mathsf w_1+y_1\mathsf w_2$, so $(I-\Pi)\mathsf Zx=y_1\mathsf c$ with $y_1=\mathbf e_1^TM_2^{-1}\mathsf W^Tx=\mathsf a^Tx$.
Hence $(I-\Pi)\mathsf Z\Pi=\mathsf c\mathsf a^T$, and by transposition $\Pi \mathsf Z(I-\Pi)=\mathsf a\mathsf c^T$. Writing $\mathsf Z\Pi=\Pi \mathsf Z\Pi+(I-\Pi)\mathsf Z\Pi$ and $\Pi \mathsf Z=\Pi \mathsf Z\Pi+\Pi \mathsf Z(I-\Pi)$
gives $[\mathsf Z,\Pi]=\mathsf c\mathsf a^T-\mathsf a\mathsf c^T$, which is the first identity of \eqref{eq:appC-interp}. Moreover
$[\mathbf v\mathbf v^T,T]=\mathbf v\mathbf u^T-\mathbf u\mathbf v^T=2(\mathsf a\mathsf c^T-\mathsf c\mathsf a^T)$, which gives the second.

To prove (c), note that a translation by $\mu$ replaces $\mathsf w_1$ by $\mathsf w_1+\mu\mathsf w_0$ and $\mathsf w_2$ by $\mathsf w_2+2\mu\mathsf w_1+\mu^2\mathsf w_0$.
It does not change $\mathcal W$, hence $\Pi$ and $T$. The vector $\mathsf a$ is the element of $\mathcal W$ with $\mathsf w_0^T\mathsf a=0$ and $\mathsf w_1^T\mathsf a=1$, so it
does not change, and $\mathsf c=(I-\Pi)\mathsf w_2$ does not change. Replacing $\omega$ by $c\omega$ multiplies each $\mathsf w_b$ by $c^{1/2}$, so it
leaves $\Pi$ unchanged and replaces $\mathsf a,\mathsf c$ by $c^{-1/2}\mathsf a,c^{1/2}\mathsf c$. With $c\gamma$ in place of $\gamma$, $\mathbf u$ and $\mathbf v$ do not change.
\end{proof}

\begin{proposition}[the dictionary]\label{prop:appC-dict}
Let $\omega$, $\gamma$ be as in \eqref{eq:appC-TUV}, $\alpha=\alpha(\omega,\gamma)$, $\mathbf n_j=\mathbf n_j(\omega,\gamma)$ as in
\eqref{eq:appC-n}, and $\mathbf q_j$ as in \eqref{eq:appC-rods}. Then
\begin{equation}\label{eq:appC-dict}
\mathbf n_j=B_\alpha\mathbf q_j,\qquad B_\alpha=\begin{pmatrix}1&0\\ \alpha/2&1\end{pmatrix},\qquad 0\le j\le4,\qquad\text{and}\qquad
-2\,\mathbf u^T\mathbf v=\alpha .
\end{equation}
\end{proposition}

\begin{proof}
We use Sylvester's identity
$\det(I+\mathsf A\mathsf B)=\det(I+\mathsf B\mathsf A)$ and the following observation. For a symmetric $4\times4$ matrix $S$, the first $j$ rows
of $I-2E_jS$ are those of $-(2S-I)$ and its other rows are those of $I$, so that
\begin{equation}\label{eq:appC-sylv}
\det(I-2E_jS)=(-1)^j\det\bigl((2S-I)_j\bigr),
\end{equation}
where $(\cdot)_j$ is the leading $j\times j$ block. The moment matrices of $\omega_{(j)}$ are $M_s(\omega_{(j)})=\mathsf W_s^T(I-2E_j)\mathsf W_s$
for $\mathsf W_2=\mathsf W$ and $\mathsf W_3=\mathsf W'$. With $S=\Pi$,
\begin{equation}\label{eq:appC-dict1}
x_j=\frac{\det(M_2-2\mathsf W^TE_j\mathsf W)}{\det M_2}=\det(I-2M_2^{-1}\mathsf W^TE_j\mathsf W)=\det(I-2E_j\Pi)=(-1)^j\det T_j ,
\end{equation}
which is the first component of $\mathbf n_j=B_\alpha\mathbf q_j$.

For the second component let $\Pi_{\mathbf v}=\Pi+\mathbf v\mathbf v^T$, so that $2\Pi_{\mathbf v}-I=T+2\mathbf v\mathbf v^T$. Let
$\mathcal W'$ be the span of $\mathsf w_0,\mathsf w_1,\mathsf w_2$, of dimension $3$ since $\det M_3=H_3>0$. It is spanned by $\mathcal W$ and
$\mathsf c\ne0$. The symmetric matrix $\Pi_{\mathbf v}$ vanishes on $\mathcal W'^\perp$, and it is positive definite on $\mathcal W'$, since
$x^T\Pi_{\mathbf v}x=\|\Pi x\|^2+(\mathbf v^Tx)^2$ and $\mathbf v^T\mathsf c=-\gamma^{-1/2}\|\mathsf c\|^2\ne0$. Put
\begin{equation}\label{eq:appC-Qplus}
\Pi^+_{\mathbf v}=\Pi+\frac\gamma{\|\mathsf c\|^2}\bigl(\mathsf a\mathsf c^T+\mathsf c\mathsf a^T\bigr)+\frac{\gamma(1+\gamma\|\mathsf a\|^2)}{\|\mathsf c\|^4}\,\mathsf c\mathsf c^T .
\end{equation}
From $\Pi\mathsf a=\mathsf a$, $\Pi\mathsf c=0$, and $\mathsf a\perp\mathsf c$ one finds $\mathbf v^T\Pi^+_{\mathbf v}=-\gamma^{1/2}\|\mathsf c\|^{-2}\mathsf c^T$ and then
$\Pi_{\mathbf v}\Pi^+_{\mathbf v}=\Pi+\|\mathsf c\|^{-2}\mathsf c\mathsf c^T$, the orthogonal projection onto $\mathcal W'$. So $\Pi^+_{\mathbf v}$, which is
symmetric with range in $\mathcal W'$, is the inverse of $\Pi_{\mathbf v}$ on $\mathcal W'$, extended by $0$. Let
$\Gamma=\mathsf W'^T\Pi^+_{\mathbf v}\mathsf W'$, which is positive definite since $\mathsf W'$ has rank $3$ and range $\mathcal W'$. For $y\in\R^3$,
$\mathsf W'\Gamma^{-1}\mathsf W'^T\Pi^+_{\mathbf v}\mathsf W'y=\mathsf W'y$. Since $\Pi^+_{\mathbf v}$ maps $\mathcal W'$ onto itself, this shows
$\mathsf W'\Gamma^{-1}\mathsf W'^T=\Pi_{\mathbf v}$ on $\mathcal W'$, and both sides vanish on $\mathcal W'^\perp$. So
\begin{equation}\label{eq:appC-Qrep}
\Pi_{\mathbf v}=\mathsf W'\Gamma^{-1}\mathsf W'^T .
\end{equation}
By Lemma~\ref{lem:appC-interp}(a), $\mathsf W'^T\mathsf a=\mathbf e_1+(\Sigma_2/H_2)\mathbf e_2$, $\mathsf W'^T\mathsf c=\|\mathsf c\|^2\mathbf e_2$, and
$\mathsf W'^T\Pi\mathsf W'=M_3-\|\mathsf c\|^2\mathbf e_2\mathbf e_2^T$, and \eqref{eq:appC-Qplus} gives
\begin{equation}\label{eq:appC-Gamma}
\Gamma=M_3+\gamma\bigl(\mathbf e_1\mathbf e_2^T+\mathbf e_2\mathbf e_1^T\bigr)+\varkappa\,\mathbf e_2\mathbf e_2^T,\qquad
\varkappa=\gamma+\frac{2\gamma\Sigma_2}{H_2}+\gamma^2\|\mathsf a\|^2-\|\mathsf c\|^2=\gamma\Bigl(1+\frac{2\Sigma_2}{H_2}-\frac\alpha2\Bigr),
\end{equation}
the last expression by \eqref{eq:appC-xichi} and \eqref{eq:alpha}. By \eqref{eq:appC-sylv} with $S=\Pi_{\mathbf v}$, \eqref{eq:appC-Qrep}, and
Sylvester's identity,
\begin{equation}\label{eq:appC-dict2}
(-1)^j\det\bigl((T+2\mathbf v\mathbf v^T)_j\bigr)=\det(I-2E_j\mathsf W'\Gamma^{-1}\mathsf W'^T)=\frac{\det(\Gamma-2\mathsf W'^TE_j\mathsf W')}{\det\Gamma} .
\end{equation}
The matrix $\Gamma-2\mathsf W'^TE_j\mathsf W'$ is $M_3(\omega_{(j)})$ with $\gamma$ added to its two entries $h_3$ and $\varkappa$ added to its
entry $h_4$. The first change turns its determinant into $H^{(\gamma)}_3(\omega_{(j)})$, by the definition of $H^{(\gamma)}_3$
(\S\ref{subsec:hankel}). The second adds $\varkappa$ times the cofactor of the entry $h_4$, which is $H_2(\omega_{(j)})=H_2x_j$, since
the first change does not meet the leading $2\times2$ block. So the numerator of \eqref{eq:appC-dict2} is $H^{(\gamma)}_3(\omega_{(j)})+\varkappa H_2x_j$.
For $j=0$ it is $\det\Gamma$, and by \eqref{eq:H3gamma}, \eqref{eq:appC-xichi}, and \eqref{eq:appC-Gamma},
$\det\Gamma=H_3-2\gamma\Sigma_2-\gamma^2H_1+\varkappa H_2=\gamma H_2$. With $\varkappa=2\gamma\Sigma_2/H_2+\gamma(1-\frac\alpha2)$ and the definition of
$\hat Y_j$ in \eqref{eq:appC-n}, \eqref{eq:appC-dict2} becomes
\begin{equation}\label{eq:appC-dict3}
(-1)^j\det\bigl((T+2\mathbf v\mathbf v^T)_j\bigr)=\frac{\hat Y_j}{\gamma H_2}+\Bigl(1-\frac\alpha2\Bigr)x_j .
\end{equation}
By the matrix determinant lemma, that is, the polynomial identity
\begin{equation}\label{eq:appC-mdl}
\det(\mathsf A+2\mathbf x\mathbf x^T)=\det\mathsf A+2\mathbf x^T\operatorname{adj}(\mathsf A)\mathbf x,
\end{equation}
for $\mathbf x\in\R^j$ and $j\times j$ matrices $\mathsf A$, the left side of \eqref{eq:appC-dict3} is
$(-1)^j\det T_j+(-1)^j2\mathbf v_{[j]}^T\operatorname{adj}(T_j)\mathbf v_{[j]}$. With \eqref{eq:appC-dict1} this gives
\begin{equation}\label{eq:appC-dict4}
(-1)^j\,2\,\mathbf v_{[j]}^T\operatorname{adj}(T_j)\mathbf v_{[j]}=\frac{\hat Y_j}{\gamma H_2}-\frac\alpha2\,x_j,
\end{equation}
which is the second component of $\mathbf n_j=B_\alpha\mathbf q_j$. Finally, by \eqref{eq:appC-xichi},
$\mathbf u^T\mathbf v=\gamma\|\mathsf a\|^2-\gamma^{-1}\|\mathsf c\|^2=(\gamma^2H_1-H_3)/(\gamma H_2)=-\frac\alpha2$.
\end{proof}

For $j=4$, $\operatorname{adj}(T)=T$ by Lemma~\ref{lem:appC-interp}(b), so $\mathbf q_4=(1,2\mathbf u^T\mathbf v)=(1,-\alpha)$, and
\eqref{eq:appC-dict} reproduces \eqref{eq:appC-nouter}. By \eqref{eq:appC-dict1}, the leading minor $(-1)^j\det T_j$ is the first component
$x_j$ of $\mathbf n_j$. It vanishes when $[\mathbf n_j]$ is the line of $(0,1)$, and the arguments below do not assume that it is
non-zero.

\subsection{First-order variations}\label{subsec:appC-var}
The convexity rests on the algebra of first-order variations of the data. Its tools include the transport of a variation of $T$ by a lower
triangular matrix, in a pencil form that needs no hypothesis on the leading minors of $T$, and two identities of Aksteiner, the first law
and the sum of squares \cite{Aksteiner2026}.

A \emph{variation} is a real-analytic one-parameter family $s\mapsto(w(s),\omega(s),\gamma(s))$, for $s$
near $0$, of data as in \eqref{eq:appC-TUV}, and $\delta$ denotes the derivative at $s=0$. A line $[\mathbf n(s)]$, with
$\mathbf n(0)\ne0$, is \emph{preserved to first order} if $\delta\mathbf n\in\R\,\mathbf n(0)$, that is, if the derivative of
$s\mapsto[\mathbf n(s)]\in\mathbb{RP}^1$ vanishes at $0$.

\begin{lemma}[transport]\label{lem:appC-transport}
Let $T$ be a real symmetric $4\times4$ matrix with $\det T_j\ne0$ for $1\le j\le4$. Then for every symmetric $4\times4$ matrix $\dot T$ there
is a unique lower triangular matrix $\mathcal K$ with $\dot T=\mathcal KT+T\mathcal K^T$. If $T$ is the matrix \eqref{eq:appC-TUV} of a
variation and $\mathcal K$ is any lower triangular matrix with $\delta T=\mathcal KT+T\mathcal K^T$, then $Q=\mathcal K+\mathcal K^T$
satisfies
\begin{equation}\label{eq:appC-ST}
QT+TQ=0,\qquad \operatorname{tr}\mathcal K=0 .
\end{equation}
\end{lemma}

\begin{proof}
The lower triangular and the symmetric $4\times4$ matrices both form spaces of dimension $10$, so it suffices to show injectivity. Since the
leading minors $D_j=\det T_j$ do not vanish, $T=\mathsf E\Delta\mathsf E^T$ with $\mathsf E$ unit lower triangular and $\Delta$ diagonal with the
entries $\Delta_j=D_j/D_{j-1}\ne0$. For lower triangular $\mathcal K$ put $\mathcal K'=\mathsf E^{-1}\mathcal K\mathsf E$, which is lower triangular. Then
$\mathcal KT+T\mathcal K^T=\mathsf E(\mathcal K'\Delta+\Delta \mathcal K'^T)\mathsf E^T$. If this vanishes, the entries $(i,l)$ of $\mathcal K'\Delta+\Delta \mathcal K'^T$ with $i>l$ give
$\mathcal K'_{il}\Delta_l=0$, and the diagonal entries give $2\mathcal K'_{ii}\Delta_i=0$. So $\mathcal K'=0$ and $\mathcal K=0$. For the second
statement, which does not use the leading minors, differentiating
$T^2=I$ (Lemma~\ref{lem:appC-interp}(b)) gives $\delta T\,T+T\,\delta T=0$, that is, $Q+TQT=0$. Multiplying by $T$ gives $QT+TQ=0$, and
$\operatorname{tr}Q=-\operatorname{tr}(TQT)=-\operatorname{tr}Q$ gives $\operatorname{tr}\mathcal K=\frac12\operatorname{tr}Q=0$.
\end{proof}

The transport of $\delta T$ alone need not exist where a leading minor of $T$ vanishes. The lines $[\mathbf n_j]$ are, however, determined
by the pair of matrices $T$ and $T+2\mathbf v\mathbf v^T$, and a generic member of the pencil spanned by them has a transport. For a real
symmetric $4\times4$ matrix $T$, a vector $\mathbf x\in\R^4$, and $\lambda\in\R$ let $T^\lambda=T+2\lambda\mathbf x\mathbf x^T$, with leading
$j\times j$ blocks $T^\lambda_j$, let $\mathbf x_{[j]}$ be the vector of the first $j$ entries of $\mathbf x$, and put
\begin{equation}\label{eq:appD-pencil}
m_j(\lambda)=\det T^\lambda_j=\det T_j+2\lambda\,\mathbf x_{[j]}^T\operatorname{adj}(T_j)\,\mathbf x_{[j]},\qquad
\hat m_j=\bigl(m_j(0),m_j(1)\bigr),\qquad 0\le j\le4,
\end{equation}
with $m_0=1$. The second expression is \eqref{eq:appC-mdl}. Thus $m_j(\lambda)=(1-\lambda)m_j(0)+\lambda m_j(1)$.

\begin{lemma}[pencil transport]\label{lem:appD-pencil}
Let $T$ and $\mathbf x$ be as above with $\det(\hat m_{j-1},\hat m_j)\ne0$ for $1\le j\le4$. Let $s\mapsto(T(s),\mathbf x(s))$ be a
differentiable family of real symmetric matrices and vectors with $(T(0),\mathbf x(0))=(T,\mathbf x)$, along which the lines $[\hat m_j]$
are preserved to first order for $1\le j\le4$. Then there is exactly one lower triangular $4\times4$ matrix $\mathcal K$ with
\begin{equation}\label{eq:appD-pencilK}
\delta T=\mathcal KT+T\mathcal K^T,\qquad \delta\mathbf x=\mathcal K\mathbf x,
\end{equation}
and for it $\delta m_j(\lambda)=2\operatorname{tr}(\mathcal K_j)\,m_j(\lambda)$ for $0\le j\le4$ and all $\lambda\in\R$, where $\mathcal K_j$ is
the leading $j\times j$ block of $\mathcal K$.
\end{lemma}

\begin{proof}
Since $\det(\hat m_{j-1},\hat m_j)\ne0$, no $\hat m_j$ vanishes, so no $m_j$ vanishes identically, and we can fix $\lambda\in\R$ with
$m_j(\lambda)\ne0$ for $1\le j\le4$. The leading minors of $T^\lambda$ are the numbers $m_j(\lambda)$, so by Lemma~\ref{lem:appC-transport}
there is exactly one lower triangular $\mathcal K$ with $\delta T^\lambda=\mathcal KT^\lambda+T^\lambda\mathcal K^T$. Put
$\mathbf z=\delta\mathbf x-\mathcal K\mathbf x$, fix $\lambda'\ne\lambda$, and let $\varphi_j=\mathbf x_{[j]}^T(T^\lambda_j)^{-1}\mathbf x_{[j]}$ for
$1\le j\le4$ and $\varphi_0=0$. By \eqref{eq:appC-mdl} with $T^\lambda_j$ in place of $T_j$,
\begin{equation}\label{eq:appD-pencilphi}
m_j(\lambda')=m_j(\lambda)\bigl(1+2(\lambda'-\lambda)\varphi_j\bigr),\qquad 0\le j\le4 .
\end{equation}
The pair $(m_j(\lambda),m_j(\lambda'))$ is the image of $\hat m_j$ under a fixed linear map of determinant $\lambda'-\lambda$, and
$m_j(\lambda)\ne0$. So $[\hat m_j]$ is preserved to first order if and only if the quotient $m_j(\lambda')/m_j(\lambda)$ is constant to
first order, that is, if and only if $\delta\varphi_j=0$. Since $\mathcal K$ is lower triangular,
$\delta T^\lambda_j=\mathcal K_jT^\lambda_j+T^\lambda_j\mathcal K_j^T$ and $(\mathcal K\mathbf x)_{[j]}=\mathcal K_j\mathbf x_{[j]}$, so
\begin{equation}\label{eq:appD-pencildphi}
\delta\varphi_j=2\bigl(\mathcal K_j\mathbf x_{[j]}+\mathbf z_{[j]}\bigr)^T(T^\lambda_j)^{-1}\mathbf x_{[j]}
-\mathbf x_{[j]}^T\bigl((T^\lambda_j)^{-1}\mathcal K_j+\mathcal K_j^T(T^\lambda_j)^{-1}\bigr)\mathbf x_{[j]}
=2\,\mathbf z_{[j]}^T(T^\lambda_j)^{-1}\mathbf x_{[j]} .
\end{equation}
Hence $\mathbf z_{[j]}^T(T^\lambda_j)^{-1}\mathbf x_{[j]}=0$ for $1\le j\le4$. This is a triangular system of equations for
the entries of $\mathbf z$, in which the $j$-th equation has the coefficient $g_j=((T^\lambda_j)^{-1}\mathbf x_{[j]})_j$ at the $j$-th entry. Write
$T^\lambda_j$ with the leading block $T^\lambda_{j-1}$, the last column $(\mathbf b,\mu)$, and $\mathbf x_{[j]}=(\mathbf x_{[j-1]},\xi)$ (for $j=1$ the leading block is empty). Let
$\varsigma=\mu-\mathbf b^T(T^\lambda_{j-1})^{-1}\mathbf b=m_j(\lambda)/m_{j-1}(\lambda)$ be the Schur complement. The formula for the inverse of a
block matrix gives $g_j=\varsigma^{-1}\bigl(\xi-\mathbf b^T(T^\lambda_{j-1})^{-1}\mathbf x_{[j-1]}\bigr)$ and $\varphi_j-\varphi_{j-1}=\varsigma g_j^2$.
On the other hand, by \eqref{eq:appD-pencilphi} for $j$ and $j-1$,
\begin{equation}\label{eq:appD-pencilg}
\varphi_j-\varphi_{j-1}=\frac{m_{j-1}(\lambda)m_j(\lambda')-m_j(\lambda)m_{j-1}(\lambda')}{2(\lambda'-\lambda)\,m_{j-1}(\lambda)m_j(\lambda)}
=\frac{\det(\hat m_{j-1},\hat m_j)}{2\,m_{j-1}(\lambda)m_j(\lambda)} ,
\end{equation}
so $g_j^2=\det(\hat m_{j-1},\hat m_j)/(2m_j(\lambda)^2)\ne0$. Hence $\mathbf z=0$, that is, $\delta\mathbf x=\mathcal K\mathbf x$, and
$\delta T=\delta T^\lambda-2\lambda(\delta\mathbf x\,\mathbf x^T+\mathbf x\,\delta\mathbf x^T)=\mathcal KT+T\mathcal K^T$. Conversely, a lower triangular
matrix with \eqref{eq:appD-pencilK} satisfies $\delta T^\lambda=\mathcal KT^\lambda+T^\lambda\mathcal K^T$, so it is $\mathcal K$. Finally,
\eqref{eq:appD-pencilK} gives $\delta T^{\lambda''}=\mathcal KT^{\lambda''}+T^{\lambda''}\mathcal K^T$ for every $\lambda''\in\R$. Since $\mathcal K$ is lower
triangular, $\delta T^{\lambda''}_j=\mathcal K_jT^{\lambda''}_j+T^{\lambda''}_j\mathcal K_j^T$. As
$\operatorname{adj}(T^{\lambda''}_j)T^{\lambda''}_j=T^{\lambda''}_j\operatorname{adj}(T^{\lambda''}_j)=\det(T^{\lambda''}_j)I$, this gives
$\delta\det T^{\lambda''}_j=\operatorname{tr}\bigl(\operatorname{adj}(T^{\lambda''}_j)\,\delta T^{\lambda''}_j\bigr)
=2\operatorname{tr}(\mathcal K_j)\det T^{\lambda''}_j$. This last step uses only that $\mathcal K$ is lower triangular and satisfies
\eqref{eq:appD-pencilK}.
\end{proof}

\begin{lemma}[direction-preserving variations]\label{lem:appC-fibre}\label{lem:appD-fibre}
Let $(w(s),\omega(s),\gamma(s))$ be a variation along which $\alpha(\omega(s),\gamma(s))$ does not depend on $s$, and let $T$, $\mathbf v$,
and $\mathbf n_j$ be as in \eqref{eq:appC-TUV} and \eqref{eq:appC-n}. Then the lines $[\mathbf n_1],[\mathbf n_2],[\mathbf n_3]$ are
preserved to first order if and only if there is a lower triangular matrix $\mathcal K$ with $\delta T=\mathcal KT+T\mathcal K^T$ and
$\delta\mathbf v=\mathcal K\mathbf v$. In that case $\mathcal K$ is unique, \eqref{eq:appC-ST} holds with $Q=\mathcal K+\mathcal K^T$, and,
with $\mathcal K_j$ the leading $j\times j$ block of $\mathcal K$,
\begin{equation}\label{eq:appC-dn}
\delta\mathbf n_j=2\operatorname{tr}(\mathcal K_j)\,\mathbf n_j,\qquad 1\le j\le3 .
\end{equation}
\end{lemma}

\begin{proof}
Let $m_j$ and $\hat m_j$ be as in \eqref{eq:appD-pencil} with $\mathbf x=\mathbf v$. By \eqref{eq:appC-rods}, $\hat m_j=(-1)^j(q^{(1)}_j,q^{(1)}_j+q^{(2)}_j)$,
where $q^{(1)}_j,q^{(2)}_j$ are the components of $\mathbf q_j$, and the map $(p_1,p_2)\mapsto(p_1,p_1+p_2)$ has determinant $1$. So
$\det(\hat m_{j-1},\hat m_j)=-\det(\mathbf q_{j-1},\mathbf q_j)=-\det(\mathbf n_{j-1},\mathbf n_j)>0$ for $1\le j\le4$, by \eqref{eq:appC-dict},
$\det B_\alpha=1$, and \S\ref{subsec:appC-rods}. Since $\alpha$ is constant, so is $B_\alpha$, and by \eqref{eq:appC-dict} the line $[\hat m_j]$ is
preserved to first order if and only if $[\mathbf n_j]$ is. For $j=4$ this holds since $\mathbf n_4=(1,-\frac\alpha2)$ is constant
\eqref{eq:appC-nouter}. If the lines $[\mathbf n_1],[\mathbf n_2],[\mathbf n_3]$ are preserved, Lemma~\ref{lem:appD-pencil} gives exactly
one $\mathcal K$ as in the statement. Conversely, if $\mathcal K$ is lower triangular with $\delta T=\mathcal KT+T\mathcal K^T$ and
$\delta\mathbf v=\mathcal K\mathbf v$, then $\delta T^\lambda=\mathcal KT^\lambda+T^\lambda\mathcal K^T$ for every $\lambda\in\R$. As in
the last step of the proof of Lemma~\ref{lem:appD-pencil}, this gives $\delta m_j(\lambda)=2\operatorname{tr}(\mathcal K_j)m_j(\lambda)$ for all
$\lambda$. For $\lambda=0,1$ this gives $\delta\mathbf q_j=2\operatorname{tr}(\mathcal K_j)\mathbf q_j$, so the lines $[\mathbf q_j]$, and with
them the lines $[\mathbf n_j]$, are preserved, and \eqref{eq:appC-dn} follows by \eqref{eq:appC-dict}. Finally \eqref{eq:appC-ST} holds by
Lemma~\ref{lem:appC-transport}.
\end{proof}

\begin{lemma}[the first law and the sum of squares]\label{lem:appC-identities}
Let $(w(s),\omega(s),\gamma(s))$ be a variation, let $\mathcal K$ be a lower triangular matrix with $\delta T=\mathcal KT+T\mathcal K^T$
and $\delta\mathbf v=\mathcal K\mathbf v$, for instance that of Lemma~\textup{\ref{lem:appC-fibre}}, and put $Q=\mathcal K+\mathcal K^T$
and $L_j=w_{j+1}-w_j$ for $j=1,2,3$. Then
\begin{equation}\label{eq:appC-firstlaw}
\sum_{j=1}^3L_j\cdot2\operatorname{tr}(\mathcal K_j)=\delta\|\mathbf v\|^2,
\end{equation}
\begin{equation}\label{eq:appC-sos}
\sum_{j=1}^3\delta L_j\cdot2\operatorname{tr}(\mathcal K_j)=\sum_{i>l}(w_i-w_l)\mathcal K_{il}^2+\|Q\mathbf v\|^2 .
\end{equation}
\end{lemma}

\begin{proof}
Summation by parts with $\operatorname{tr}\mathcal K=0$ \eqref{eq:appC-ST} gives, for any numbers $c_1,\dots,c_4$,
\begin{equation}\label{eq:appC-parts}
\sum_{j=1}^3(c_{j+1}-c_j)\cdot2\operatorname{tr}(\mathcal K_j)=2\sum_{i=1}^4\mathcal K_{ii}(c_4-c_i)=-2\sum_{i=1}^4c_i\mathcal K_{ii}=-\operatorname{tr}\bigl(\mathrm{diag}(c_i)\,Q\bigr).
\end{equation}
Let $\mathsf Z_{\mathbf v}=\mathsf Z+\mathbf v\mathbf v^T$, which commutes with $T$ by \eqref{eq:appC-interp}. By \eqref{eq:appC-ST}, $Q=-TQT$, so
$\operatorname{tr}(\mathsf Z_{\mathbf v}Q)=-\operatorname{tr}(\mathsf Z_{\mathbf v}TQT)=-\operatorname{tr}(T\mathsf Z_{\mathbf v}QT)=-\operatorname{tr}(\mathsf Z_{\mathbf v}Q)$, and $\operatorname{tr}(\mathsf Z_{\mathbf v}Q)=0$. With
$c_i=w_i$, \eqref{eq:appC-parts} gives $\sum_jL_j\cdot2\operatorname{tr}(\mathcal K_j)=-\operatorname{tr}(\mathsf ZQ)=\mathbf v^TQ\mathbf v=2\mathbf v^T\mathcal K\mathbf v=\delta\|\mathbf v\|^2$,
which is \eqref{eq:appC-firstlaw}. With $c_i=\delta w_i$, \eqref{eq:appC-parts} gives $\sum_j\delta L_j\cdot2\operatorname{tr}(\mathcal K_j)=-\operatorname{tr}(\delta \mathsf Z\,Q)$.
Differentiating $[\mathsf Z_{\mathbf v},T]=0$, with $\delta\mathsf Z_{\mathbf v}=\delta \mathsf Z+\mathcal K\mathbf v\mathbf v^T+\mathbf v\mathbf v^T\mathcal K^T$ and
$[\mathsf Z_{\mathbf v},\mathcal KT+T\mathcal K^T]=[\mathsf Z_{\mathbf v},\mathcal K]T+T[\mathsf Z_{\mathbf v},\mathcal K^T]$, and expanding, one finds that
\begin{equation}\label{eq:appC-F}
\mathcal F=\delta \mathsf Z+[\mathsf Z,\mathcal K]+\mathbf v\mathbf v^TQ\qquad\text{satisfies}\qquad \mathcal FT=T\mathcal F^T .
\end{equation}
Indeed, $\mathcal FT-T\mathcal F^T=[\delta \mathsf Z,T]+[\mathsf Z,\mathcal K]T+T[\mathsf Z,\mathcal K^T]+\mathbf v\mathbf v^TQT-TQ\mathbf v\mathbf v^T$. The derivative of $[\mathsf Z_{\mathbf v},T]=0$ is
$[\delta \mathsf Z,T]+[\mathsf Z,\mathcal K]T+T[\mathsf Z,\mathcal K^T]+[\mathbf v\mathbf v^T,\mathcal K]T+T[\mathbf v\mathbf v^T,\mathcal K^T]+[\mathcal K\mathbf v\mathbf v^T+\mathbf v\mathbf v^T\mathcal K^T,T]=0$, and in it the last three terms
add up to $\mathbf v\mathbf v^TQT-TQ\mathbf v\mathbf v^T$. So $\mathcal FT-T\mathcal F^T$ is this derivative and vanishes. By \eqref{eq:appC-F} and \eqref{eq:appC-ST},
$\operatorname{tr}(\mathcal FQ)=\operatorname{tr}(T\mathcal F^TTQ)=\operatorname{tr}(\mathcal F^TTQT)=-\operatorname{tr}(\mathcal F^TQ)=-\operatorname{tr}(\mathcal FQ)$, so
$\operatorname{tr}(\mathcal FQ)=0$. Here $\operatorname{tr}(\mathbf v\mathbf v^TQ^2)=\|Q\mathbf v\|^2$ and, since $\operatorname{tr}(\mathsf Z\mathcal K^2)=\operatorname{tr}(\mathcal K\mathsf Z\mathcal K)$,
\begin{equation}\label{eq:appC-ZK}
\operatorname{tr}([\mathsf Z,\mathcal K]Q)=\operatorname{tr}(\mathsf Z\mathcal K\mathcal K^T)-\operatorname{tr}(\mathcal K\mathsf Z\mathcal K^T)=\sum_{i,l}(w_i-w_l)\mathcal K_{il}^2=\sum_{i>l}(w_i-w_l)\mathcal K_{il}^2,
\end{equation}
as $\mathcal K$ is lower triangular. So $-\operatorname{tr}(\delta \mathsf Z\,Q)=\sum_{i>l}(w_i-w_l)\mathcal K_{il}^2+\|Q\mathbf v\|^2$, which is \eqref{eq:appC-sos}.
\end{proof}

\begin{lemma}[irreducibility]\label{lem:appC-irred}
Let $\omega=\sum_ka_k\delta_{w_k}$ be a positive measure on points $w_1<\dots<w_4$ and $\Pi$ as in \eqref{eq:appC-TUV}. Then the
following hold.
\begin{enumerate}[leftmargin=*,label=\textup{(\alph*)}]
\item A diagonal $4\times4$ matrix that commutes with $\Pi$ is a multiple of $I$.
\item If a variation with $\delta w=0$ has $\delta\Pi=0$, then $\delta a_k/a_k$ does not depend on $k$.
\end{enumerate}
\end{lemma}

\begin{proof}
To prove (a), suppose that a diagonal matrix with at least two distinct entries commutes with $\Pi$. Then so does the diagonal projection onto the coordinates
in a set $J'$ on which it takes one of its values. Let $J''$ be the complement, so that $J',J''\ne\emptyset$. Then $\mathcal W$ is the direct
sum of its intersections with $\R^{J'}$ and $\R^{J''}$. The elements of $\mathcal W$ are the vectors $(a_k^{1/2}p(w_k))_k$ with polynomials
$p$ of degree at most $1$, and such a vector vanishes on $J''$ only if $p$ has $|J''|$ zeros. So $\mathcal W\cap\R^{J'}$ has dimension
$\max(0,2-|J''|)$, and similarly for $J''$. Since $|J'|+|J''|=4$, these dimensions add up to at most $1<2$, a contradiction.

To prove (b), put $\upsilon_k=\delta a_k/(2a_k)$, so that $\delta\mathsf w_b=\mathrm{diag}(\upsilon_k)\mathsf w_b$. Differentiating $\Pi\mathsf w_b=\mathsf w_b$ for
$b=0,1$ gives $(I-\Pi)\delta\mathsf w_b=\delta\Pi\,\mathsf w_b=0$, so $(\upsilon_ka_k^{1/2})_k$ and $(\upsilon_ka_k^{1/2}w_k)_k$ lie in $\mathcal W$. Hence there
are polynomials $p,p'$ of degree at most $1$ with $\upsilon_k=p(w_k)$ and $\upsilon_kw_k=p'(w_k)$ for all $k$. The polynomial $p(x)x-p'(x)$ has
degree at most $2$ and four zeros, so it vanishes, and $p$ is constant.
\end{proof}

\subsection{Convexity}\label{subsec:appC-convex}
The cone errors are the gradient of a strictly convex function of the rod lengths (Proposition~\ref{prop:appC-convex}). The first step of
the proof is the local invertibility of the map from the weights to the lines at fixed turning points.

Let $\alpha\ge0$. For a positive measure $\omega$ on points $w_1<\dots<w_4$ let $\gamma_\alpha(\omega)$ be the unique $\gamma>0$ with
$\alpha(\omega,\gamma)=\alpha$, the positive root of $H_1\gamma^2+\frac\alpha2H_2\gamma-H_3=0$ by \eqref{eq:alpha}. It is a real-analytic
function of the points and the weights, and $\gamma_\alpha(c\omega)=c\gamma_\alpha(\omega)$. For $\alpha=0$,
$\gamma_0(\omega)=\bigl(H_3(\omega)/H_1(\omega)\bigr)^{1/2}$. By \eqref{eq:appC-nscale}, the data
$(c\omega,c\gamma)$ have the same vectors $\mathbf n_j$ as $(\omega,\gamma)$, so we may normalise $a_1=1$. Put
$\mathcal A=\{a\in\R^4_{>0}:a_1=1\}$, write $\omega_a=\sum_ka_k\delta_{w_k}$, and for $a\in\mathcal A$ let
\begin{equation}\label{eq:appC-Psi}
\Psi^\alpha_w(a)=\bigl([\mathbf n_1],[\mathbf n_2],[\mathbf n_3]\bigr)\in(\mathbb{RP}^1)^3,\qquad \mathbf n_j=\mathbf n_j\bigl(\omega_a,\gamma_\alpha(\omega_a)\bigr).
\end{equation}
The vectors $\mathbf n_j$ do not vanish, since consecutive ones have negative determinant. So $\Psi^\alpha_w$ is real-analytic, and it
depends real-analytically on $w$ as well.

\begin{lemma}\label{lem:appC-local}\label{lem:appD-local}\label{lem:appD-localall}
Let $\alpha\ge0$ and $w_1<\dots<w_4$. Then the differential of the map $\Psi^\alpha_w$ of \eqref{eq:appC-Psi} is invertible at every
$a\in\mathcal A$.
\end{lemma}

\begin{proof}
Both $\mathcal A$ and $(\mathbb{RP}^1)^3$ have dimension $3$, so it suffices to show that the kernel is zero. Let $\dot a$ be in the
kernel at $a$, so that $\dot a_1=0$, and consider the variation $s\mapsto(w,\omega_{a+s\dot a},\gamma_\alpha(\omega_{a+s\dot a}))$. It has
$\delta w=0$ and constant $\alpha$, and it preserves the lines $[\mathbf n_1],[\mathbf n_2],[\mathbf n_3]$ to first order. Let $\mathcal K$ be
as in Lemma~\ref{lem:appC-fibre}. By \eqref{eq:appC-sos} with $\delta L=0$, $\mathcal K_{il}=0$ for $i>l$ and $Q\mathbf v=0$. So $\mathcal K$ is
diagonal, $Q=2\mathcal K$, and \eqref{eq:appC-ST} gives $\mathcal KT+T\mathcal K=0$, that is, $\delta T=0$. Hence $\delta\Pi=\frac12\delta T=0$,
and Lemma~\ref{lem:appC-irred}(b) gives that $\delta a_k/a_k$ does not depend on $k$. As $\delta a_1=\dot a_1=0$, it vanishes. So
$\dot a=0$.
\end{proof}

Lemma~\ref{lem:appC-local} does not assume that the leading minors of $T$ are non-zero. Where one of them vanishes, the transport of
$\delta T$ alone (Lemma~\ref{lem:appC-transport}) need not exist, and the pencil transport (Lemma~\ref{lem:appD-pencil}) takes its place.

\begin{proposition}[convexity]\label{prop:appC-convex}\label{prop:appD-convex}
Let $\alpha\ge0$, let $\ell_0,\dots,\ell_4$ be as in \eqref{eq:appC-adm}, and let $\mathfrak b$ be as in \eqref{eq:appC-b}. Then
$\mathfrak b\colon\R^3_{>0}\to\R^3$ is real-analytic, its Jacobian matrix is symmetric and positive definite at every point, and there is a
strictly convex real-analytic function $\mathfrak F_\ell$ on $\R^3_{>0}$ with $\mathfrak b=\nabla\mathfrak F_\ell$. In particular at most one
$L\in\R^3_{>0}$ is balanced.
\end{proposition}

\begin{proof}
Let $L^0\in\R^3_{>0}$. By Proposition~\ref{prop:appC-ident}(i), the second identity of \eqref{eq:appC-nscale}, and
$\gamma_\alpha(c\omega)=c\gamma_\alpha(\omega)$, there is $a^0\in\mathcal A$ such that $(\omega_{a^0},\gamma_\alpha(\omega_{a^0}))$, on the
points $w(L^0)$, is direction matched for $L^0$. The map $(L,a)\mapsto\Psi^\alpha_{w(L)}(a)$ is real-analytic, it takes the value
$([\ell_1],[\ell_2],[\ell_3])$ at $(L^0,a^0)$, and its differential in $a$ is invertible there by Lemma~\ref{lem:appC-local}. By the
real-analytic implicit function theorem there are a neighbourhood $U$ of $L^0$ and a real-analytic map $L\mapsto a(L)\in\mathcal A$ on $U$
with $a(L^0)=a^0$ and $\Psi^\alpha_{w(L)}(a(L))=([\ell_1],[\ell_2],[\ell_3])$. Put $\omega_L=\omega_{a(L)}$ on the points $w(L)$ and
$\gamma_L=\gamma_\alpha(\omega_L)$. These data are direction matched for $L$, so by Proposition~\ref{prop:appC-ident}(ii) and
\eqref{eq:appC-nb}, $\mathfrak b_j(L)=\log\bigl(\mathbf n_j\cdot\ell_j/|\ell_j|^2\bigr)$ with $\mathbf n_j=\mathbf n_j(\omega_L,\gamma_L)$.
This is real-analytic on $U$. As $L^0$ is arbitrary, $\mathfrak b$ is real-analytic.

Let $L\in U$, $\dot L\in\R^3$, and consider the variation $s\mapsto(w(L+s\dot L),\omega_{L+s\dot L},\gamma_{L+s\dot L})$. It has
constant $\alpha$ and it preserves the lines $[\mathbf n_j]=[\ell_j]$. Let $\mathcal K$ be as in Lemma~\ref{lem:appC-fibre}. By
\eqref{eq:appC-dn}, $\delta\mathbf n_j=2\operatorname{tr}(\mathcal K_j)\mathbf n_j$. Since $\mathbf n_j=e^{\mathfrak b_j}\ell_j$ with $\ell_j$
fixed, $\delta\mathbf n_j=\delta\mathfrak b_j\,\mathbf n_j$, so
\begin{equation}\label{eq:appC-db}
\delta\mathfrak b_j=\bigl(D\mathfrak b(L)\dot L\bigr)_j=2\operatorname{tr}(\mathcal K_j),\qquad 1\le j\le3 .
\end{equation}
Let $\tilde{\mathfrak F}(L)=\sum_jL_j\mathfrak b_j(L)-\|\mathbf v\|^2$ on $U$, where $L_j=w_{j+1}(L)-w_j(L)$ and
$\mathbf v$ is formed with $(\omega_L,\gamma_L)$. By \eqref{eq:appC-firstlaw}, $d\tilde{\mathfrak F}=\sum_j\mathfrak b_j\,dL_j+\sum_jL_j\,d\mathfrak b_j-d\|\mathbf v\|^2=\sum_j\mathfrak b_j\,dL_j$.
So the Jacobian matrix of $\mathfrak b$ is the Hessian of $\tilde{\mathfrak F}$ on $U$, and it is symmetric. Hence the $1$-form
$\sum_j\mathfrak b_j(L)\,dL_j$ is closed on $\R^3_{>0}$, and since $\R^3_{>0}$ is convex it is the differential of the real-analytic function
\begin{equation}\label{eq:appC-potential}
\mathfrak F_\ell(L)=\int_0^1\mathfrak b\bigl(L^\circ+s(L-L^\circ)\bigr)\cdot(L-L^\circ)\,ds ,
\end{equation}
for any fixed $L^\circ\in\R^3_{>0}$. By \eqref{eq:appC-db} and \eqref{eq:appC-sos},
\begin{equation}\label{eq:appC-hessian}
\dot L^T\,D\mathfrak b(L)\,\dot L=\sum_j\delta L_j\,\delta\mathfrak b_j=\sum_{i>l}(w_i-w_l)\mathcal K_{il}^2+\|Q\mathbf v\|^2\ge0 ,
\end{equation}
since $w_i>w_l$ for $i>l$. If the left side vanishes, then, as in the proof of Lemma~\ref{lem:appC-local}, $\mathcal K$ is diagonal with
$Q\mathbf v=0$, $Q=2\mathcal K$, and $\delta T=0$, so $\delta\mathbf v=\mathcal K\mathbf v=\frac12Q\mathbf v=0$. Differentiating
$[\mathsf Z+\mathbf v\mathbf v^T,T]=0$ then gives $[\delta \mathsf Z,T]=0$, so $\delta \mathsf Z$ commutes with $\Pi$, and by Lemma~\ref{lem:appC-irred}(a) all
$\delta w_k$ are equal. As $w_1(L)=0$ for all $L$, $\delta w=0$ and $\dot L=0$. So the Hessian $D\mathfrak b$ of $\mathfrak F_\ell$ is
positive definite, and $\mathfrak F_\ell$ is strictly convex. If $\mathfrak b(L)=\mathfrak b(L')=0$ with $L\ne L'$, the function
$s\mapsto\mathfrak F_\ell(L+s(L'-L))$ on $[0,1]$ would have vanishing derivative at both ends and a positive second derivative, which is
impossible.
\end{proof}

\subsection{The ALF rod structures of the chain $(-1,1,-1)$}\label{subsec:appC-ALF}
Together with a scaling property of the cone errors (Lemma~\ref{lem:appC-scale}), Proposition~\ref{prop:appC-convex} shows that the
balanced rod lengths of the rod structures of type ALF$^-$ of the chain $\kappa=(-1,1,-1)$ lie on one ray for all values of $\alpha$
(Corollary~\ref{cor:appC-ray}).

Let $v_0,\dots,v_4\in\Z^2$ be the vectors \eqref{eq:ALF-rods}, $v_0=(1,0)$, $v_1=(0,1)$, $v_2=(-1,-1)$,
$v_3=(-1,-2)$, $v_4=(2,3)$. They satisfy $v_{j+1}+v_{j-1}=\kappa_jv_j$ for $j=1,2,3$ and
$\det(v_{j-1},v_j)=1$ for $1\le j\le4$, and $\det(v_0,v_4)=3$, as determinants of the coordinates in the negatively oriented basis
$(v_0,v_1)$ (Section~\ref{sec:explicit}). Let $\alpha>0$, let $F_\alpha$ be the linear map of $\R^2$ with $F_\alpha(v_0)=(1,\frac\alpha2)$ and
$F_\alpha(v_4)=(1,-\frac\alpha2)$, and put $\ell_j=F_\alpha(v_j)$. Explicitly,
\begin{equation}\label{eq:appC-ell}
\ell_0=\bigl(1,\tfrac\alpha2\bigr),\quad \ell_1=\bigl(-\tfrac13,-\tfrac\alpha2\bigr),\quad \ell_2=\bigl(-\tfrac23,0\bigr),\quad
\ell_3=\bigl(-\tfrac13,\tfrac\alpha2\bigr),\quad \ell_4=\bigl(1,-\tfrac\alpha2\bigr),
\end{equation}
so that $\ell_{j+1}+\ell_{j-1}=\kappa_j\ell_j$ and
\begin{equation}\label{eq:appC-detell}
\det(\ell_{j-1},\ell_j)=\det F_\alpha=-\tfrac\alpha3,\qquad 1\le j\le4 .
\end{equation}
(The lines $[\ell_j]$ are the rod lines of \cite[\S9]{Gerochcf}.) The first components $1,-\frac13,-\frac23,-\frac13,1$ are the numbers of Lemma~\ref{lem:seq}(b)
with $\varepsilon=-1$. By \eqref{eq:appC-ell} and \eqref{eq:appC-detell} the vectors $\ell_j$ satisfy \eqref{eq:appC-adm}. We write
$R_{L,\alpha}$, $\Phi_{L,\alpha}$, $\nu_{L,\alpha}$, $\tilde v_j(L,\alpha)$, and $\mathfrak b(L,\alpha)$ for the objects $R_L$, $\Phi_L$, $\nu_L$,
$\tilde v_j(L)$, and $\mathfrak b(L)$ of \S\ref{subsec:appC-rods} formed with these vectors. We call the pair $(L,\alpha)$ \emph{balanced} if
$L$ is, and call data \emph{direction matched} for $(L,\alpha)$ if they are direction matched for $L$. The rod structure $R_{L,\alpha}$ has
the asymptotic type ALF$^-$ and the model $\Phi^\sharp_{1,\alpha}$, with the model augmentation $\frac12\log(1+\alpha/(2R))$. For the lattice
$F_\alpha(\Z^2)$, in which every $\ell_j$ is primitive and every pair $(\ell_{j-1},\ell_j)$ is a basis, the cone angle along the rod $I_j$
is $2\pi e^{-\mathfrak b_j(L,\alpha)}$ \cite[Def.~4.19]{LS2025}. By
Proposition~\ref{prop:appC-convex}, for every $\alpha>0$ the map $\mathfrak b(\cdot,\alpha)$ is the gradient of a strictly convex
real-analytic function $\mathfrak F_\alpha$ on $\R^3_{>0}$, the function $\mathfrak F_\ell$ of \eqref{eq:appC-potential} for these vectors. Hence there is at most one $L\in\R^3_{>0}$ such that $(L,\alpha)$ is balanced.

\begin{lemma}\label{lem:appC-scale}
Let $\alpha>0$, $L\in\R^3_{>0}$, and $c>0$, and let $\mathfrak b(\cdot,\alpha)$ be as above. Then $\mathfrak b(cL,c\alpha)=\mathfrak b(L,\alpha)$.
\end{lemma}

\begin{proof}
Let $(\omega,\gamma)$, $\omega=\sum_ka_k\delta_{w_k(L)}$, be direction matched for $(L,\alpha)$ (Proposition~\ref{prop:appC-ident}(i)),
and put $\omega'=\sum_ka_k\delta_{cw_k(L)}$ and $\gamma'=c^3\gamma$. Then $cw(L)=w(cL)$, $\alpha(\omega',\gamma')=c\alpha$, and by
\eqref{eq:appC-nscale} and \eqref{eq:appC-nb}, $\mathbf n_j(\omega',\gamma')=e^{\mathfrak b_j(L,\alpha)}\mathrm{diag}(1,c)\,\ell_j$. Since
$\mathrm{diag}(1,c)F_\alpha$ and $F_{c\alpha}$ agree on $v_0$ and $v_4$, $\mathrm{diag}(1,c)\ell_j$ is the vector
$\ell_j$ of \eqref{eq:appC-ell} for $c\alpha$. So $(\omega',\gamma')$ is direction matched for $(cL,c\alpha)$, and
Proposition~\ref{prop:appC-ident}(ii) gives $e^{\mathfrak b_j(cL,c\alpha)}=e^{\mathfrak b_j(L,\alpha)}$.
\end{proof}

\begin{corollary}\label{cor:appC-ray}
Let $\alpha,\alpha'>0$ and $L,L'\in\R^3_{>0}$ be such that the pairs $(L,\alpha)$ and $(L',\alpha')$ are balanced. Then
$(L',\alpha')=c\,(L,\alpha)$ with $c=\alpha'/\alpha$.
\end{corollary}

\begin{proof}
By Lemma~\ref{lem:appC-scale}, $\mathfrak b(cL,\alpha')=\mathfrak b(L,\alpha)=0$, and by Proposition~\ref{prop:appC-convex} at $\alpha'$, $L'=cL$.
\end{proof}

\subsection{Positive solutions are balanced}\label{subsec:appC-balanced}

\begin{proposition}\label{prop:appC-balanced}
Let $\mathbf d=(L_2,L_3,a_2,a_3,a_4,\gamma)$ be a solution of \eqref{eq:ALF-lin} with positive coordinates, with the normalised data
\eqref{eq:ALF-unknowns}, $\omega=\sum_ka_k\delta_{w_k}$ and $a_1=1$, and put $L=(1,L_2,L_3)$. Then $\alpha=\alpha(\omega,\gamma)>0$, the data
$(\omega,\gamma)$ are direction matched for $(L,\alpha)$, $\mathbf n_j(\omega,\gamma)=\ell_j$ for $0\le j\le4$, and the pair $(L,\alpha)$ is balanced.
\end{proposition}

\begin{proof}
By \eqref{eq:ALF-lin}, $V_4=2V_0+3V_1$, so $\det(V_0,V_4)=3\det(V_0,V_1)<0$ \cite[Prop.~5.8(i)]{Gerochcf}. On the other hand, by
\eqref{eq:H3gamma}, $V_0=(H_2,H_3-2\gamma\Sigma_2-\gamma^2H_1)$ and, since $\omega_{(4)}=-\omega$, $V_4=(H_2,-H_3-2\gamma\Sigma_2+\gamma^2H_1)$, so
that $\det(V_0,V_4)=-2H_2(H_3-\gamma^2H_1)$. Since $H_2>0$ and $\gamma>0$, $H_3-\gamma^2H_1>0$, and $\alpha>0$ by \eqref{eq:alpha}.

By \eqref{eq:appC-A}, $\mathbf n_j=A_{\omega,\gamma}V_j$ with one invertible matrix $A_{\omega,\gamma}$, so by \eqref{eq:ALF-lin} the vectors $\mathbf n_j$
satisfy $\mathbf n_{j+1}+\mathbf n_{j-1}=\kappa_j\mathbf n_j$ for $j=1,2,3$, as do the $v_j$. Let $G$ be the linear map of $\R^2$ with
$G(v_0)=\mathbf n_0$ and $G(v_1)=\mathbf n_1$. By the two recursions, $G(v_j)=\mathbf n_j$ for $0\le j\le4$. By \eqref{eq:appC-nouter},
$G(v_0)=(1,\frac\alpha2)=F_\alpha(v_0)$ and $G(v_4)=(1,-\frac\alpha2)=F_\alpha(v_4)$. Since $v_0,v_4$ form a basis of
$\R^2$, $G=F_\alpha$ and $\mathbf n_j=\ell_j$ for all $j$. The turning points $(0,1,1+L_2,1+L_2+L_3)$ of $\mathbf d$ are the points $w(L)$ of
\eqref{eq:appC-wL}. So $(\omega,\gamma)$ is direction matched for $(L,\alpha)$, and \eqref{eq:appC-nb} gives $e^{\mathfrak b_j(L,\alpha)}\ell_j=\ell_j$, that is,
$\mathfrak b(L,\alpha)=0$.
\end{proof}
\bibliographystyle{amsplain}
\bibliography{n4}
\end{document}